\documentclass[english]{article}
\usepackage[T1]{fontenc}
\usepackage[latin9]{inputenc}
\usepackage{geometry}
\usepackage{color}
\usepackage{mathtools}
\usepackage{amsmath}
\usepackage{amsthm}
\usepackage{amssymb}
\usepackage{authblk}
\usepackage{hyperref}
\usepackage{setspace}
\usepackage{url}

\makeatletter
\theoremstyle{definition}
\newtheorem{defn}{\protect\definitionname}
\theoremstyle{remark}
\newtheorem{rem}{\protect\remarkname}
\theoremstyle{plain}
\newtheorem{thm}{\protect\theoremname}
\theoremstyle{remark}
\newtheorem{notation}{\protect\notationname}
\theoremstyle{plain}
\newtheorem{cor}{\protect\corollaryname}
\theoremstyle{plain}
\newtheorem{lem}{\protect\lemmaname}
\theoremstyle{plain}
\newtheorem*{thm*}{\protect\theoremname}
\theoremstyle{plain}
\newtheorem{prop}{\protect\propositionname}
\theoremstyle{remark}

\theoremstyle{remark}
\newtheorem*{claim*}{\protect\claimname}

\makeatother

\usepackage{babel}
\providecommand{\corollaryname}{Corollary}
\providecommand{\definitionname}{Definition}
\providecommand{\lemmaname}{Lemma}
\providecommand{\notationname}{Notation}
\providecommand{\remarkname}{Remark}
\providecommand{\theoremname}{Theorem}
\providecommand{\propositionname}{Proposition}
\providecommand{\claimname}{Claim}

\providecommand{\keywords}[1]
{
  \small	
  \textbf{\textbf{Keywords}} #1
}

\begin{document}
\title{A Linear Structure from Magnetic-Dipole Systems and Its Geometry}

\author[1,2]{Bohuan Lin\thanks{\url{bohuan.lin@xjtlu.edu.cn}}}
\author[3]{Fengping Li\thanks{wzulfp@126.com}}
\author[4]{Zhengya Zhang\thanks{\url{zhangzhengyachina@gmail.com}, \textbf{corresponding author}.}}

\affil[1]{Department of Foundational Mathematics, School of Mathematics and Physics, Xi'an Jiaotong-Liverpool University, 215123, Suzhou, China.}
\affil[2]{Department of Applied Mathematics, School of Mathematics and Physics, Xi'an Jiaotong-Liverpool University, 215123, Suzhou, China.}
\affil[3]{Oujiang Laboratory (Zhejiang Lab for Regenerative Medicine, Vision and Brain Health), 325000, Wenzhou, China}
\affil[4]{School of Robot Engineering, Wenzhou University of Technology, 325000, Wenzhou, China.}

\maketitle

\begin{abstract}
    This work has been motivated by the study of magnetic-dipole systems in micro-robotics. We define and investigate a class of algebras on $\mathbb{R}^3$ arising and generalized from the algebraic structure of magnetic gradient fields induced by systems of synchronous magnets with identical dipole moments (i.e., $\mathbf{M}_i=\mathbf{M},\,\forall i$). We show that when there is a $2$ dimensional sub-algebra, the linear map associated to the multiplication of the algebra admits a certain type of decompositions, which can be applied to an actual dipole system to locate the dipole moment $\bar{\mathbf{M}}$ that yields the strongest translational force(s) on a test magnet $\mathfrak{m}$. Upper and lower bounds to the largest translational force(s) are established, which can be applied in practice for fast estimation.
\end{abstract}



\section{Introduction}
\label{sec:introduction}

In this paper we present a study on a class of possibly non-associative algebras whose definition has been motivated by the specific structure of magnetic gradient fields generated by a group of magnets. This mainly concerns a map
\[
\mathbb{R}^{3}\ni\mathbf{M}\xmapsto{\ \ \mathbb{F}\ \ }\mathbb{F}_{\mathbf{M}}\in\mathrm{SyM}_{3}^{tr=0}
\]
with the following three properties: 
\begin{enumerate}
\item Linearity: $\mathbb{F}$ is a (nontrivial) linear map from $\mathbb{R}^{3}$ to
the space $\mathrm{SyM}_{3}^{tr=0}$ of all $3\times3$ real symmetric
matrices with the zero trace;
\item Reciprocity/Commutativity: $\mathbb{F}_{\mathbf{M}}\mathbf{m}=\mathbb{F}_{\mathbf{m}}\mathbf{M}$
for any $\mathbf{M},\mathbf{m}\in\mathbb{R}^{3}$;
\item Planarity: there is subspace $\mathcal{P}$ of dimension $2$ such
that $\mathbb{F}_{\mathbf{M}}\mathbf{m}\in\mathcal{P}$ holds for
any $\mathbf{M},\mathbf{m}\in\mathcal{P}$.
\end{enumerate}
This is an abstraction and generalization for the specific structure given below by \eqref{eq:F(M)}, which expresses the magnetic gradient field of a system of magnets ($\mathfrak{M}_i$) with identical dipole moments $\mathbf{M}_i=\mathbf{M}$. The interpretation in physics of the quantity $\mathbb{F}_{\mathbf{M}}\mathbf{m}$ is then the translational force experienced by a test magnet ($\mathfrak{m}$) with dipole moment $\mathbf{m}$.
To be accurate, with the
test magnet $\mathfrak{m}$ placed at $p\in\mathbb{R}^{3}$ and the
system's magnets $\mathfrak{M}_{i}$ placed at $p_{i}\in\mathbb{R}^{3}$,
the translational force experienced by the test magnet $\mathfrak{m}$
is $\mathbf{F}=\nabla\mathbf{B}\cdot\mathbf{m}=\frac{3\mu_{0}}{4\pi}\mathbb{F}_{\mathbf{M}}\mathbf{m}$ 
with (e.g., see \cite{pittiglio2022collaborative,davy2025utilizing})
\begin{equation}
\mathbb{F}_{\mathbf{M}}=\mathrm{P}\mathbf{M}^{\intercal}+\mathbf{M}\mathrm{P}^{\intercal}+\big(\mathbf{M}^{\intercal}\mathrm{P}\big)\cdot[\mathbf{I}]-5\sum_{i=1}^{n}\frac{\big(\mathbf{M}^{\intercal}\hat{\mathrm{p}}_{i}\big)\cdot\hat{\mathrm{p}}_{i}\hat{\mathrm{p}}_{i}^{\intercal}}{||\mathrm{p}_{i}||^{4}},\label{eq:F(M)}
\end{equation}
where $\hat{\mathrm{p}}_{i}=\frac{\mathrm{p}_{i}}{||\mathrm{p}_{i}||}$
with $\mathrm{p}_{i}=p-p_{i}$, $\mathrm{P}=\underset{i=1}{\overset{n}{\sum}}\frac{\hat{\mathrm{p}}_{i}}{||\mathrm{p}_{i}||^{4}}$,
and, $[\mathbf{I}]$ denotes the identity matrix. To ease the reading we also provide the derivation of Eq.\eqref{eq:F(M)} in the appendix. 
The linearity of the mapping $\mathbf{M}\mapsto\mathbb{F}_{\mathbf{M}}$ and the reciprocity $\mathbb{F}_{\mathbf{M}}\mathbf{m}=\mathbb{F}_{\mathbf{m}}\mathbf{M}$ follow directly from
(\ref{eq:F(M)}), while planarity is not a property demonstrated in
all the "actual" systems of synchronous magnets. One circumstance under which the $\mathbb{F}$ given by \eqref{eq:F(M)} has the
property of planarity is when all the magnets, $\mathfrak{M}_{i}$
and $\mathfrak{m}$, are placed in the same plane. In particular,
if the system consists of only two magnets, $\mathfrak{M}_{\pm}$
at $p_{\pm}$, then all the three magnets, $\mathfrak{M}_{\pm}$ and
$\mathfrak{m}$ are always in the same plane. 

This study has been
is an extension in mathematics of the previous works on the largest magnitude of the gradient force produced by a pair of rotating magnets with identical (but time-varying) dipole moments \cite{Zhengya2025AccurateWorstMagGradForce,Zhengya2025EstimatedWorstMagGradForce,Zhengya2025onPlaneWorstMagGradForce}, which have been motivated by the development of micro-robots actuated by a system of rotating magnets ($\mathfrak{M}_{i}$) \cite{mahoney2012control, hosney2015propulsion,zhang2025design}. It is a typical design that a helical micro-robot carries
a magnet ($\mathfrak{m}$) and propels itself through the rotation
caused by the rotating magnetic field generated by the rotating magnets
($\mathfrak{M}_{i}$)\cite{dong2022magnetic}. In such a robotic system, the torque exerted
on the robot by the magnetic field serves as the main control of the
robot's motion, while the translational force due to the gradient
of the field is undesirable since it disturbs the control\cite{alshafeei2014magnetic}. This has
led to the investigation on the strongest translational
force exerted by the actuating dipole $\mathbf{M}$ on the test dipole $\mathbf{m}$ at an arbitrary field point $p$, which, thanks to the symmetry of the matrices $\mathbb{F}_{\mathbf{M}}$, corresponds to the maximum of all the eigenvalues
of $\mathbb{F}_{\mathbf{M}}$ for all $\mathbf{M}\in S^{2}$ \cite{mahoney2013managing}. To be
accurate, the quantity below is of concern
\begin{equation}\label{eq:max-Lambda}
    \bar{\lambda}:=\underset{\mathbf{M,m}\in S^2}{\max}||\mathbb{F}_{\mathbf{M}}\mathbf{m}||=\max\{|\lambda|\geq0|\ \exists\mathbf{M}\in S^{2}\text{ such that }\lambda\text{ is an eigenvalue of }\mathbb{F}_{\mathbf{M}}\},
\end{equation}
with which the strongest translational force(s) produced at the field point $p$ by the actuating magnets has (have) the magnitude $\max||\mathbf{F}||=\frac{3\mu_{0}||\mathbf{M}||\cdot||\mathbf{m}||}{4\pi}\bar{\lambda}$.  The corresponding configurations $(\bar{\mathbf{M}},\bar{\mathbf{m}})$ of the actuating and the test dipoles are also of importance in practice, and the set of these configurations is:
\begin{equation}\label{eq:argmaxForce}
    \underset{(\mathbf{M},\mathbf{m})\in S^2\times S^2}{\arg\max}||\mathbb{F_{\mathbf{M}}\mathbf{m}}||=\{(\bar{\mathbf{M}},\bar{\mathbf{m}})\in S^2\times S^2\big|\;||\mathbb{F_{\bar{\mathbf{M}}}\bar{\mathbf{m}}}||=\bar{\lambda}\}.
\end{equation}
\begin{rem}\label{rem:LargestEigenLambda}
    Now that  $\mathbb{F}_{\bar{\mathbf{M}}}$ is a symmetric matrix, either $\bar{\lambda}$ or $-\bar{\lambda}$ is the principal eigenvalue of $\mathbb{F}_{\bar{\mathbf{M}}}$.
\end{rem}

In order to mitigate the undesirable translational magnetic force and enhance the controllability, developers of magnetic-torque-driven micro-robots deployed a pair
of synchronized magnets ($\mathfrak{M}_{\pm}$) 
for the actuation\cite{alshafeei2014magnetic, hosney2015propulsion,pittiglio2022collaborative,davy2025utilizing,zhang2025design}, and then the corresponding linear structure
$\mathbb{F}$ has the property of planarity. It turns out that this
property of $\mathbb{F}$ plays an important role in the search for
the strongest translational force, as well as for locating the dipole
moment(s) $\bar{\mathbf{M}}$. Moreover, as is to be shown later, the
property of planarity is not limited to systems with all the magnets
($\mathfrak{M}_{i}$ and $\mathfrak{m}$) placed in the same plane,
but is also possessed by those non-planar systems in which the distribution
of the actuating magnets ($\mathfrak{M}_{i}$) demonstrates certain
symmetry. One example of such systems is an infinite cubic lattice
of magnets with identical dipole moments. These rationalize the study
of this paper and promise a broad scope of the research. 

As a linear map from $\mathbb{R}^{3}$ to $\mathrm{SyM}_{3}^{tr=0}$,
$\mathbb{F}$ induces a multiplication on $\mathbb{R}^{3}$ which
obeys the distributive law with the addition of vectors, making
$\mathbb{R}^{3}$ into an algebra with possible lack of
associativity for the multiplication.  Indeed, non-associative structures appear from many different places in physics, and for a recent review we refer to \cite{szabo2019introduction}. However, to the best of our knowledge, there is still a lack of investigation on this specific type of algebras with the same perspective as we do in this paper.

\vspace{0.2cm}
We adopt the following definitions of an abstract magnetic algebra and its planarity.
\begin{defn}
\label{def:MagneticAlgebra}A magnetic algebra is a pair $(\mathbb{R}^{3},\mathbb{F})$
in which $\mathbb{F}:\mathbb{R}^{3}\rightarrow\mathrm{SyM}_{3}^{tr=0}$
is a (nontrivial) linear map with reciprocity, i.e., $\mathbb{F}_{\mathbf{M}}\mathbf{m}=\mathbb{F}_{\mathbf{m}}\mathbf{M}$
for all $\mathbf{M},\mathbf{m}\in\mathbb{R}^{3}$. 
\end{defn}
\begin{defn}
\label{def:PlanarConfiguration}A magnetic algebra $(\mathbb{R}^{3},\mathbb{F})$
has planar configuration if $\mathbb{F}$ has the property of planarity,
and the plane $\mathcal{P}$ is called an $\mathbb{F}$-invariant
plane or an $\mathbb{F}$-subalgebra.
\end{defn}
\begin{rem}
It is worth mentioning the difference between the unit of the outputs $\mathbb{F}_{\mathbf{M}}\mathbf{m}$
(Newton) and the unit of the inputs $\mathbf{M},\mathbf{m}$ (Joule
per Tesla) in the "actual" algebra given by \eqref{eq:F(M)}.
\end{rem}
It is observed that the general setting of a magnetic algebra with
planar configuration given in Definitions \ref{def:MagneticAlgebra}
and \ref{def:PlanarConfiguration} already captures certain important
features of an ``actual'' magnetic algebra given explicitly by (\ref{eq:F(M)}).
The main purpose of the paper is to demonstrate this discovery,  
and the main results are presented in Sections \ref{sec:AlgebraicStructure} to \ref{sec:GeometricBound}, following an introduction to basic properties about magnetic algebras in Section \ref{sec:Fundamentals}.

Section \ref{sec:AlgebraicStructure} concerns the algebraic structure of an
abstract magnetic algebra $(\mathbb{R}^{3},\mathbb{F})$ with planar
configuration. Observe from (\ref{eq:F(M)}) that, for an actual magnet
system with all the magnets ($\mathfrak{M}_{i}$ and $\mathfrak{m}$)
lying in a single plane (parallel to the subspace $\mathcal{P}$),
the operator $\mathbb{F}$ admits a decomposition $\mathbb{F}=\mathbb{E}-\mathbb{P}$
with
\[
\mathbb{E}_{\mathbf{M}}=\big[\mathrm{P}\mathbf{M}^{\intercal}+\mathbf{M}\mathrm{P}^{\intercal}\big]+\big(\mathbf{M}\cdot\mathrm{P}\big)\mathbf{Id}
\]
and $\mathbb{P}_{\mathbf{M}}=\mathbb{E}_{\mathbf{M}}-\mathbb{F}_{\mathbf{M}}$
for all $\mathbf{M}\in\mathbb{R}^{3}$. The main features of $\mathbb{E}_{\mathbf{M}}$
and $\mathbb{P}_{\mathbf{M}}$ are that, $\mathrm{Im}\mathbb{P}_{\mathbf{M}}\subset\mathcal{P}$,
and, $\mathbb{E}_{\mathrm{C}\mathbf{M}}=\mathrm{C}\mathbb{E}_{\mathbf{M}}\mathrm{C}^{\mathrm{T}}$
for any rotation $\mathrm{C}$ around $\mathrm{P}$. This structure is important
for locating $\bar{\mathbf{M}}\in S^{2}$ and finding $\bar{\lambda}$
for an actual system, and it is natural to ask if an abstract magnetic algebra
also has such a decomposition. 

The answer is affirmative. In fact, for any (abstract) magnetic algebra
$(\mathbb{R}^{3},\mathbb{F})$ with planar configuration, this kind
of decompositions of $\mathbb{F}$ exist and can be parameterized
by $\mathbb{R}$ in a particular way. For clearity, we denote by $\mathrm{Rot}_{\mathrm{P}}$
the subgroup of $\mathbb{SO}(3)$ with all the rotations fixing vector
$\mathrm{P}$, and by $\mathrm{SyM}_{3}$ the set of $3\times3$ real
symmetric matrices, and introduce the following definition:
\begin{defn}
\label{def:PolarDecomposition}For a(n) (abstract) magnetic algebra
$(\mathbb{R}^{3},\mathbb{F})$ with an $\mathbb{F}$-invariant plane
$\mathcal{P}$, $\mathbb{F}=\mathbb{E}-\mathbb{P}$ is called a $\mathcal{P}$-decomposition
if $\mathbb{E}$, $\mathbb{P}$ are linear maps from $\mathbb{R}^{3}$
to $\mathrm{SyM}_{3}$, such that $\mathrm{Im}\mathbb{P}_{\mathbf{M}}\subset\mathcal{P}$,
and, $\mathbb{E}_{\mathrm{C}\mathbf{M}}=\mathrm{C}\mathbb{E}_{\mathbf{M}}\mathrm{C}^{\mathrm{T}}$
for any $\mathrm{C}\in\mathrm{Rot}_{\mathrm{P}}$ and $\mathbf{M}\in\mathbb{R}^3$, where $\mathrm{P}=\mathbb{F}_{\bar{\mathbf{n}}}\bar{\mathbf{n}}$
with $\bar{\mathbf{n}}$ being an unit normal vector of $\mathcal{P}$.
\end{defn}
Note that in Definition \ref{def:PolarDecomposition}, the vector
$\mathrm{P}$ explicitly given by $\mathrm{P}=\underset{i}{\sum}\frac{\bar{\mathrm{p}}_{i}}{||\mathrm{p}_{i}||^{4}}$
in an actual magnetic system is defined to be $\mathrm{P}=\mathbb{F}_{\bar{\mathbf{n}}}\bar{\mathbf{n}}$
for an abstract algebra. The main result about the structure of $\mathbb{F}$ given in Section \ref{sec:AlgebraicStructure} is the theorem below:
\begin{thm}
\label{thm:=00005B=00005CP-decomposition=00005D}Let $(\mathbb{R}^{3},\mathbb{F})$
be a magnetic algebra with an $\mathbb{F}$-invariant plane $\mathcal{P}$. If $||\mathrm{P}||\neq0$, then $\mathbb{F}$ admits $\mathcal{P}$-decompositions. Moreover, every $\mathcal{P}$-decomposition of $\mathbb{F}$ takes the form of $\mathbb{F}=\mathbb{E}^{\xi}-\mathbb{P}^{\xi}$ for some $\xi\in\mathbb{R}$, in which
\begin{equation}\label{eq:=00005Bparameterization=00005D-=00005CE}
\mathbb{E}_{\mathbf{M}}^{\xi}=\big[\mathrm{P}\mathbf{M}^{\intercal}+\mathbf{M}\mathrm{P}^{\intercal}\big]+\big(\mathbf{M}\cdot\mathrm{P}\big)\mathbf{Id}-\xi\big(\mathbf{M}\cdot\mathrm{P}\big)\mathrm{P}\mathrm{P}^{\intercal},
\end{equation}
and,
\begin{equation}\label{eq:P_n=[0]}
\forall\mathbf{n}\perp\mathcal{P},\;\;
    \mathbb{P}^{\xi}_{\mathbf{n}}=[\mathbf{0}]\;(\text{the zero matrix}).
\end{equation}

\end{thm}

Based on Theorem \ref{thm:=00005B=00005CP-decomposition=00005D} as well as the other structural results about $(\mathbb{R}^3,\mathbb{F})$ obtained in Section \ref{sec:AlgebraicStructure}, we establish an analytical theory in Section \ref{sec:argmax-Locating} for calculating $\bar{\lambda}$ and locating the corresponding configurations $(\bar{\mathbf{M}},\bar{\mathbf{m}})$. While in practice it may be convenient to apply numerical methods, the analytical approach yields more insightful results and provides an accurate description of the configurations $(\bar{\mathbf{M}},\bar{\mathbf{m}})$ with $||\mathbb{F}_{\bar{\mathbf{M}}}\bar{\mathbf{m}}||=\bar{\lambda}$. For convenience, we adopt the following notation:
\begin{notation}
For any $\mathbf{M}\in\mathbb{R}^{3}$, denote by $\lambda_{\mathbf{M}}$
the principal eigenvalue of $\mathbb{F}_{\mathbf{M}}$. The other
eigenvalues are then $-\frac{\lambda_{\mathbf{M}}}{2}\pm\delta_{\mathbf{M}}$
with $\delta_{\mathbf{M}}\in[0,|\frac{\lambda_{\mathbf{M}}}{2}|]$. 
\end{notation}

$\;\;\;$
As is mentioned in Remark \ref{rem:LargestEigenLambda}, $\bar{\lambda}$ equals the magnitude of the principal eigenvalue of $\mathbb{F}_{\bar{\mathbf{M}}}$, i.e., $\bar{\lambda}=|\lambda_{\bar{\mathbf{M}}}|$, and then $\bar{\mathbf{m}}$ is either an associated eigenvector with $\lambda_{\bar{\mathbf{M}}}$ (if $\delta_{\bar{\mathbf{M}}}<\frac{\bar{\lambda}}{2}$), or, lies in the plane spanned by eigenvectors associated to the eigenvalues $\pm\bar{\lambda}$ (if $\delta_{\bar{\mathbf{M}}}=\frac{\bar{\lambda}}{2}$). Therefore, for understanding the set of the pairs $(\bar{\mathbf{M}},\bar{\mathbf{m}})$ in \eqref{eq:argmaxForce}, it suffices to look into the set of dipoles below:
\begin{equation}\label{eq:argmax|lambda|}
    \overline{\mathcal{M}}:=\underset{\mathbf{M}\in S^2}{\arg\max}{|\lambda_{\mathbf{M}}|}=\{\bar{\mathbf{M}}\in S^2\big|\;\bar{\lambda}=|\lambda_{\bar{\mathbf{M}}}|\}.
\end{equation}
Note that, in addition to its meaning in physics, the quantity $\bar{\lambda}$ also characterizes the ``geometry'' of the algebra $(\mathbb{R}^3,\mathbb{F})$: treating $\mathbb{F}_{\mathbf{M}}\mathbf{m}$ as the multiplication ($\mathbf{M}*\mathbf{m}$) on the algebra, $\bar{\lambda}=\underset{\mathbf{M,m}\in S^2}{\max}||\mathbf{M}*\mathbf{m}||$ shows the extent to which elements on $S^2$ can be ``stretched'' through the multiplication, and then $\overline{\mathcal{M}}$ contains exactly those the most stretched. Another subset of $S^2$ reflecting the geometry of $\mathbb{F}$  is
\[
\mathcal{M}_\mathbb{F}:=\underset{\mathbf{M}\in S^2}{\arg\max}||\mathbb{F}_{\mathbf{M}}||_{\mathbb{R}^9}.
\]
Here $||\mathbb{F}_{\mathbf{M}}||_{\mathbb{R}^9}$ is the Frobenius norm of the $3\times3$ matrix $\mathbb{F}_{\mathbf{M}}$, which is to treat  $\mathbb{F}_{\mathbf{M}}$
as an element in $\mathbb{R}^9$ and take the standard Euclidean norm. Viewing  $\mathbb{F}$ as a linear map from $\mathbb{R}^{3}$
to $\mathbb{R}^{9}$, it can be identified with a $9\times3$ matrix $[\mathbb{F}]$, and then for
$\mathbf{M}_0,\mathbf{M}_1\in\mathbb{R}^{3}$, the inner product of $\mathbb{F}_{\mathbf{M}_0}=[\mathbb{F}]\mathbf{M}_0$ and $\mathbb{F}_{\mathbf{M}_1}=[\mathbb{F}]\mathbf{M}_1$ in $\mathbb{R}^9$ is
\begin{equation}\label{eq:FrobeniusInnerProduct}
    \mathbf{M}_0^{\intercal}[\mathbb{F}^{\intercal}\mathbb{F}]\mathbf{M}_1=\mathrm{tr}[\mathbb{F}^\intercal_{\mathbf{M}_0}\mathbb{F}_{\mathbf{M}_1}]=\mathrm{tr}[\mathbb{F}_{\mathbf{M}_0}\mathbb{F}_{\mathbf{M}_1}].
\end{equation}
Therefore, the Frobenius norm
is just
\[
||\mathbb{F}_{\mathbf{M}}||_{\mathbb{R}^{9}}^{2}=\mathbf{M}^{\intercal}[\mathbb{F}^{\intercal}\mathbb{F}]\mathbf{M}=\text{tr}\mathbb{F}_{\mathbf{M}}^{2}.
\]
Here $\mathbb{F}^{\intercal}\mathbb{F}$ is a $3\times3$ semi-positive
definite matrix, and its principal eigenvalue is then 
\[\lambda_{\mathbb{F}}=\underset{\mathbf{M}\in S^2}{\max}||\mathbb{F}_{\mathbf{M}}||_{\mathbb{R}^{9}}^{2}.\]
$\mathcal{M}_{\mathbb{F}}$ is then just the set of all unit eigenvectors of $\mathbb{F}^\intercal\mathbb{F}$ associated to $\lambda_{\mathbb{F}}$, and therefore
\[
\lambda_{\mathbb{F}}=||\mathbb{F}_{\mathbf{M}_{\mathbb{F}}}||_{\mathbb{R}^{9}}^{2}=\text{tr}\mathbb{F}_{\mathbf{M}_{\mathbb{F}}}^{2},\;\forall\mathbf{M}_{\mathbb{F}}\in\mathcal{M}_\mathbb{F}.
\]
With the help of $\mathcal{M}_{\mathbb{F}}$ as well as the $\mathcal{P}$ decomposition $\mathbb{F}=\mathbb{E}-\mathbb{P}$ (with $\xi=0$), we obtain in Section \ref{sec:argmax-Locating} the following analytical result, which serves as the foundation for an algebraic approach to compute the set  $\overline{\mathcal{M}}$: 

\begin{thm}\label{thm:analytic-argmax}
Let $\mathrm{Det}^{\mathbb{F}}_0=\{\mathbf{M}\in\mathbb{R}^3\big|\;\det\mathbb{F}_{\mathbf{M}}=0\}$ and $\mathrm{Q}=\bar{\mathbf{n}}\times\mathrm{P}$. 
If $\mathrm{Det}^{\mathbb{F}}_0\cap\overline{\mathcal{M}}\neq\emptyset$, then $\overline{\mathcal{M}}=\mathcal{M}_{\mathbb{F}}=\mathrm{span}\{\bar{\mathbf{n}},\mathrm{P}\}$, $\bar{\lambda}=||\mathrm{P}||$, and, 
    \begin{equation}\label{eq:FwhenDet=0}
         \mathbb{F}_{\mathbf{M}}=(\mathbf{M}\cdot\bar{\mathbf{n}})[\mathrm{P}\bar{\mathbf{n}}^\intercal+\bar{\mathbf{n}}\mathrm{P}^\intercal]+(\mathrm{M}\cdot\mathrm{P})[\bar{\mathbf{n}}\bar{\mathbf{n}}^\intercal-\bar{\mathrm{P}}\bar{\mathrm{P}}^\intercal],\;\;\forall\mathbf{M}\in\mathbb{R}^3.
    \end{equation}
Otherwise, $\overline{\mathcal{M}}\subset\mathcal{A}\cup\mathcal{B}$ and $\mathbb{F}_{\bar{\mathbf{M}}}\bar{\mathbf{M}}=\pm\bar{\lambda}\bar{\mathbf{M}}$ for all $\bar{\mathbf{M}}\in\overline{\mathcal{M}}$. Set $\mathcal{A}$ contains elements at most four elements, each of which takes the form of
\begin{equation}\label{eq:Solution-to-[P_M]M=-P}
\hat{\mathbf{M}}=\cos\hat{\theta}\cdot\bar{\mathbf{n}}+\sin\hat{\theta}\cdot\frac{\mathbf{m}_+\pm\mathbf{m}_-}{||\mathbf{m}_+\pm\mathbf{m}_-||},
\end{equation}
where $\mathbf{m}_{\pm}\in\mathcal{P}\cap S^2$ are orthogonal eigenvectors of $\mathbb{P}_{\bar{\mathrm{Q}}}$, and, the value of $\hat\theta$ and the sign in $\pm$ can be determined via
\begin{equation}\label{eq:coeff-[P_M]M=-P}
\sin^2\hat{\theta}\cdot\big(\frac{5}{2}||\mathrm{P}||\pm\mathbf{m}^\intercal_+\mathbb{P}_{\bar{\mathrm{P}}}\mathbf{m}_-\big)=||\mathrm{P}||.
\end{equation}
Set $\mathcal{B}$ is the finite subset of $\mathcal{P}\cap S^2$ whose elements are $\mathbf{M}=\pm\frac{\mathbf{m}}{||\mathbf{m}||}$ with $\mathbf{m}\in\mathcal{P}$ being the solutions to the quadratic equations
\begin{equation}\label{eq:ConicCurves-P&Q}
    \mathbf{m}^\intercal\mathbb{F}_{\mathrm{P}}\mathbf{m}-\mathrm{P}^\intercal\mathbf{m}=0 \;\;\;\text{and}\;\;\;
    \mathbf{m}^\intercal\mathbb{F}_{\mathrm{Q}}\mathbf{m}-\mathrm{Q}^\intercal\mathbf{m}=0.
\end{equation}
Furthermore, when $\bar{\mathbf{M}}\in\mathcal{A}\cap\overline{\mathcal{M}}$, it holds $\bar{\lambda}=2|\bar{\mathbf{M}}\cdot\mathrm{P}|$; when $\bar{\mathbf{M}}\in\mathcal{B}\cap\overline{\mathcal{M}}$, $\bar{\mathbf{M}}=\pm\frac{\hat{\mathbf{m}}}{||\hat{\mathbf{m}}||}$ with $\hat{\mathbf{m}}$ solving \eqref{eq:ConicCurves-P&Q} and $\bar{\lambda}=\frac{1}{||\hat{\mathbf{m}}||}$.
\end{thm}

$\;\;\;$
Theorem \ref{thm:analytic-argmax} shows in a certain case $\overline{\mathcal{M}}=\mathcal{M}_\mathbb{F}$ and then $\lambda_{\mathbf{M}_\mathbb{F}}=\pm\bar{\lambda}$ for any $\mathbf{M}_{\mathbb{F}}\in\mathcal{M}_{\mathbb{F}}$. Indeed, while for any $\bar{\mathbf{M}}\in\overline{\mathcal{M}}$ the matrix $\mathbb{F}_{\bar{\mathbf{M}}}$ has the largest eigenvalue among all $\mathbb{F}_{\mathbf{M}}$, $\mathbb{F}_{\mathbf{M}_\mathbb{F}}$ has the largest Frobenius norm $\sqrt{\lambda_{\mathbb{F}}}$. Now that $||\mathbb{F}_{\mathbf{M}}||^2$ equals the sum of the squares of the eigenvalues, $\lambda_{\mathbf{M}_{\mathbb{F}}}$ of $\mathbb{F}_{\mathbf{M}_{\mathbb{F}}}$ should be ``close'' to $\lambda_{\bar{\mathbf{M}}}=\pm\bar{\lambda}$ of $\mathbb{F}_{\bar{\mathbf{M}}}$ to some extent. In Section \ref{sec:GeometricBound} we continue to explore the relation between $\lambda_{\mathbf{M}_{\mathbf{F}}}$ and $\lambda_{\bar{\mathbf{M}}}$, which leads to an analytic estimation (given as Theorem \ref{thm:=00005BGeometric-Bound=00005D} below) for $\bar{\lambda}$ in terms of the "algebraic characteristics" of $\mathbb{F}$: $\lambda_{\mathbb{F}}$, $\mathbb{F}$,  $||\mathrm{P}||$, $\lambda_{\mathbf{M}_{\mathbb{F}}}$, and, 
\begin{notation}
$\lambda_{\mathcal{P}}:=\underset{\mathbf{M}\in S^{2}\cap\mathcal{P}}{\max}\lambda_{\mathbf{M}}$. 
\end{notation}

%
%

%
\begin{thm}
\label{thm:=00005BGeometric-Bound=00005D}Let
$(\mathbb{R}^{3},\mathbb{F})$ be a magnetic algebra with $\mathbb{F}$-invariant plane $\mathcal{P}$. It always holds
\begin{equation}
||\mathrm{P}||\leq|\lambda_{\mathbf{M}_{\mathbb{F}}}|\leq\lambda_{\mathcal{P}}\leq\bar{\lambda}\leq|\lambda_{\mathbf{M}_{\mathbb{F}}}|+\frac{||\mathrm{P}||}{2},\label{eq:BasicGeoInequation=00005BwithPlanarity=00005D}
\end{equation}
in which the upper bound to $\bar{\lambda}$ can be further refined: 

when $\lambda_{\mathcal{P}}\leq2||\mathrm{P}||$, 
\begin{equation}
\bar{\lambda}\leq\min\bigg\{|\lambda_{\mathbf{M}_{\mathbb{F}}}|+\frac{||\mathrm{P}||}{3},\,2||\mathrm{P}||\cdot\sqrt{\frac{||\mathrm{P}||}{3||\mathrm{P}||-\lambda_{\mathcal{P}}}}\bigg\};\label{eq:RefinedUpperBound-=00005B=00005Clambda<2|P|=00005D}
\end{equation}
when $\lambda_{\mathcal{P}}\geq2||\mathrm{P}||$, 
\begin{equation}
\bar{\lambda}=\lambda_{\mathcal{P}}<\frac{||\mathrm{P}||+\sqrt{2\lambda_{\mathbb{F}}-3||\mathrm{P}||^{2}}}{2}<\sqrt{\frac{2}{3}\lambda_{\mathbb{F}}}.\label{eq:RefinedUpperBound-=00005B=00005Clambda>2|P|=00005D}
\end{equation}
\end{thm}
When $\mathbb{F}$ is given with an explicitly formula, e.g., (\ref{eq:F(M)}),
the quantities $||\mathrm{P}||$, $\lambda_{\mathbb{F}}$
and $\lambda_{\mathbf{M}_{\mathbb{F}}}$ can be calculated arithmetically.
This means that the estimation given in Theorem \ref{thm:=00005BGeometric-Bound=00005D} can be computed efficiently, which is desirable in practice. As is mentioned before, for the control
of micro-robots by rotating magnets, the translational force due to
the magnetic gradient field is undesirable, and hence efficient methods
for estimating the maximal magnitude of this force will be of great
help for feedback control of the robot. 

In addition to Sections \ref{sec:AlgebraicStructure}, \ref{sec:argmax-Locating} and \ref{sec:GeometricBound} for establishing the main theorems mentioned above, we discuss in detail a specific example in Section \ref{sec:Example-F*F=uniform}, showing a certain case in which we can directly deduce $\overline{\mathcal{M}}$ and $\bar{\lambda}$ from $\mathbb{F}^\intercal\mathbb{F}$.

\section{Fundamentals of Magnetic Algebras}
\label{sec:Fundamentals}

\subsection{Preliminaries about Eigenvalues of $\mathbb{F}_{\mathbf{M}}$}

Consider a linear map $\mathbb{F}:\mathbb{R}^{3}\rightarrow\mathrm{SyM}_{3}^{tr=0}$.
It is straightforward to check that $\frac{\text{tr}\mathbb{F}_{\mathbf{M}}^{3}}{3}=\det\mathbb{F}_{\mathbf{M}}$
holds for each $\mathbf{M}\in\mathbb{R}^{3}$. Also check that
\begin{equation}
\begin{aligned}\text{tr}\mathbb{F}_{\mathbf{M}}^{2} &  & = & \lambda_{\mathbf{M}}^{2}+\bigg(\frac{\lambda_{\mathbf{M}}}{2}+\delta_{\mathbf{M}}\bigg)^{2}+\bigg(\frac{\lambda_{\mathbf{M}}}{2}-\delta_{\mathbf{M}}\bigg)^{2}\\
 &  & = & \frac{3}{2}\lambda_{\mathbf{M}}^{2}+2\delta_{\mathbf{M}}^{2},
\end{aligned}
\label{eq:trF^2_M}
\end{equation}
and
\[
\begin{aligned} & \bigg(\frac{\lambda_{\mathbf{M}}}{2}+\delta_{\mathbf{M}}\bigg)\cdot\bigg(\frac{\lambda_{\mathbf{M}}}{2}-\delta_{\mathbf{M}}\bigg)-\lambda_{\mathbf{M}}\bigg(\frac{\lambda_{\mathbf{M}}}{2}+\delta_{\mathbf{M}}\bigg)-\lambda_{\mathbf{M}}\bigg(\frac{\lambda_{\mathbf{M}}}{2}-\delta_{\mathbf{M}}\bigg)\\
= & -\frac{3}{4}\lambda_{\mathbf{M}}^{2}-\delta_{\mathbf{M}}^{2}\ \ =-\frac{\text{tr}\mathbb{F}_{\mathbf{M}}^{2}}{2}.
\end{aligned}
\]
As a result, the characteristic equation of $\mathbb{F}_{\mathbf{M}}\in\mathrm{SyM}_{3}^{tr=0}$
takes the form
\begin{equation}
t^{3}-\frac{\text{tr}\mathbb{F}_{\mathbf{M}}^{2}}{2}t-\frac{\text{tr}\mathbb{F}_{\mathbf{M}}^{3}}{3}=0.\label{eq:characteristicEQ}
\end{equation}

Let $r_{\mathbf{M}}\in[0,1]$ be the ratio for $\delta_{\mathbf{M}}=r_{\mathbf{M}}\frac{|\lambda_{\mathbf{M}}|}{2}$.
The following identity will be used in later discussion 
\begin{equation}
\text{tr}\mathbb{F}_{\mathbf{M}}^{2}=||\mathbb{F}_{\mathbf{M}}||_{\mathbb{R}^{9}}^{2}=\frac{3+r_{\mathbf{M}}^{2}}{2}\lambda_{\mathbf{M}}^{2},\ \forall\mathbf{M}\in\mathbb{R}^{3}.\label{eq:=00005Clamda^2=000026|F_M|^2}
\end{equation}

It follows directly from (\ref{eq:trF^2_M}) and $r_{\mathbf{M}}^{2}\leq1$
that 
\begin{equation}
\frac{\text{tr}\mathbb{F}_{\mathbf{M}}^{2}}{2}=\frac{3+r_{\mathbf{M}}^{2}}{4}\lambda_{\mathbf{M}}^{2}\leq\lambda_{\mathbf{M}}^{2}=\frac{2\cdot\text{tr}\mathbb{F}_{\mathbf{M}}^{2}}{3+r_{\mathbf{M}}^{2}}\leq\frac{2}{3}\text{tr}\mathbb{F}_{\mathbf{M}}^{2}.\label{eq:BasicInequation-=00005Clambda_M}
\end{equation}

\subsection{First Analysis for $\bar{\lambda}$}

Let $(\mathbb{R}^{3},\mathbb{F})$ be a magnetic algebra (not necessarily
with planarity). First of all, note that $\bar{\lambda}=\underset{\mathbf{M},\mathbf{m}\in S^{2}}{\max}||\mathbb{F}_{\mathbf{M}}\mathbf{m}||$
and it is achieved at some point $(\bar{\mathbf{M}},\bar{\mathbf{m}})\in S^{2}\times S^{2}$
due to the continuity of $(\mathbf{M},\mathbf{m})\mapsto||\mathbb{F}_{\mathbf{M}}\mathbf{m}||$
and the compactness of $S^{2}\times S^{2}$. Now that $\mathbb{F}_{\bar{\mathbf{M}}}$
is symmetric, we can further take $\bar{\mathbf{m}}$ to be an eigenvector
of $\mathbb{F}_{\bar{\mathbf{M}}}$, and hence $\bar{\lambda}=|\lambda_{\bar{\mathbf{M}}}|$.
In fact, by reciprocity, $\big|\big|\mathbb{F}_{\bar{\mathbf{m}}}\bar{\mathbf{M}}\big|\big|=\big|\big|\mathbb{F}_{\bar{\mathbf{M}}}\bar{\mathbf{m}}\big|\big|=\bar{\lambda}$,
and hence $|\lambda_{\bar{\mathbf{m}}}|=|\lambda_{\bar{\mathbf{M}}}|=\bar{\lambda}$.

We have the following basic estimation for $\bar{\lambda}$ 
\begin{thm}
\label{thm:BasicIneq}For any magnetic algebra $(\mathbb{R}^{3},\mathbb{F})$,
it holds
\begin{equation}
|\lambda_{\mathbf{M}_{\mathbb{F}}}|^{2}\leq\bar{\lambda}^{2}\leq\frac{2}{3}\big|\big|\mathbb{F}_{\mathbf{M}_{\mathbb{F}}}\big|\big|_{\mathbb{R}^{9}}^{2}=\frac{2}{3}\lambda_{\mathbb{F}}.\label{eq:(1st estimate)<=00005Cbar=00007B=00005Clambda=00007D<}
\end{equation}
\end{thm}
To see the upper bound in (\ref{eq:(1st estimate)<=00005Cbar=00007B=00005Clambda=00007D<}),
note that by definition, $||\mathbb{F}_{\mathbf{M}_{\mathbb{F}}}||$
is the largest among all the $\big|\big|\mathbb{F}_{\mathbf{M}}\big|\big|_{\mathbb{R}^{9}}^{2}$,
and in particular, $||\mathbb{F}_{\mathbf{M}_{\mathbb{F}}}||_{\mathbb{R}^{9}}^{2}\geq||\mathbb{F}_{\bar{\mathbf{M}}}||_{\mathbb{R}^{9}}^{2}$.
As a result of (\ref{eq:BasicInequation-=00005Clambda_M}), it holds
\[
\bar{\lambda}^{2}=|\lambda_{\bar{\mathbf{M}}}|^{2}\leq\frac{2}{3}||\mathbb{F}_{\bar{\mathbf{M}}}||_{\mathbb{R}^{9}}^{2}\leq\frac{2}{3}||\mathbb{F}_{\mathbf{M}_{\mathbb{F}}}||_{\mathbb{R}^{9}}^{2}.
\]
$|\lambda_{\mathbf{M}_{\mathbb{F}}}|^{2}$ being a lower bound for
$\bar{\lambda}^{2}$ follows directly from the definition of $\bar{\lambda}$.
To see the difference between $|\lambda_{\mathbf{M}_{\mathbb{F}}}|$
and $\bar{\lambda}$, we resort to (\ref{eq:=00005Clamda^2=000026|F_M|^2})
and calculate
\[
\bar{\lambda}^{2}-\lambda_{\mathbf{M}_{\mathbb{F}}}^{2}\leq\frac{2}{3}\big|\big|\mathbb{F}_{\mathbf{M}_{\mathbb{F}}}\big|\big|_{\mathbb{R}^{9}}^{2}-\lambda_{\mathbf{M}_{\mathbb{F}}}^{2}=\frac{r_{\mathbf{M}_{\mathbb{F}}}^{2}}{3}\lambda_{\mathbf{M}_{\mathbb{F}}}^{2},
\]
i.e., $(\bar{\lambda}-|\lambda_{\mathbf{M}_{\mathbb{F}}}|)(\bar{\lambda}+|\lambda_{\mathbf{M}_{\mathbb{F}}}|)\leq\frac{r_{\mathbf{M}_{\mathbb{F}}}^{2}}{3}\lambda_{\mathbf{M}_{\mathbb{F}}}^{2}$,
and then 
\begin{equation}
0\leq\bar{\lambda}-|\lambda_{\mathbf{M}_{\mathbb{F}}}|\leq\frac{r_{\mathbf{M}_{\mathbb{F}}}^{2}}{3}\frac{\lambda_{\mathbf{M}_{\mathbb{F}}}^{2}}{\bar{\lambda}+|\lambda_{\mathbf{M}_{\mathbb{F}}}|}\leq\frac{r_{\mathbf{M}_{\mathbb{F}}}^{2}}{6}|\lambda_{\mathbf{M}_{\mathbb{F}}}|.\label{Ineq:1stErrEsti-GenericPosition}
\end{equation}
From (\ref{Ineq:1stErrEsti-GenericPosition}) it is easy to see that
when $r_{\mathbf{M}_{\mathbb{F}}}=0$, i.e., $\delta_{\mathbf{M}_{\mathbb{F}}}=0$,
$\bar{\lambda}=|\lambda_{\mathbf{M}_{\mathbb{F}}}|$.

Another scenario in which $\bar{\lambda}=|\lambda_{\mathbf{M}_{\mathbb{F}}}|$
is when $\det\mathbb{F}_{\bar{\mathbf{M}}}\cdot\det\mathbb{F}_{\bar{\mathbf{m}}}=0$.
In this case, either $\det\mathbb{F}_{\bar{\mathbf{M}}}=0$ or $\det\mathbb{F}_{\bar{\mathbf{m}}}=0$.
Let's say, $\det\mathbb{F}_{\bar{\mathbf{M}}}=0$, and then it follows
from (\ref{eq:characteristicEQ}) that
\[
\lambda_{\bar{\mathbf{M}}}^{3}-\frac{\text{tr}\mathbb{F}_{\bar{\mathbf{M}}}^{2}}{2}\lambda_{\bar{\mathbf{M}}}=0,
\]
which together with (\ref{eq:BasicInequation-=00005Clambda_M}) implies
\[
|\lambda_{\mathbf{M}_{\mathbb{F}}}|\leq\bar{\lambda}=|\lambda_{\bar{\mathbf{M}}}|=\sqrt{\frac{\text{tr}\mathbb{F}_{\bar{\mathbf{M}}}^{2}}{2}}\leq\sqrt{\frac{\text{tr}\mathbb{F}_{\mathbf{M}_{\mathbb{F}}}^{2}}{2}}\leq|\lambda_{\mathbf{M}_{\mathbb{F}}}|.
\]
In fact, we have the following general result:
\begin{thm}
\label{thm:=00005Cbar=00007B=00005Clambda=00007D=00003D=00005Clambda_=00007BM_F=00007D}For
any $\hat{\mathbf{M}}\in S^{2}$ with $\lambda_{\hat{\mathbf{M}}}^{2}\geq\lambda_{\mathbf{M}_{\mathbb{F}}}^{2}$,
it holds $0\leq r_{\hat{\mathbf{M}}}\leq r_{\mathbf{M}_{\mathbb{F}}}\leq1$.
If $r_{\hat{\mathbf{M}}}=r_{\mathbf{M}_{\mathbb{F}}}$, then $|\lambda_{\hat{\mathbf{M}}}|=|\lambda_{\mathbf{M}_{\mathbb{F}}}|$ and $[\mathbb{F}^\intercal\mathbb{F}]\hat{\mathbf{M}}=\lambda_{\mathbb{F}}\hat{\mathbf{M}}$.
In particular, if $r_{\bar{\mathbf{M}}}=1$ (i.e. $\det\mathbb{F}_{\bar{\mathbf{M}}}=0$) or $r_{\mathbf{M}_{\mathbb{F}}}=0$ holds,
$\bar{\lambda}=|\lambda_{\mathbf{M}_{\mathbb{F}}}|$ and $[\mathbb{F}^\intercal\mathbb{F}]\bar{\mathbf{M}}=\lambda_{\mathbb{F}}\bar{\mathbf{M}}$. 
\end{thm}
\begin{proof}
Since $\lambda_{\hat{\mathbf{M}}}^{2}\geq\lambda_{\mathbf{M}_{\mathbb{F}}}^{2}$
and
\[
\frac{3+r_{\hat{\mathbf{M}}}^{2}}{2}\lambda_{\hat{\mathbf{M}}}^{2}=\text{tr}\mathbb{F}_{\hat{\mathbf{M}}}^{2}\leq\text{tr}\mathbb{F}_{\mathbf{M}_{\mathbb{F}}}^{2}=\frac{3+r_{\mathbf{M}_{\mathbb{F}}}^{2}}{2}\lambda_{\mathbf{M}_{\mathbb{F}}}^{2},
\]
it requires $r_{\hat{\mathbf{M}}}\leq r_{\mathbf{M}_{\mathbb{F}}}$.
Suppose that $r_{\hat{\mathbf{M}}}=r_{\mathbf{M}_{\mathbb{F}}}$.
Then $\lambda_{\hat{\mathbf{M}}}^{2}\geq\lambda_{\mathbf{M}_{\mathbb{F}}}^{2}$
implies
\[
\text{tr}\mathbb{F}_{\hat{\mathbf{M}}}^{2}=\frac{3+r_{\hat{\mathbf{M}}}^{2}}{2}\lambda_{\hat{\mathbf{M}}}^{2}\geq\frac{3+r_{\mathbf{M}_{\mathbb{F}}}^{2}}{2}\lambda_{\mathbf{M}_{\mathbb{F}}}^{2}=\text{tr}\mathbb{F}_{\mathbf{M}_{\mathbb{F}}}^{2},
\]
and as a result, $\hat{\mathbf{M}}^\intercal[\mathbb{F}^\intercal\mathbb{F}]\hat{\mathbf{M}}=\text{tr}\mathbb{F}_{\hat{\mathbf{M}}}^{2}=\text{tr}\mathbb{F}_{\mathbf{M}_{\mathbb{F}}}^{2}=\lambda_{\mathbb{F}}$,
which implies $\lambda_{\hat{\mathbf{M}}}^{2}=\lambda_{\mathbf{M}_{\mathbb{F}}}^{2}$ and that $\hat{\mathbf{M}}$ is also an eigenvector of $\mathbb{F}^\intercal\mathbb{F}$ with eigenvalue $\lambda_{\mathbb{F}}$.
\end{proof}
\begin{rem}\label{rem:detF=0}
In general, there can be more than one $\mathbf{M}_{\mathbb{F}}$ (unit eigenvector(s) of $\mathbb{F}^\intercal\mathbb{F}$ associated to $\lambda_{\mathbb{F}}$) with different $|\lambda_{\bar{\mathbf{M}}_\mathbb{F}}|$. However, the choice of $\mathbf{M}_\mathbb{F}$ in Theorem \ref{thm:=00005Cbar=00007B=00005Clambda=00007D=00003D=00005Clambda_=00007BM_F=00007D} is arbitrary, and hence the result means that when $\det\mathbb{F}_{\bar{\mathbf{M}}}=0$, all the $|\lambda_{\bar{\mathbf{M}}_\mathbb{F}}|$ should equal to $\bar{\lambda}$, and also, all $\det\mathbb{F}_{\mathbf{M}_{\mathbb{F}}}$ equal to $0$.
\end{rem}

\subsection{(Actual) Magnet Systems with Invariant Planes}

At each field point $p$, the magnetic algebra of an actual system
given by (\ref{eq:F(M)}) with only two synchronized magnets ($\mathfrak{M}_{\pm}$,
placed at $\mathrm{o}_{\pm}$ and both with dipole moment $\mathbf{M}$)
has an $\mathbb{F}$-invariant plane, which is spanned by $\mathbf{p}_{\pm}:=p-\mathrm{o}_{\pm}$.
For a system with $n$ magnets ($\mathfrak{M}_{i}$, $i=1,...,n$)
arbitrarily distributed in space, there is no garantee for the existence
of an $\mathbb{F}$-invariant plane. However, as is to be shown below,
when the distribution of the synchronized magnets demonstrates mirror
symmetry with respect to a plane $\mathcal{P}$, the algebras at certain
field points will have $\mathcal{P}$ as an invariant plane. 
\begin{notation}
By mirror symmetry w.r.t. $\mathcal{P}$ we mean the magnets are distributed
in such way that, if there is a magnet placed at $u+v$ with $u\in\mathcal{P}$
and $v\perp\mathcal{P}$, then correpondingly there is another magnet
at $u-v$.
\end{notation}
\begin{thm}
\label{thm:MagnetArray-with-MirrorSymmetry}Let $\mathcal{P}$ be
a two dimensional subspace of $\mathbb{R}^{3}$ with normal direction
$\bar{\mathbf{n}}$. Suppose that the synchronized magnets (with dipole
moment $\mathbf{M}$) are distributed in $\mathbb{R}^{3}$ with mirror
symmetry with respect to $\mathcal{P}$. Then, at any field point
$p\in\mathcal{P}$, $\mathcal{P}$ is an $\mathbb{F}$-invariant plane.
\end{thm}
\begin{proof}
Without loss of generality, we take the field point to be $p=\mathbf{0}$.
Due to the mirror symmetry with respect to $\mathcal{P}$, the magnets
can be classified into three groups: 1) the magnets placed at $\mathbf{q}_{j}\in\mathcal{P}$;
2) the magnets placed at $\mathbf{p}_{i}^{+}=\mathbf{w}_{i}+t_{i}\bar{\mathbf{n}}$
with $\mathbf{w}_{i}\in\mathcal{P}$ and $t_{i}>0$; 3) the magnets
at $\mathbf{p}_{i}^{-}=\mathbf{w}_{i}-t_{i}\bar{\mathbf{n}}$. The
last two groups contain the same number of magnets, and $||\mathbf{p}_{i}^{+}||=||\mathbf{p}_{i}^{-}||$
(denoted by $||\mathbf{p}_{i}||$ below).

With $\bar{\mathbf{p}}_{i}^{+}=\frac{\mathbf{p}_{i}^{+}}{||\mathbf{p}_{i}^{+}||}$
and $\bar{\mathbf{q}}_{j}=\frac{\mathbf{q}_{j}}{||\mathbf{q}_{j}||}$,
define $\mathrm{P}=\underset{i}{\sum}\frac{\bar{\mathbf{p}}_{i}^{+}}{||\mathbf{p}_{i}||^{4}}+\underset{i}{\sum}\frac{\bar{\mathbf{p}}_{i}^{-}}{||\mathbf{p}_{i}||^{4}}+\underset{j}{\sum}\frac{\bar{\mathbf{q}}_{j}}{||\mathbf{q}_{j}||^{4}}$.
The mapping $\mathbb{F}$ at the field point $p=\mathbf{0}$ is then
\[
\mathbb{F}_{\mathbf{M}}=\mathbf{M}\mathrm{P}^{\intercal}+\mathrm{P}\mathbf{M}^{\intercal}+(\mathbf{M}\cdot\mathrm{P})\mathbf{Id}-5\bigg[\underset{j}{\sum}\frac{\mathbf{M}\cdot\bar{\mathbf{q}}_{j}}{||\mathbf{q}_{j}||^{4}}[\bar{\mathbf{q}}_{j}][\bar{\mathbf{q}}_{j}]^{\intercal}+\underset{i}{\sum}\underset{\mathrm{l}=+,-}{\sum}\frac{\mathbf{M}\cdot\bar{\mathbf{p}}_{i}^{\mathrm{l}}}{||\mathbf{p}_{i}||^{4}}[\bar{\mathbf{p}}_{i}^{\mathrm{l}}][\bar{\mathbf{p}}_{i}^{\mathrm{l}}]^{\intercal}\bigg].
\]

Now we show that $\mathcal{P}$ is an $\mathbb{F}$-invariant plane.
First of all, check that
\[
\begin{aligned}\mathrm{P} &  & = & \underset{i}{\sum}\underset{\mathrm{l}=\pm1}{\sum}\frac{\mathbf{u}_{i}+\mathrm{l}\cdot\mathbf{v}_{i}}{||\mathbf{p}_{i}||^{4}}+\underset{j}{\sum}\frac{\bar{\mathbf{q}}_{j}}{||\mathbf{q}_{j}||^{4}}\\
 &  & = & \underset{i}{\sum}\underset{\mathrm{l}=\pm1}{\sum}\frac{\mathbf{u}_{i}}{||\mathbf{p}_{i}||^{4}}+\underset{j}{\sum}\frac{\bar{\mathbf{q}}_{j}}{||\mathbf{q}_{j}||^{4}}\\
 &  & \in & \mathcal{P}.
\end{aligned}
\]
Based on this, we split $\mathbb{F}_{\mathbf{M}}$ into two parts:
$\mathbb{F}_{\mathbf{M}}=\mathbb{A}_{\mathbf{M}}+\mathbb{B}_{\mathbf{M}}$
with
\[
\mathbb{A}_{\mathbf{M}}=\mathbf{M}\mathrm{P}^{\intercal}+\mathrm{P}\mathbf{M}^{\intercal}+(\mathbf{M}\cdot\mathrm{P})\mathbf{Id}-5\underset{j}{\sum}\frac{\mathbf{M}\cdot\bar{\mathbf{q}}_{j}}{|\mathbf{q}_{j}|^{4}}[\bar{\mathbf{q}}_{j}][\bar{\mathbf{q}}_{j}]^{\intercal}
\]
and 
\[
\mathbb{B}_{\mathbf{M}}=-5\underset{i}{\sum}\underset{\mathrm{l}=+,-}{\sum}\frac{\mathbf{M}\cdot\bar{\mathbf{p}}_{i}^{\mathrm{l}}}{||\mathbf{p}_{i}||^{4}}[\bar{\mathbf{p}}_{i}^{\mathrm{l}}][\bar{\mathbf{p}}_{i}^{\mathrm{l}}]^{\intercal}.
\]
It is easy to see that, for any $\mathbf{M},\mathbf{m}\in\mathcal{P}$,
$\mathbb{A}_{\mathbf{M}}\mathbf{m}\in\mathcal{P}$ holds, and then
it remains to show $\mathbb{B}_{\mathbf{M}}\mathbf{m}$ is also in
$\mathcal{P}$. To this end we shall first rewrite $\underset{i}{\sum}\underset{\mathrm{l}=+,-}{\sum}\frac{\mathbf{M}\cdot\bar{\mathbf{p}}_{i}^{\mathrm{l}}}{||\mathbf{p}_{i}||^{4}}[\bar{\mathbf{p}}_{i}^{\mathrm{l}}][\bar{\mathbf{p}}_{i}^{\mathrm{l}}]^{\intercal}$. 

For each $i$, let $\mathbf{u}_{i}=\frac{\mathbf{w}_{i}}{||\mathbf{p}_{i}||}$
and $\mathbf{v}_{i}=\frac{t_{i}\bar{\mathbf{n}}}{||\mathbf{p}_{i}||}$,
and then $\bar{\mathbf{p}}_{i}^{\pm}=\mathbf{u}_{i}\pm\mathbf{v}_{i}$
is the orthogonal decomposition with $\mathbf{u}_{i}\in\mathcal{P}$
and $\mathbf{v}_{i}\perp\mathcal{P}$. As a result,
\[
\begin{aligned}\underset{\mathrm{l}=+,-}{\sum}\frac{\mathbf{M}\cdot\bar{\mathbf{p}}_{i}^{\mathrm{l}}}{||\mathbf{p}_{i}||^{4}}[\bar{\mathbf{p}}_{i}^{\mathrm{l}}][\bar{\mathbf{p}}_{i}^{\mathrm{l}}]^{\intercal} &=\underset{\mathrm{l}=+,-}{\sum}\frac{\mathbf{M}\cdot\bar{\mathbf{p}}_{i}^{\mathrm{l}}}{||\mathbf{p}_{i}||^{4}}[\bar{\mathbf{p}}_{i}^{\mathrm{l}}][\bar{\mathbf{p}}_{i}^{\mathrm{l}}]^{\intercal}\\
 &=\underset{\mathrm{l}=\pm1}{\sum}\frac{\mathbf{M}\cdot(\mathbf{u}_{i}+\mathrm{l}\cdot\mathbf{v}_{i})}{||\mathbf{p}_{i}||^{4}}\bigg[\big(\mathbf{u}_{i}\mathbf{u}_{i}^{\intercal}+\mathbf{v}_{i}\mathbf{v}_{i}^{\intercal}\big)+\mathrm{l}\cdot\big(\mathbf{u}_{i}\mathbf{v}_{i}^{\intercal}+\mathbf{v}_{i}\mathbf{u}_{i}^{\intercal}\big)\bigg]\\
 &=2\bigg(\frac{\mathbf{M}\cdot\mathbf{u}_{i}}{||\mathbf{p}_{i}||^{4}}\big[\mathbf{u}_{i}\mathbf{u}_{i}^{\intercal}+\mathbf{v}_{i}\mathbf{v}_{i}^{\intercal}\big]+\frac{\mathbf{M}\cdot\mathbf{v}_{i}}{||\mathbf{p}_{i}||^{4}}\big[\mathbf{u}_{i}\mathbf{v}_{i}^{\intercal}+\mathbf{v}_{i}\mathbf{u}_{i}^{\intercal}\big]\bigg).
\end{aligned}
\]
Now, with $\mathbf{M},\mathbf{m}\in\mathcal{P}$, $\frac{\mathbf{M}\cdot\mathbf{v}_{i}}{||\mathbf{p}_{i}||^{4}}\big[\mathbf{u}_{i}\mathbf{v}_{i}^{\intercal}+\mathbf{v}_{i}\mathbf{u}_{i}^{\intercal}\big]=[\mathbf{0}]$
since $\mathbf{M}\cdot\mathbf{v}_{i}=0$, and also,
\[
\frac{\mathbf{M}\cdot\mathbf{u}_{i}}{||\mathbf{p}_{i}||^{4}}\big[\mathbf{u}_{i}\mathbf{u}_{i}^{\intercal}+\mathbf{v}_{i}\mathbf{v}_{i}^{\intercal}\big]\mathbf{m}=\frac{(\mathbf{M}\cdot\mathbf{u}_{i})\cdot(\mathbf{u}_{i}\cdot\mathbf{m})}{||\mathbf{p}_{i}||^{4}}\mathbf{u}_{i}\in\mathcal{P},
\]
which concludes the proof.
\end{proof}
\begin{rem}
As a consequence from Theorem \ref{thm:MagnetArray-with-MirrorSymmetry},
an actual magnetic algebra may have more than one $\mathbb{F}$-invariant
planes. For example, consider a system given by (\ref{eq:F(M)}) with
$8$ synchronized dipoles, each of which is placed at a vertex of
a cube. The distribution has mirror symmetry with respect to more
than one plane, each of which by Theorem \ref{thm:MagnetArray-with-MirrorSymmetry}
is (parallel to) an invariant plane of the algebra.
\end{rem}

\section{Structures of Magnetic Algebras with Planarity}
\label{sec:AlgebraicStructure}

\subsection{Basic Structures of $(\mathbb{R}^{3},\mathbb{F})$}

Equip $\mathbb{R}^{3}$ with the standard inner product $\mathbf{x}\cdot\mathbf{y}:=\mathbf{x}^{\intercal}\mathbf{y}$.
The orthogonal complement $\mathcal{P}^{\perp}$ of $\mathcal{P}$
has dimension $1$ and is spanned by $\bar{\mathbf{n}}$. It turns
out that $\bar{\mathbf{n}}$ is very special to the structure of $\mathbb{F}$.
Since $\mathcal{P}$ is an invariant subspace of the symmetric linear
map $\mathbb{F}_{\mathbf{M}}$ for $\mathbf{M}\in\mathcal{P}$, $\mathcal{P}^{\perp}$
is also an invariant subspace of $\mathbb{F}_{\mathbf{M}}$, and hence
$\bar{\mathbf{n}}$ is eigenvalue of $\mathbb{F}_{\mathbf{M}}$ because
of $\dim\mathcal{P}^{\perp}=1$. By the reciprocity of $\mathbb{F}$,
$\mathbb{F}_{\bar{\mathbf{n}}}\mathbf{M}=\mathbb{F}_{\mathbf{M}}\bar{\mathbf{n}}$
lies in $\mathcal{P}^{\perp}$, that is $\mathbb{F}_{\bar{\mathbf{n}}}\big(\mathcal{P}\big)\subset\mathcal{P}^{\perp}$.
As is shown below, this relation is actually an equivalent characterization
of the planarity of $\mathbb{F}$. 
\begin{thm}
\label{thm:=00005BEquiv-Planarity=00005D}Let $\mathbb{F}:\mathbb{R}^{3}\rightarrow\mathrm{SyM}_{3}^{tr=0}$
be a linear map with reciprocity $\mathbb{F}_{\mathbf{M}}\mathbf{m}=\mathbb{F}_{\mathbf{m}}\mathbf{M}$.
Given an unit vector $\bar{\mathbf{n}}$ with $\mathcal{P}=\big(\mathbb{R}\cdot\bar{\mathbf{n}}\big)^{\perp}$,
the following statements are equivalent: 1) $\mathcal{P}$ is an $\mathbb{F}$-invariant
plane; 2) $\mathbb{F}_{\bar{\mathbf{n}}}\big(\mathcal{P}\big)\subset\mathbb{R}\cdot\bar{\mathbf{n}}$;
3) $\exists\bar{\mathbf{n}}\in\mathbb{R}^{3}$ with $\bar{\mathbf{n}}\neq0$
such that $[\bar{\mathbf{n}}]\mathbb{F}_{\bar{\mathbf{n}}}[\bar{\mathbf{n}}]$
is the zero matrix $[\mathbf{0}]$.
\end{thm}
\begin{notation}
By $\mathbb{R}\cdot\bar{\mathbf{n}}$ we denote the subspace spanned
by $\bar{\mathbf{n}}$ and by $\big(\mathbb{R}\cdot\bar{\mathbf{n}}\big)^{\perp}$
its orthogonal complement. With $\bar{\mathbf{n}}=(n_{1},n_{2},n_{3})$,
$[\bar{\mathbf{n}}]$ stands for the matrix 
\[
\left[\begin{array}{ccc}
0 & -n_{3} & n_{2}\\
n_{3} & 0 & -n_{1}\\
-n_{2} & n_{1} & 0
\end{array}\right].
\]
\end{notation}
\begin{proof}
The equivalence between 2) and 3) is straightforward.

For 1)$\iff$2), we have shown in the previous discussion that $\mathcal{P}$
being an $\mathbb{F}$-invariant plane implies $\mathbb{F}_{\bar{\mathbf{n}}}\big(\mathcal{P}\big)\subset\mathbb{R}\cdot\bar{\mathbf{n}}$
for its normal vector $\bar{\mathbf{n}}$. Conversely, suppose that
$\mathbb{F}_{\bar{\mathbf{n}}}\big(\mathcal{P}\big)\subset\mathbb{R}\cdot\bar{\mathbf{n}}$.
Given any $\mathbf{M}\in\mathcal{P}$, $\mathbb{F}_{\mathbf{M}}\bar{\mathbf{n}}=\mathbb{F}_{\bar{\mathbf{n}}}\mathbf{M}\in\mathbb{R}\cdot\bar{\mathbf{n}}$,
i.e., $\mathbb{R}\cdot\bar{\mathbf{n}}$ is an invariant subspace
of $\mathbb{F}_{\mathbf{M}}$. Since $\mathbb{F}_{\mathbf{M}}$ is
a symmetic matrix, the orthogonal complement $\mathcal{P}$ is also
invariant under $\mathbb{F}_{\mathbf{M}}$, and hence it is an $\mathbb{F}$-subalgebra. 
\end{proof}
\begin{thm}
\label{thm:StruturalThm-PlaneConfig}Let $\bar{\mathbf{n}}\in\mathbb{R}^{3}$
be an unit vector and suppose that $[\bar{\mathbf{n}}]\mathbb{F}_{\bar{\mathbf{n}}}[\bar{\mathbf{n}}]=[\mathbf{0}]$.
Then, $\mathrm{P}:=\mathbb{F}_{\bar{\mathbf{n}}}\bar{\mathbf{n}}$
is perpendicular to $\mathcal{\bar{\mathbf{n}}}$, and 
\begin{equation}\label{eq:[F_n]}
\mathbb{F}_{\bar{\mathbf{n}}}=\mathrm{P}\bar{\mathbf{n}}^{\intercal}+\bar{\mathbf{n}}\mathrm{P}^{\intercal}.
\end{equation}
Moreover, for any $\mathbf{M}\perp\mathcal{\bar{\mathbf{n}}}$ it
holds $\mathbb{F}_{\mathbf{M}}\bar{\mathbf{n}}=(\mathrm{P}\cdot\mathbf{M})\bar{\mathbf{n}}$.
\end{thm}
\begin{proof}
From Theorem \ref{thm:=00005BEquiv-Planarity=00005D} we know that
$\mathcal{P}=\big(\mathbb{R}\cdot\bar{\mathbf{n}}\big)^{\perp}$ is
an $\mathbb{F}$-invariant plane. For the special case $\mathbb{F}_{\bar{\mathbf{n}}}\big(\mathcal{P}\big)=\{\mathbf{0}\}$,
$\text{tr}\mathbb{F}_{\bar{\mathbf{n}}}=0$ implies $\mathbb{F}_{\bar{\mathbf{n}}}$
to be the zero map, i.e., $\forall\mathbf{M}\in\mathbb{R}^{3}$, $\mathbb{F}_{\bar{\mathbf{n}}}\mathbf{M}=\mathbf{0}$,
and due to reciprocity, $\mathbb{F}_{\mathbf{M}}\bar{\mathbf{n}}=\mathbf{0}\in\mathbb{R}\cdot\bar{\mathbf{n}}$.
As a result, $\mathrm{P}:=\mathbb{F}_{\bar{\mathbf{n}}}\bar{\mathbf{n}}=\mathbf{0}$
and $\mathbb{F}_{\mathbf{M}}\bar{\mathbf{n}}=(\mathrm{P}\cdot\mathbf{M})\bar{\mathbf{n}}$.

Now consider the case $\mathbb{F}_{\bar{\mathbf{n}}}\big(\mathcal{P}\big)\neq\{\mathbf{0}\}$,
which means $\mathbb{F}_{\bar{\mathbf{n}}}\big(\mathcal{P}\big)=\mathbb{R}\cdot\bar{\mathbf{n}}$.
In this case, $\ker\mathbb{F}_{\bar{\mathbf{n}}}\cap\mathcal{P}$
has dimension $1$. Note that $\mathbb{F}_{\bar{\mathbf{n}}}\in\mathrm{SyM}_{3}^{tr=0}$
with $\mathbb{F}_{\bar{\mathbf{n}}}\neq[\mathbf{0}]$ implies $\dim\ker\mathbb{F}_{\bar{\mathbf{n}}}=1$,
and hence $\ker\mathbb{F}_{\bar{\mathbf{n}}}=\mathbb{R}\cdot\bar{\mathrm{Q}}$
with some unit vector $\bar{\mathrm{Q}}\in\mathcal{P}$ since . Let
$\bar{\mathrm{P}}$ be an unit vector in $\mathcal{P}$ perpendicular
to $\bar{\mathrm{Q}}$. Since $\mathbb{F}_{\bar{\mathbf{n}}}\big(\mathcal{P}\big)$
is nontrivial, $\mathbb{F}_{\bar{\mathbf{n}}}\bar{\mathrm{P}}=\kappa\cdot\bar{\mathbf{n}}$
for some $\kappa\neq0$. We show that $\mathbb{F}_{\bar{\mathbf{n}}}\bar{\mathbf{n}}=\kappa\cdot\bar{\mathrm{P}}$.
The triple $(\bar{\mathbf{n}},\bar{\mathrm{P}},\bar{\mathrm{Q}})$
constitutes an orthonormal basis, with respect to which the linear
map $\mathbb{F}_{\bar{\mathbf{n}}}$ is represented by some symmetric
matrix. Since $\mathbb{F}_{\bar{\mathbf{n}}}\bar{\mathrm{Q}}=\mathbf{0}$
and $\mathbb{F}_{\bar{\mathbf{n}}}\bar{\mathrm{P}}=\kappa\cdot\bar{\mathbf{n}}$,
it takes the form
\[
\left[\begin{array}{ccc}
a & \kappa & 0\\
\kappa & 0 & 0\\
0 & 0 & 0
\end{array}\right]
\]
with $a=\text{tr}\mathbb{F}_{\bar{\mathbf{n}}}=0$. As a result, $\mathrm{P}:=\mathbb{F}_{\bar{\mathbf{n}}}\bar{\mathbf{n}}=\kappa\cdot\bar{\mathrm{P}}\in\mathcal{P}$
and $||\mathrm{P}||=|\kappa|$, and hence $\mathbb{F}_{\bar{\mathbf{n}}}=\mathrm{P}\bar{\mathbf{n}}^{\intercal}+\bar{\mathbf{n}}\mathrm{P}^{\intercal}$. 

For any $\mathbf{M}\in\mathcal{P}$, 
\[
\mathbf{M}=||\mathbf{M}||\big(\cos\beta\cdot\bar{\mathrm{P}}+\sin\beta\cdot\bar{\mathrm{Q}}\big),
\]
and then the reciprocity of $\mathbb{F}$ implies $\mathbb{F}_{\mathbf{M}}\bar{\mathbf{n}}=(\mathrm{P}\cdot\mathbf{M})\bar{\mathbf{n}}$:
\[
\mathbb{F}_{\mathbf{M}}\bar{\mathbf{n}}=\cos\beta\cdot\mathbb{F}_{\bar{\mathbf{n}}}\bar{\mathrm{P}}+\sin\beta\cdot\mathbb{F}_{\bar{\mathbf{n}}}\bar{\mathrm{Q}}=||\mathbf{M}||\cos\beta\bar{\mathbf{n}}.
\]
\end{proof}
\begin{cor}
\label{col:||=00005CF_=00005Cn||=000026(=00005Cn,=00005CP,=00005CQ)}If
an unit vector $\bar{\mathbf{n}}\in\mathbb{R}^{3}$ satisfies $[\bar{\mathbf{n}}]\mathbb{F}_{\bar{\mathbf{n}}}[\bar{\mathbf{n}}]=[\mathbf{0}]$,
then $||\mathbb{F}_{\bar{\mathbf{n}}}||_{\mathbb{R}^{9}}^{2}=2||\mathrm{P}||^{2}$
with $\mathrm{P}:=\mathbb{F}_{\bar{\mathbf{n}}}\bar{\mathbf{n}}$.
Furthermore, if $\mathbb{F}_{\bar{\mathbf{n}}}$ is nontrivial, then
$\bar{\mathbf{n}},\bar{\mathrm{P}},\bar{\mathrm{Q}}$ constitute an
orthonormal basis, where $\bar{\mathrm{P}}=\frac{\mathrm{P}}{||\mathrm{P}||}$,
and, $\bar{\mathrm{Q}}=\bar{\mathbf{n}}\times\bar{\mathrm{P}}$ and
it is in $\ker\mathbb{F}_{\bar{\mathbf{n}}}$.
\end{cor}
\begin{proof}
These results follows directly from $\mathbb{F}_{\bar{\mathbf{n}}}=\mathrm{P}\bar{\mathbf{n}}^{\intercal}+\bar{\mathbf{n}}\mathrm{P}^{\intercal}$
and $\mathrm{P}\perp\bar{\mathbf{n}}$. 
\end{proof}
From Theorem \ref{thm:StruturalThm-PlaneConfig} we can prove that,
in this general context, $\bar{\mathbf{n}}$ is also an eigenvector
of $\mathbb{F}^{\intercal}\mathbb{F}$ with eigenvalue $2||\mathrm{P}||^{2}$.
For clarity, recall that as a linear map from $\mathbb{R}^{3}$ to
$\text{\ensuremath{\mathrm{SyM}_{3}^{tr=0}}}\subset\mathbb{R}^{9}$,
$\mathbb{F}$ has a canonical representation as a $9\times3$ matrix,
based on which the $3\times3$ matrix $\mathbb{F}^{\intercal}\mathbb{F}$
is constructed. Also, 
\[
\mathbf{m}^{\intercal}[\mathbb{F}^{\intercal}\mathbb{F}]\mathbf{M}=\text{tr}\mathbb{F}_{\mathbf{m}}\mathbb{F}_{\mathbf{M}}.
\]

\begin{thm}
\label{thm:=00005BStructuralThm-PlaneConfig=00005D=00005Cn-eigenVectorF*F}As
a normal vector of an $\mathbb{F}$-invariant plane, $\bar{\mathbf{n}}$
is an eigenvector of $\mathbb{F}^{\intercal}\mathbb{F}$ with eigenvalue
$2||\mathrm{P}||^{2}$.
\end{thm}
\begin{proof}
As usual, let $\bar{\mathbf{n}}$ be an unit vector. If it is indeed
an eigenvector, then following Theorem \ref{col:||=00005CF_=00005Cn||=000026(=00005Cn,=00005CP,=00005CQ)}
the eigenvalue would be
\[
\bar{\mathbf{n}}^{\intercal}\mathbb{F}^{\intercal}\mathbb{F}\bar{\mathbf{n}}=||\mathbb{F}_{\bar{\mathbf{n}}}||_{\mathbb{R}^{9}}^{2}=2||\mathrm{P}||^{2}.
\]

To prove that $\bar{\mathbf{n}}$ is an eigenvector of $\mathbb{F}^{\intercal}\mathbb{F}$,
it suffices to show for any $\mathbf{M}\perp\bar{\mathbf{n}}$, $\mathbf{M}^{\intercal}\mathbb{F}^{\intercal}\mathbb{F}\bar{\mathbf{n}}=0$.
Take the orthonormal basis $(\bar{\mathbf{n}},\bar{\mathrm{P}},\bar{\mathrm{Q}})$
with $\bar{\mathrm{P}}$ and $\bar{\mathrm{Q}}$ as in the proofs
of Theorems \ref{thm:StruturalThm-PlaneConfig} and \ref{col:||=00005CF_=00005Cn||=000026(=00005Cn,=00005CP,=00005CQ)}.
The matrix representations of $\mathbb{F}_{\mathbf{M}}$ and $\mathbb{F}_{\bar{\mathbf{n}}}$
w.r.t. $(\bar{\mathbf{n}},\bar{\mathrm{P}},\bar{\mathrm{Q}})$ then
take the form 
\[
\left[\begin{array}{ccc}
a & 0 & 0\\
0 & b & c\\
0 & c & d
\end{array}\right]\ \text{ and }\ \left[\begin{array}{ccc}
0 & \kappa & 0\\
\kappa & 0 & 0\\
0 & 0 & 0
\end{array}\right],\ \text{ respectively}.
\]
Direct computation shows $\mathbf{M}^{\intercal}\mathbb{F}^{\intercal}\mathbb{F}\bar{\mathbf{n}}=\text{tr}\mathbb{F}_{\mathbf{M}}\mathbb{F}_{\bar{\mathbf{n}}}=0$.
\end{proof}
From Theorem \ref{thm:StruturalThm-PlaneConfig} we know that $\bar{\mathbf{n}}$
is an eigenvector of $\mathbb{F}^{\intercal}\mathbb{F}$ with eigenvalue
$2||\mathrm{P}||^{2}$.  Therefore, as the largest eigenvalue of $\mathbb{F}^{\intercal}\mathbb{F}$,
\[
\lambda_{\mathbb{F}}=||\mathbb{F}_{\mathbf{M}_{\mathbb{F}}}||_{\mathbb{R}^{9}}^{2}\geq2||\mathrm{P}||^{2},
\]
combining which with (\ref{eq:(1st estimate)<=00005Cbar=00007B=00005Clambda=00007D<})
yields
\begin{equation}
\lambda_{\mathbf{M}_{\mathbb{F}}}^{2}\geq\frac{||\mathbb{F}_{\mathbf{M}_{\mathbb{F}}}||_{\mathbb{R}^{9}}^{2}}{2}\geq||\mathrm{P}||^{2}\ \big(=\lambda_{\bar{\mathbf{n}}}^{2}\big).\label{Ineq:LowBound-=00005CLambda_=00005CMF}
\end{equation}

Due to the symmetry, $\mathcal{P}$ is also an invariant subspace
of $\mathbb{F}^{\intercal}\mathbb{F}$. As an eigenvector of $\mathbb{F}^{\intercal}\mathbb{F}$
associated to its largest eigenvalue $\lambda_{\mathbb{F}}$, $\mathbf{M}_{\mathbb{F}}$
can always be taken from the set $\{\pm\bar{\mathbf{n}}\}\cup\mathcal{P}$.
(Of course, when $\lambda_{\mathbb{F}}=2||\mathrm{P}||^{2}$ and its
multiplicity is $2$ or more, $\mathbf{M}_{\mathbb{F}}$ can also
be taken from other place.)

\subsection{$\mathcal{P}$-Decompositions}

For an actual magnetic system we can see from (\ref{eq:F(M)}) that
each matrix $\mathbb{F}_{\mathbf{M}}$ breaks into two parts: $\mathbb{F}_{\mathbf{M}}=\mathbb{E}_{\mathbf{M}}-\mathbb{P}_{\mathbf{M}}$
with
\begin{equation}
\mathbb{E}_{\mathbf{M}}=\mathbf{M}\mathrm{P}^{\intercal}+\text{\ensuremath{\mathrm{P}}}\mathbf{M}^{\intercal}+(\mathbf{M}\cdot\mathrm{P})\mathbf{Id}\label{eq:E-=00005BStandard=00005D}
\end{equation}
and
\[
\mathbb{P}_{\mathbf{M}}=5\sum_{i=1}^{n}\frac{\big(\mathbf{M}^{\intercal}\hat{\mathrm{p}}_{i}\big)\cdot\hat{\mathrm{p}}_{i}\hat{\mathrm{p}}_{i}^{\intercal}}{||\mathrm{p}_{i}||^{4}}.
\]
While $\mathbb{P}_{\mathbf{M}}$ maps the whole $\mathbb{R}^{3}$
to the invariant plane $\mathcal{P}$, $\mathbb{E}:\text{\ensuremath{\mathbf{M}}}\mapsto\mathbb{E}_{\mathbf{M}}$
is equivariant under the circle action by $\mathrm{Rot}_{\mathrm{P}}$,
the subgroup of $\mathbb{SO}(3)$ that keeps $\mathrm{P}$ still:
\begin{equation}
\mathbb{E}_{\mathrm{C}\mathbf{M}}=\mathrm{C}\mathbb{E}_{\mathbf{M}}\mathrm{C}^{\intercal},\;\forall\mathrm{C}\in\mathrm{Rot}_{\mathrm{P}}.\label{eq:Rot_=00005CP-Equivariance}
\end{equation}

In general, given any (abstract) magnetic algebra $(\mathbb{R}^{3},\mathbb{F})$
with invariant plane $\mathcal{P}$, we have specified with Definition
\ref{def:PolarDecomposition} this type of decompositions (if they
exist) as $\mathcal{P}$-decompositions. As is stated in Theorem \ref{thm:=00005B=00005CP-decomposition=00005D},
every magnetic algebra with planarity admits a family of such decompositions.
We first prove Theorem \ref{thm:=00005CF=00003D=00005CE-=00005CP}
(Theorem \ref{thm:=00005Bspecific=00005D=00005CF=00003D=00005CE-=00005CP})
to confirm the existence of a (particular) $\mathcal{P}$-decomposition, based on which
we then complete the whole proof for Theorem \ref{thm:=00005B=00005CP-decomposition=00005D}.
\begin{thm}
\label{thm:=00005CF=00003D=00005CE-=00005CP}Let $(\mathbb{R}^{3},\mathbb{F})$
be a magnetic algebra with an $\mathbb{F}$-invariant plane $\mathcal{P}$.
Suppose that for $\bar{\mathbf{n}}\perp\mathcal{P}$ with $||\bar{\mathbf{n}}||=1$, and, 
$\mathrm{P}=\mathbb{F}_{\bar{\mathbf{n}}}\bar{\mathbf{n}}\neq\mathbf{0}$.
Then, there exists a $\mathcal{P}$-decomposition
\[
\mathbb{F}=\mathbb{E}-\mathbb{P}
\]
with linear maps $\mathbb{E},\mathbb{P}:\mathbb{R}^{3}\rightarrow\mathrm{SyM}_{3}$
such that for each $\mathbf{M}\in\mathbb{R}^{3}$, $\mathrm{Im}\mathbb{P}_{\mathbf{M}}\subset\mathcal{P}$
and 
\[
\mathbb{E}=\mathrm{C}\mathbb{E}_{\mathbf{M}}\mathrm{C}^{\intercal},\ \text{\ensuremath{\forall}\ensuremath{\mathrm{C}\in} \ensuremath{\mathrm{Rot}_{\mathrm{P}}}}.
\]
\end{thm}
\begin{proof}
As before, $\bar{\mathbf{n}}$ is an unit normal vector to $\mathcal{P}$.
We first show the existence of the decomposition. Each $\mathbf{M}$
and $\mathbf{m}$ in $\mathbb{R}^{3}$ can be decomposed as
\[
\mathbf{M}=(\mathbf{M}\cdot\bar{\mathbf{n}})\bar{\mathbf{n}}+\mathbf{M}_{\bot}\ \text{ and }\ \mathbf{m}=(\mathbf{m}\cdot\bar{\mathbf{n}})\bar{\mathbf{n}}+\mathbf{m}_{\bot}
\]
with $\mathbf{M}_{\bot},\mathbf{m}_{\bot}\in\mathcal{P}$ being the
orthogonal projection of $\mathbf{M}$,$\mathbf{m}$ onto $\mathcal{P}$,
respectively. As a result, 
\[
\begin{aligned}\mathbb{F}_{\mathbf{M}}\mathbf{m} &=(\mathbf{M}\cdot\bar{\mathbf{n}})\mathbb{F}_{\bar{\mathbf{n}}}\mathbf{m}+\mathbb{F}_{\mathbf{M}_{\bot}}\mathbf{m}\\
 &=(\mathbf{M}\cdot\bar{\mathbf{n}})\mathbb{F}_{\bar{\mathbf{n}}}\mathbf{m}+(\mathbf{m}\cdot\bar{\mathbf{n}})\mathbb{F}_{\mathbf{M}_{\bot}}\bar{\mathbf{n}}+\mathbb{F}_{\mathbf{M}_{\bot}}\mathbf{m}_{\bot}\\
 &=(\mathbf{M}\cdot\bar{\mathbf{n}})\big[\mathrm{P}\bar{\mathbf{n}}^{\intercal}+\bar{\mathbf{n}}\mathrm{P}^{\intercal}\big]\mathbf{m}+(\mathbf{m}\cdot\bar{\mathbf{n}})\big[\mathrm{P}\bar{\mathbf{n}}^{\intercal}+\bar{\mathbf{n}}\mathrm{P}^{\intercal}\big]\mathbf{M}_{\bot}+\mathbb{F}_{\mathbf{M}_{\bot}}\mathbf{m}_{\bot}\\
 &=\big[\mathrm{P}\mathbf{M}_{\bar{\mathbf{n}}}^{\intercal}+\mathbf{M}_{\bar{\mathbf{n}}}\mathrm{P}^{\intercal}\big]\mathbf{m}+\big(\mathbf{M}\cdot\mathrm{P}\big)\bar{\mathbf{n}}\bar{\mathbf{n}}^{\intercal}\mathbf{m}+\mathbb{F}_{\mathbf{M}_{\bot}}\mathbf{m}_{\bot}\\
 &=\big[\mathrm{P}\mathbf{M}^{\intercal}+\mathbf{M}\mathrm{P}^{\intercal}\big]\mathbf{m}+\big(\mathbf{M}\cdot\mathrm{P}\big)\mathbf{m}+\mathbb{F}_{\mathbf{M}_{\bot}}\mathbf{m}_{\bot}-\mathbb{E}_{\mathbf{M}}^{\perp}\mathbf{m},
\end{aligned}
\]
in which $\mathbf{M}_{\bar{\mathbf{n}}}=(\mathbf{M}\cdot\bar{\mathbf{n}})\bar{\mathbf{n}}$,
and, with $\big[\mathbf{pr}_{\perp}\big]=\mathbf{Id}-\bar{\mathbf{n}}\bar{\mathbf{n}}^{\intercal}$,
\[
\mathbb{E}_{\mathbf{M}}^{\perp}:=\big[\mathrm{P}\mathbf{M}_{\perp}^{\intercal}+\mathbf{M}_{\perp}\mathrm{P}^{\intercal}\big]+\big(\mathbf{M}\cdot\mathrm{P}\big)\big[\mathbf{pr}_{\perp}\big].
\]
Define $\mathbb{P}_{\mathbf{M}}:=\mathbb{E}_{\mathbf{M}}^{\perp}-\mathbb{F}_{\mathbf{M}_{\perp}}\big[\mathbf{pr}_{\perp}\big]$
and check that both $\mathbb{E}_{\mathbf{M}}^{\perp}$ and $\mathbb{F}_{\mathbf{M}_{\perp}}\big[\mathbf{pr}_{\perp}\big]$
have their images in $\mathcal{P}$. Let $\mathbb{E}_{\mathbf{M}}=\mathbb{F}_{\mathbf{M}}+\mathbb{P}_{\mathbf{M}}$
and then
\[
\mathbb{E}_{\mathbf{M}}=\big[\mathrm{P}\mathbf{M}^{\intercal}+\mathbf{M}\mathrm{P}^{\intercal}\big]+\big(\mathbf{M}\cdot\mathrm{P}\big)\mathbf{Id}.
\]
It is straightforward to check that $\mathbb{E}$ satisfies (\ref{eq:Rot_=00005CP-Equivariance}).
\end{proof}
Note that in the proof above the $\mathbb{E}_{\mathbf{M}}$ has been
constructed as in (\ref{eq:E-=00005BStandard=00005D}), and hence
what is proved therein is actually something more specific:
\begin{thm}
\label{thm:=00005Bspecific=00005D=00005CF=00003D=00005CE-=00005CP}Let
$(\mathbb{R}^{3},\mathbb{F})$ be a magnetic algebra with an $\mathbb{F}$-invariant
plane $\mathcal{P}$. For each $\mathbf{M}\in\mathbb{R}^{3}$, with
$\mathbb{E}_{\mathbf{M}}$ defined in (\ref{eq:E-=00005BStandard=00005D}),
$\mathbb{P}_{\mathbf{M}}:=\mathbb{E}_{\mathbf{M}}-\mathbb{F}_{\mathbf{M}}$
maps $\mathbb{R}^{3}$ to $\mathcal{P}$.
\end{thm}
Based on Theorem \ref{thm:=00005Bspecific=00005D=00005CF=00003D=00005CE-=00005CP},
the proof for Theorem \ref{thm:=00005B=00005CP-decomposition=00005D}
will be completed with the help of the lemmas given below.
The key of the theorem is the structure of the operator $\mathbb{E}^{\xi}$.
Due to linearity, for each $\mathbf{M}\in\mathbb{R}^{3}$ with the
orthogonal decomposition $\mathbf{M}=(\mathbf{M}\cdot\bar{\mathrm{P}})\bar{\mathrm{P}}+\mathbf{M}_{\bar{\mathrm{P}}}^{\perp}$,
$\mathbb{E}_{\mathbf{M}}^{\xi}=(\mathbf{M}\cdot\bar{\mathrm{P}})\mathbb{E}_{\bar{\mathrm{P}}}^{\xi}+\mathbb{E}_{\mathbf{M}_{\bar{\mathrm{P}}}^{\perp}}^{\xi}$,
which suggests to look into $\mathbb{E}_{\bar{\mathrm{P}}}^{\xi}$
and $\mathbb{E}_{\mathbf{m}}^{\xi}$ with $\mathbf{m}\perp\bar{\mathrm{P}}$.
Moreover, since $\mathbb{E}_{\mathbf{m}}^{\xi}=\mathbb{E}_{\mathrm{C}\bar{\mathbf{n}}}^{\xi}=\mathrm{C}\mathbb{E}_{\bar{\mathbf{n}}}^{\xi}\mathrm{C}^{\intercal}$
for some $\mathrm{C}\in\mathrm{Rot}_{\mathrm{P}}$, it suffices to
look at $\mathbb{E}_{\bar{\mathbf{n}}}^{\xi}$. Lemma \ref{lem:Structure-=00005CE}
and Lemma \ref{lem:=00005CE_=00005Cn} reveal basic structures of
$\mathbb{E}_{\bar{\mathrm{P}}}^{\xi}$ and $\mathbb{E}_{\bar{\mathbf{n}}}^{\xi}$,
respectively.
\begin{lem}
\label{lem:Structure-=00005CE}Let $\mathbb{E}:\mathbb{R}^{3}\rightarrow\mathrm{SyM}_{3}$
be a linear map. Suppose that $\exists\bar{\mathrm{P}}\neq\mathbf{0}$
such that $\mathbb{E}_{\mathrm{C}\mathbf{M}}=\mathrm{C}\mathbb{E}_{\mathbf{M}}\mathrm{C}^{\intercal}$
for any $\mathrm{C}\in\mathrm{Rot}_{\bar{\mathrm{P}}}$ and $\mathbf{M}\in\mathbb{R}^3$. Then, there exist $\xi,\zeta\in\mathbb{R}$ such that
\begin{equation}
\mathbb{E}_{\bar{\mathrm{P}}}=\xi\bar{\mathrm{P}}\bar{\mathrm{P}}^{\intercal}+\zeta\big(\mathbf{Id}-\bar{\mathrm{P}}\bar{\mathrm{P}}^{\intercal}\big),\label{eq:=00005CE_=00005CP-general}
\end{equation}
and, for any $\mathbf{M}\perp\bar{\mathrm{P}}$, it holds $\text{tr}\mathbb{E}_{\mathbf{M}}=0$
and $\text{tr}\mathbb{E}_{\bar{\mathrm{P}}}\mathbb{E}_{\mathbf{M}}=0$. 
\end{lem}
\begin{proof}
For proving the identity for $\mathbb{E}_{\bar{\mathrm{P}}}$, it
suffices to show that, $\bar{\mathrm{P}}$ is an eigenvector of $\mathbb{E}_{\bar{\mathrm{P}}}$
(with eigenvalue $\xi$), and, the other two eigenvalues are the same
($=\zeta$). 

Since $\mathrm{C}\bar{\mathrm{P}}=\bar{\mathrm{P}}$ for all $\mathrm{C}\in\mathrm{Rot}_{\bar{\mathrm{P}}}$,
we have
\[
\mathbb{E}_{\bar{\mathrm{P}}}\bar{\mathrm{P}}=\mathbb{E}_{\mathrm{C}\bar{\mathrm{P}}}\bar{\mathrm{P}}=\mathrm{C}\mathbb{E}_{\bar{\mathrm{P}}}\mathrm{C}^{\intercal}\bar{\mathrm{P}}=\mathrm{C}\mathbb{E}_{\bar{\mathrm{P}}}\bar{\mathrm{P}}.
\]
That is, the vector $\mathbb{E}_{\bar{\mathrm{P}}}\bar{\mathrm{P}}$
is invariant under any rotation around $\bar{\mathrm{P}}$, which
means $\mathbb{E}_{\bar{\mathrm{P}}}\bar{\mathrm{P}}=\xi\bar{\mathrm{P}}$
for some $\xi\in\mathbb{R}^{3}$. Consequently, there are at least
other two eigenvectors, $\mathbf{e}_{0}$ and $\mathbf{e}_{1}$, of
$\mathbb{E}_{\bar{\mathrm{P}}}$ such that $\mathbf{e}_{0}\perp\mathbf{e}_{1}$
and $\mathbf{e}_{0},\mathbf{e}_{1}\perp\bar{\mathrm{P}}$. Note that
for any $\mathbf{m}\perp\bar{\mathrm{P}}$ and $\mathrm{C}\in\mathrm{Rot}_{\bar{\mathrm{P}}}$,
it holds
\[
\mathbf{m}^{\intercal}\mathbb{E}_{\bar{\mathrm{P}}}\mathbf{m}=\mathbf{m}^{\intercal}\mathbb{E}_{\mathrm{C}^{\intercal}\bar{\mathrm{P}}}\mathbf{m}=\mathbf{m}^{\intercal}\mathrm{C}^{\intercal}\mathbb{E}_{\bar{\mathrm{P}}}\mathrm{C}\mathbf{m},
\]
and hence $\mathbb{E}_{\bar{\mathrm{P}}}\mathbf{e}_{0}=\zeta\mathbf{e}_{0}$
and $\mathbb{E}_{\bar{\mathrm{P}}}\mathbf{e}_{1}=\zeta\mathbf{e}_{1}$
for some $\zeta\in\mathbb{R}$. As a result, (\ref{eq:=00005CE_=00005CP-general})
holds.

Note that for any $\mathbf{M}\in\mathbb{R}^{3}$ and $\mathrm{C}\in\mathrm{Rot}_{\bar{\mathrm{P}}}$,
$\text{tr}\mathbb{E}_{\mathrm{C}\mathbf{M}}=\text{tr}\mathrm{C}\mathbb{E}_{\mathbf{M}}\mathrm{C}^{\intercal}=\text{tr}\mathbb{E}_{\mathbf{M}}$
and
\[
\text{tr}\mathbb{E}_{\bar{\mathrm{P}}}\mathbb{E}_{\mathrm{C}\mathbf{M}}=\text{tr}\mathbb{E}_{\mathrm{C}\bar{\mathrm{P}}}\mathbb{E}_{\mathrm{C}\mathbf{M}}=\text{tr}\mathrm{C}\mathbb{E}_{\bar{\mathrm{P}}}\mathbb{E}_{\mathbf{M}}\mathrm{C}^{\intercal}=\text{tr}\mathbb{E}_{\bar{\mathrm{P}}}\mathbb{E}_{\mathbf{M}}
\]
and hence the linear functions $\mathbf{M}\mapsto\text{tr}\mathbb{E}_{\mathbf{M}}$
and $\mathbf{M}\mapsto\text{tr}\mathbb{E}_{\bar{\mathrm{P}}}\mathbb{E}_{\mathbf{M}}$
are invariant under any rotation around $\bar{\mathrm{P}}$. Therefore,
their Riez representations are vectors parallel to $\bar{\mathrm{P}}$,
i.e.,
\[
\exists a,b\in\mathbb{R},\ \text{s.t.}\ \text{tr}\mathbb{E}_{\mathbf{M}}=a\bar{\mathrm{P}}\cdot\mathbf{M}\ \text{ and }\ \text{tr}\mathbb{E}_{\bar{\mathrm{P}}}\mathbb{E}_{\mathbf{M}}=b\bar{\mathrm{P}}\cdot\mathbf{M}.
\]
As a result, $\text{tr}\mathbb{E}_{\mathbf{M}}=0$ and $\text{tr}\mathbb{E}_{\bar{\mathrm{P}}}\mathbb{E}_{\mathbf{M}}=0$
if $\mathbf{M}\perp\bar{\mathrm{P}}$.
\end{proof}
For $\mathbb{E}_{\mathbf{M}}$ with $\mathbf{M}\in\mathcal{P}$, note
that due to the equivairant $\mathbb{E}_{\mathrm{C}\mathbf{M}}=\mathrm{C}\mathbb{E}_{\mathbf{M}}\mathrm{C}^{\intercal}$,
it suffices to look at $\mathbb{E}_{\bar{\mathbf{n}}}$. We first
demonstrate a result obtained with a weaker condition: instead of
requiring $\mathrm{Im}\mathbb{P}_{\mathbf{M}}\subset\mathcal{P}$
for every $\mathbf{M}\in\mathbb{R}^{3}$, only $\mathrm{Im}\mathbb{P}_{\bar{\mathbf{n}}}\subset\mathcal{P}$
is assumed:
\begin{lem}
\label{lem:=00005CE_=00005Cn} With the assumptions about $\mathbb{F}$
in Theorem \ref{thm:=00005CF=00003D=00005CE-=00005CP}, if $\mathbb{F}=\mathbb{E}-\mathbb{P}$
with the linear maps $\mathbb{E},\mathbb{P}:\mathbb{R}^{3}\rightarrow\mathrm{SyM}_{3}$
satisfying $\mathbb{E}_{\mathrm{C}\mathbf{M}}=\mathrm{C}\mathbb{E}_{\mathbf{M}}\mathrm{C}^{\intercal}$
for all $\mathrm{C}\in\mathrm{Rot}_{\bar{\mathrm{P}}}$, and, $\mathrm{Im}\mathbb{P}_{\bar{\mathbf{n}}}\subset\mathcal{P}$.
Then, $\exists\eta\in\mathbb{R}$
such that: 
\begin{equation}
\mathbb{E}_{\bar{\mathbf{n}}}=\big[\mathrm{P}\bar{\mathbf{n}}^{\intercal}+\bar{\mathbf{n}}\mathrm{P}^{\intercal}\big]+\eta\cdot[\mathrm{P}\bar{\mathrm{Q}}^{\intercal}+\bar{\mathrm{Q}}\mathrm{P}^{\intercal}],\label{eq:=00005CE_=00005Cn-WeakCondition}
\end{equation}
where $\bar{\mathrm{Q}}:=\bar{\mathbf{n}}\times\frac{\mathrm{P}}{||\mathrm{P}||}$. 
\end{lem}
\begin{proof}
From Corollary \ref{col:||=00005CF_=00005Cn||=000026(=00005Cn,=00005CP,=00005CQ)}
we know that $\bar{\mathrm{Q}}\in\ker\mathbb{F}_{\bar{\mathbf{n}}}$,
and also, $(\bar{\mathbf{n}},\bar{\mathrm{P}},\bar{\mathrm{Q}})$
is an orthonormal basis of $\mathbb{R}^{3}$. We prove the theorem
by looking at the matrix representation of $\mathbb{E}_{\bar{\mathbf{n}}}$
w.r.t. $(\bar{\mathbf{n}},\bar{\mathrm{P}},\bar{\mathrm{Q}})$. 

Since $\mathrm{Im}\mathbb{P}_{\bar{\mathbf{n}}}\subset\mathcal{P}$,
it holds $\bar{\mathbf{n}}\mathbb{F}_{\bar{\mathbf{n}}}\mathbf{m}=\bar{\mathbf{n}}\mathbb{E}_{\bar{\mathbf{n}}}\mathbf{m}$
for any $\mathbf{m}\in\mathbb{R}^{3}$. As a result, 
\[
\bar{\mathbf{n}}\cdot\mathbb{E}_{\bar{\mathbf{n}}}\bar{\mathbf{n}}=\bar{\mathbf{n}}\cdot\mathbb{F}_{\bar{\mathbf{n}}}\bar{\mathbf{n}}=\bar{\mathbf{n}}\cdot\mathrm{P}=0,
\]
\[
\bar{\mathbf{n}}\cdot\mathbb{E}_{\bar{\mathbf{n}}}\bar{\mathrm{Q}}=\bar{\mathbf{n}}\cdot\mathbb{F}_{\bar{\mathbf{n}}}\bar{\mathrm{Q}}=0,
\]
and $\bar{\mathbf{n}}\cdot\mathbb{E}_{\bar{\mathbf{n}}}\bar{\mathrm{P}}=\bar{\mathbf{n}}\cdot\mathbb{F}_{\bar{\mathbf{n}}}\bar{\mathrm{P}}=||\mathrm{P}||$. 

To compute $\bar{\mathrm{P}}^{\intercal}\mathbb{E}_{\bar{\mathbf{n}}}\bar{\mathrm{P}}$,
consider the linear function on $\mathbb{R}^{3}:$ $f(\mathbf{M})=\bar{\mathrm{P}}^{\intercal}\mathbb{E}_{\mathbf{M}}\bar{\mathrm{P}}$
and note that
\[
\bar{\mathrm{P}}^{\intercal}\mathbb{E}_{\mathrm{C}\mathbf{M}}\bar{\mathrm{P}}=\bar{\mathrm{P}}^{\intercal}\mathrm{C}\mathbb{E}_{\mathbf{M}}\mathrm{C}^{\intercal}\bar{\mathrm{P}}=\bar{\mathrm{P}}^{\intercal}\mathbb{E}_{\mathbf{M}}\bar{\mathrm{P}},\ \forall\mathrm{C}\in\mathrm{Rot}_{\mathrm{P}}.
\]
This means that the Riez representation of $f$ is invariant under
any rotation around $\mathrm{P}$, and hence it is parallel to $\mathrm{P}$,
i.e., $\exists r\in\mathbb{R}$ s.t. $f(\mathbf{M})=r\mathrm{P}\cdot\mathbf{M}$.
As a result, 
\[
\bar{\mathrm{P}}^{\intercal}\mathbb{E}_{\bar{\mathbf{n}}}\bar{\mathrm{P}}=f(\bar{\mathbf{n}})=0.
\]
According to Lemma \ref{lem:Structure-=00005CE}, $\text{tr}\mathbb{E}_{\bar{\mathbf{n}}}=0$,
and hence
\[
\bar{\mathrm{Q}}^{\intercal}\mathbb{E}_{\bar{\mathbf{n}}}\bar{\mathrm{Q}}=-\bar{\mathbf{n}}\cdot\mathbb{E}_{\bar{\mathbf{n}}}\bar{\mathbf{n}}-\bar{\mathrm{P}}^{\intercal}\mathbb{E}_{\bar{\mathbf{n}}}\bar{\mathrm{P}}=0.
\]
To this point, we already determined the matrix of $\mathbb{E}_{\bar{\mathbf{n}}}$
up to the form
\[
\big[\mathbb{E}_{\bar{\mathbf{n}}}\big]=\left[\begin{array}{ccc}
0 & ||\mathrm{P}|| & 0\\
||\mathrm{P}|| & 0 & \eta\\
0 & \eta & 0
\end{array}\right]
\]
or equivalently, $\mathbb{E}_{\bar{\mathbf{n}}}=\big[\mathrm{P}\bar{\mathbf{n}}^{\intercal}+\bar{\mathbf{n}}\mathrm{P}^{\intercal}\big]+\eta\cdot[\mathrm{P}\bar{\mathrm{Q}}^{\intercal}+\bar{\mathrm{Q}}\mathrm{P}^{\intercal}]$
(abusing the notation $\eta$). 
\end{proof}
Now we are ready to prove Theorem \ref{thm:=00005B=00005CP-decomposition=00005D}.
\begin{proof}
\textbf{{[}of Theorem \ref{thm:=00005B=00005CP-decomposition=00005D}{]}}. Theorem \ref{thm:=00005Bspecific=00005D=00005CF=00003D=00005CE-=00005CP} has confirmed that there exists at least one $\mathcal{P}$-decomposition of $\mathbb{F}$, which is the one determined by \eqref{eq:=00005Bparameterization=00005D-=00005CE} with $\xi=0$. This actually implies that any decomposition of $\mathbb{F}$ determined by \eqref{eq:=00005Bparameterization=00005D-=00005CE} is a $\mathcal{P}$-decomposition. It remains to show that, conversely, each $\mathcal{P}$-decomposition of $\mathbb{F}$ takes the form of $\mathbb{F}=\mathbb{E}^\xi-\mathbb{P}^\xi$ with $\mathbb{E}^\xi$ given by Eq. \eqref{eq:=00005Bparameterization=00005D-=00005CE}.

Let $\mathbb{F}=\mathbb{E}-\mathbb{P}$ be an arbitrary $\mathcal{P}$-decomposition. Note that by definition this implies $\mathrm{Im}\mathbb{P}_{\mathbf{M}}\subset\mathcal{P}$ for any $\mathbf{M}\in\mathbb{R}^3$, and in particular, $\mathrm{Im}\mathbb{P}_{\bar{\mathrm{Q}}},\mathrm{Im}\mathbb{P}_{\bar{\mathrm{P}}}\subset\mathcal{P}$. We first show that
the coefficient $\eta$ in Lemma \ref{lem:=00005CE_=00005Cn} should
be $0$, and then for any $\mathbf{M}=\mathrm{C}\bar{\mathbf{n}}$
with $\mathrm{C}\in\mathrm{Rot}_{\mathrm{P}}$, it holds
\[
\mathbb{E}_{\mathbf{M}}=\mathrm{P}\mathbf{M}^{\intercal}+\mathbf{M}\mathrm{P}^{\intercal}.
\]
Note that $\bar{\mathbf{n}}$ and $\bar{\mathrm{Q}}$ are related
by a rotation around $\mathrm{P}$ with angle $\frac{\pi}{2}$, that
is, $\bar{\mathrm{Q}}=\mathrm{C}_{\frac{\pi}{2}}\bar{\mathbf{n}}$
for some $\mathrm{C}_{\frac{\pi}{2}}\in\mathrm{Rot}_{\mathrm{P}}$
with $\mathrm{C}_{\frac{\pi}{2}}^{2}=-\mathbf{Id}$. Therefore, $\mathbb{E}_{\bar{\mathrm{Q}}}=\mathbb{E}_{\mathrm{C}_{\frac{\pi}{2}}\bar{\mathbf{n}}}=\mathrm{C}_{\frac{\pi}{2}}\mathbb{E}_{\bar{\mathbf{n}}}\mathrm{C}_{\frac{\pi}{2}}^{\intercal}$,
combining which with (\ref{eq:=00005CE_=00005Cn-WeakCondition}) gives
\[
\mathbb{E}_{\bar{\mathrm{Q}}}=\big[\mathrm{P}\bar{\mathrm{Q}}^{\intercal}+\bar{\mathrm{Q}}\mathrm{P}^{\intercal}\big]-\eta\cdot[\mathrm{P}\bar{\mathbf{n}}^{\intercal}+\bar{\mathbf{n}}\mathrm{P}^{\intercal}],
\]
and then
\[
\mathbb{P}_{\bar{\mathrm{Q}}}=\bigg(\big[\mathrm{P}\bar{\mathrm{Q}}^{\intercal}+\bar{\mathrm{Q}}\mathrm{P}^{\intercal}\big]-\mathbb{F}_{\bar{\mathrm{Q}}}\bigg)-\eta\cdot[\mathrm{P}\bar{\mathbf{n}}^{\intercal}+\bar{\mathbf{n}}\mathrm{P}^{\intercal}].
\]

From Theorem \ref{thm:=00005Bspecific=00005D=00005CF=00003D=00005CE-=00005CP}
we know that $\big[\mathrm{P}\bar{\mathrm{Q}}^{\intercal}+\bar{\mathrm{Q}}\mathrm{P}^{\intercal}\big]-\mathbb{F}_{\bar{\mathrm{Q}}}$
maps $\mathbb{R}^{3}$ to the invariant plane $\mathcal{P}$, and
then
\[
\mathbb{P}_{\bar{\mathrm{Q}}}\mathrm{P}=\text{ a vector in }\mathcal{P}\ +\ \eta\cdot||\mathrm{P}||^{2}\bar{\mathbf{n}}.
\]
Since $\mathrm{Im}\mathbb{P}_{\bar{\mathrm{Q}}}\subset\mathcal{P}$, this implies $\eta=0$. Secondly, we check that, the
$\zeta$ in \eqref{eq:=00005CE_=00005CP-general} should be $||\mathrm{P}||$:
\[
\zeta=\bar{\mathbf{n}}^{\intercal}\mathbb{E}_{\bar{\mathrm{P}}}\bar{\mathbf{n}}=\bar{\mathbf{n}}^{\intercal}\mathbb{F}_{\bar{\mathrm{P}}}\bar{\mathbf{n}}=\bar{\mathrm{P}}\cdot\mathbb{F}_{\bar{\mathbf{n}}}\bar{\mathbf{n}}=||\mathrm{P}||.
\]
Combining both these results with the linearity of $\mathbb{E}$,
we see that, for arbitrary $\mathbf{M}\in\mathbb{R}^{3}$, using the orthogonal decomposition $\mathbf{M}=\mathrm{(\mathbf{M}\cdot\bar{\mathrm{P}})}\bar{\mathrm{P}}+\mathbf{M}_{\bar{\mathrm{P}}}^{\perp}$, we get
\[
\begin{aligned}\mathbb{E}_{\mathbf{M}} =&(\mathbf{M}\cdot\bar{\mathrm{P}})\mathbb{E}_{\bar{\mathrm{P}}}+\mathbb{E}_{\mathbf{M}_{\bar{\mathrm{P}}}^{\perp}}\\
 =&\mathrm{(\mathbf{M}\cdot\bar{\mathrm{P}})}\bigg[\xi\bar{\mathrm{P}}\bar{\mathrm{P}}^{\intercal}+||\mathrm{P}||\big(\mathbf{Id}-\bar{\mathrm{P}}\bar{\mathrm{P}}^{\intercal}\big)\bigg]+\mathrm{P}\big(\mathbf{M}_{\bar{\mathrm{P}}}^{\perp}\big)^{\intercal}+\mathbf{M}_{\bar{\mathrm{P}}}^{\perp}\mathrm{P}^{\intercal}\\
 =&\mathrm{(\mathbf{M}\cdot\mathrm{P})\mathbf{Id}}+\mathrm{2(\mathbf{M}\cdot\bar{\mathrm{P}})}||\mathrm{P}||\bar{\mathrm{P}}\bar{\mathrm{P}}^{\intercal}+\mathrm{P}\big(\mathbf{M}_{\bar{\mathrm{P}}}^{\perp}\big)^{\intercal}+\mathbf{M}_{\bar{\mathrm{P}}}^{\perp}\mathrm{P}^{\intercal}-\\
 &2\mathrm{(\mathbf{M}\cdot\bar{\mathrm{P}})}||\mathrm{P}||\cdot\bar{\mathrm{P}}\bar{\mathrm{P}}^{\intercal}+\mathrm{(\mathbf{M}\cdot\bar{\mathrm{P}})}(\xi-||\mathrm{P}||)\bar{\mathrm{P}}\bar{\mathrm{P}}^{\intercal}.
\end{aligned}
\]
Check that
\[
\begin{aligned} & \mathrm{2(\mathbf{M}\cdot\bar{\mathrm{P}})}||\mathrm{P}||\cdot\bar{\mathrm{P}}\bar{\mathrm{P}}^{\intercal}+\mathrm{P}\big(\mathbf{M}_{\bar{\mathrm{P}}}^{\perp}\big)^{\intercal}+\mathbf{M}_{\bar{\mathrm{P}}}^{\perp}\mathrm{P}^{\intercal}\\
= & \mathrm{P}\big[(\mathbf{M}\cdot\bar{\mathrm{P}})\bar{\mathrm{P}}\big]^{\intercal}+\mathrm{P}\big(\mathbf{M}_{\bar{\mathrm{P}}}^{\perp}\big)^{\intercal}+\mathrm{\big[(\mathbf{M}\cdot\bar{\mathrm{P}})}\bar{\mathrm{P}}\big]\mathrm{P}^{\intercal}+\mathbf{M}_{\bar{\mathrm{P}}}^{\perp}\mathrm{P}^{\intercal}\\
= & \mathrm{P}\mathbf{M}^{\intercal}+\mathrm{\mathbf{M}}\mathrm{P}^{\intercal},
\end{aligned}
\]
and also,
\[
2\mathrm{(\mathbf{M}\cdot\bar{\mathrm{P}})}||\mathrm{P}||\cdot\bar{\mathrm{P}}\bar{\mathrm{P}}^{\intercal}-\mathrm{(\mathbf{M}\cdot\bar{\mathrm{P}})}(\xi-||\mathrm{P}||)\bar{\mathrm{P}}\bar{\mathrm{P}}^{\intercal}=\frac{\big(3||\mathrm{P}||-\xi\big)}{||\mathrm{P}||^3}(\mathbf{M}\cdot\mathrm{P})\mathrm{P}\mathrm{P}^\intercal.
\]
The proof for \eqref{eq:=00005Bparameterization=00005D-=00005CE} is then completed by simply denoting the constant $\frac{\big(3||\mathrm{P}||-\xi\big)}{||\mathrm{P}||^3}$ by $\xi$.

To prove \eqref{eq:P_n=[0]}, it suffices to show $\mathbb{P}^{\xi}_{\bar{\mathbf{n}}}=[\mathbf{0}]$. Combining \eqref{eq:=00005Bparameterization=00005D-=00005CE} with Theorem \ref{thm:StruturalThm-PlaneConfig} yields $\mathbb{F}_{\bar{\mathbf{n}}}=\mathrm{P}\bar{\mathbf{n}}^\intercal+\bar{\mathbf{n}}\mathrm{P}^\intercal=\mathbb{E}_{\bar{\mathbf{n}}}^{\xi}$, and hence $\mathbb{P}^{\xi}_{\bar{\mathbf{n}}}=\mathbb{E}_{\bar{\mathbf{n}}}^{\xi}-\mathbb{F}_{\bar{\mathbf{n}}}=[\mathbf{0}]$. 
\end{proof}

The result below follows directly from Theorem \ref{thm:=00005B=00005CP-decomposition=00005D}:
\begin{cor}
\label{cor:[P]-reciporcity}
    For any $\mathcal{P}$-decomposition $\mathbb{F}=\mathbb{E}^{\xi}-\mathbb{P}^{\xi}$ of a magnetic algebra $(\mathbb{R}^3,\mathbb{F})$ with planarity, both $\mathbb{E}^\xi$ and $\mathbb{P}^\xi$ are reciprocal in the sense that $\mathbb{E}^\xi_{\mathbf{M}}\mathbf{m}=\mathbb{E}^\xi_{\mathbf{m}}\mathbf{M}$ and $\mathbb{P}^\xi_{\mathbf{M}}\mathbf{m}=\mathbb{P}^\xi_{\mathbf{m}}\mathbf{M}$ for any $\mathbf{M},\mathbf{m}\in \mathbb{R}^3$.
\end{cor}
\begin{proof}
   It is straightforward to check $\mathbb{E}^\xi_{\mathbf{M}}\mathbf{m}=\mathbb{E}^\xi_{\mathbf{m}}\mathbf{M}$ with \eqref{eq:=00005Bparameterization=00005D-=00005CE}, and then $\mathbb{P}^\xi_{\mathbf{M}}\mathbf{m}=\mathbb{P}^\xi_{\mathbf{m}}\mathbf{M}$ follows since $\mathbb{P}^\xi=\mathbb{E}^\xi-\mathbb{F}$. 
\end{proof}
\vspace{0.2cm}
\begin{rem}
    Note that when taking $\xi=0$ in \eqref{eq:=00005Bparameterization=00005D-=00005CE}, the $\mathcal{P}$-decomposition $\mathbb{F}=\mathbb{E}^{0}-\mathbb{P}^{0}$ is exactly the one given in Theorem \ref{thm:=00005Bspecific=00005D=00005CF=00003D=00005CE-=00005CP}. In the following discussion, unless specified otherwise, we will denote this decomposition simply by $\mathbb{F}=\mathbb{E}-\mathbb{P}$.
\end{rem}

\section{Locate $\bar{\mathbf{M}}$ on $S^2$ and Solve $\overline{\mathcal{M}}:=\underset{\mathbf{M}\in S^2}{\arg\max}|\lambda_{\mathbf{M}}|$}\label{sec:argmax-Locating}

With the structural results obtained in Section \ref{sec:AlgebraicStructure}, we are able to locate the ``dipole moment''
$\bar{\mathbf{M}}\in S^{2}\cap\overline{\mathcal{M}}$ at which the maximal ``translational
force'' $\bar{\lambda}$ is achieved. This is done in this section by proving Theorem \ref{thm:analytic-argmax}, which fully describes the set $\overline{\mathcal{M}}$ and gives analytical solutions to the optimization problem \eqref{eq:argmax|lambda|}. 

We first introduce the following preliminary result:
\begin{thm}[also see \cite{Zhengya2025AccurateWorstMagGradForce}]
\label{thm:MaximalLargestEigenvalue-all=00005CF_M}Let $(\mathbb{R}^{3},\mathbb{F})$
be a magnetic algebra with $\mathbb{F}$-invariant plane $\mathcal{P}$ and $||\mathrm{P}||>0$, and 
suppose that $\big|\big|\mathbb{F}_{\bar{\mathbf{M}}}\bar{\mathbf{m}}\big|\big|=\bar{\lambda}$ for some $\bar{\mathbf{M}},\bar{\mathbf{m}}\in S^2$.
If $\det\mathbb{F}_{\bar{\mathbf{M}}}\cdot\det\mathbb{F}_{\bar{\mathbf{m}}}=0$,
then $\bar{\lambda}=\lambda_{\mathbf{M}_{\mathbb{F}}}$ with $\mathbf{M}_{\mathbb{F}}\in S^{2}$
being \textbf{any} eigenvector associated to the largest eigenvalue of $\mathbb{F}^{\intercal}\mathbb{F}$;
otherwise, $\bar{\mathbf{M}}=\pm\bar{\mathbf{m}}$ and $\mathbb{F}_{\bar{\mathbf{M}}}\bar{\mathbf{M}}=\pm\bar{\lambda}\bar{\mathbf{M}}$. Moreover, if $\mathbb{F}_{\bar{\mathbf{M}}}\bar{\mathbf{M}}=\pm\bar{\lambda}\bar{\mathbf{M}}$
holds, then
\begin{equation}
\label{eq:[F_M]M||M}
    \big(\lambda_{\bar{\mathbf{M}}}-2\bar{\mathbf{M}}\cdot\mathrm{P}\big)\bar{\mathbf{M}}=\mathrm{P}-\mathbb{P}_{\bar{\mathbf{M}}}\bar{\mathbf{M}},
\end{equation}
and as a result, either $\bar{\mathbf{M}}\subset\mathcal{P}$, or,
$\mathbb{P}_{\bar{\mathbf{M}}}\bar{\mathbf{M}}=\mathrm{P}$ with
$\bar{\lambda}=|\lambda_{\bar{\mathbf{M}}}|=2|\bar{\mathbf{M}}\cdot\mathrm{P}|$.
\end{thm}
\begin{proof}
Due to the reciprocity, $\mathbb{F}_{\bar{\mathbf{M}}}\bar{\mathbf{m}}=\mathbb{F}_{\bar{\mathbf{m}}}\bar{\mathbf{M}}$.
Since 
\[
\bar{\lambda}=\big|\big|\mathbb{F}_{\bar{\mathbf{M}}}\bar{\mathbf{m}}\big|\big|=\big|\big|\mathbb{F}_{\bar{\mathbf{m}}}\bar{\mathbf{M}}\big|\big|=\underset{\mathbf{M},\mathbf{m}\in S^{2}}{\max}\mathbb{F}_{\mathbf{M}}\mathbf{m},
\]
it holds $|\lambda_{\bar{\mathbf{M}}}|=|\lambda_{\bar{\mathbf{m}}}|=\bar{\lambda}\geq|\lambda_{\mathbf{M}_{\mathbb{F}}}|$.
In the case $\det\mathbb{F}_{\bar{\mathbf{M}}}\cdot\det\mathbb{F}_{\bar{\mathbf{m}}}=0$,
either $\det\mathbb{F}_{\bar{\mathbf{M}}}=0$ or $\det\mathbb{F}_{\bar{\mathbf{m}}}=0$
holds, or equivalently, either $r_{\bar{\mathbf{M}}}=1$ or $r_{\bar{\mathbf{m}}}=1$,
and then by Theorem \ref{thm:=00005Cbar=00007B=00005Clambda=00007D=00003D=00005Clambda_=00007BM_F=00007D}
it holds $\lambda_{\mathbf{M}_{\mathbb{F}}}^{2}=\bar{\lambda}^{2}$. 

Now suppose that $\det\mathbb{F}_{\bar{\mathbf{M}}}\cdot\det\mathbb{F}_{\bar{\mathbf{m}}}\neq0$.
Since $\mathbb{F}_{\bar{\mathbf{M}}}\in\mathrm{SyM}_{3}^{tr=0}$,
this implies $\lambda_{\bar{\mathbf{M}}}$ to be of multiplicity $1$.
As a result, from $\big|\big|\mathbb{F}_{\bar{\mathbf{M}}}\bar{\mathbf{m}}\big|\big|=|\lambda_{\bar{\mathbf{M}}}|$
we know that $\bar{\mathbf{m}}$ has to be an eigenvector associated
to $\lambda_{\bar{\mathbf{M}}}$, i.e., $\mathbb{F}_{\bar{\mathbf{M}}}\bar{\mathbf{m}}=\lambda_{\bar{\mathbf{M}}}\bar{\mathbf{m}}$.
Due to the same reason, it also holds $\mathbb{F}_{\bar{\mathbf{m}}}\bar{\mathbf{M}}=\lambda_{\bar{\mathbf{m}}}\bar{\mathbf{M}}$,
and hence
\[
\lambda_{\bar{\mathbf{m}}}\bar{\mathbf{M}}=\mathbb{F}_{\bar{\mathbf{m}}}\bar{\mathbf{M}}=\mathbb{F}_{\bar{\mathbf{M}}}\bar{\mathbf{m}}=\lambda_{\bar{\mathbf{M}}}\bar{\mathbf{m}},
\]
implying $\bar{\mathbf{M}}=\pm\bar{\mathbf{m}}$, and then $\mathbb{F}_{\bar{\mathbf{M}}}\bar{\mathbf{M}}=\lambda_{\bar{\mathbf{M}}}\bar{\mathbf{M}}=\pm\bar{\lambda}\bar{\mathbf{M}}$.

Now suppose that $\mathbb{F}_{\bar{\mathbf{M}}}\bar{\mathbf{M}}=\pm\bar{\lambda}\bar{\mathbf{M}}$.
Applying Theorem \ref{thm:=00005B=00005CP-decomposition=00005D}
yields
\[
\lambda_{\bar{\mathbf{M}}}\bar{\mathbf{M}}=\mathbb{E}_{\bar{\mathbf{M}}}\bar{\mathbf{M}}-\mathbb{P}_{\bar{\mathbf{M}}}\bar{\mathbf{M}}=\mathrm{P}+2\big(\bar{\mathbf{M}}\cdot\mathrm{P}\big)\bar{\mathbf{M}}-\mathbb{P}_{\bar{\mathbf{M}}}\bar{\mathbf{M}},
\]
from which it is straightforward to complete the proof. 
 
\end{proof}
%
%

Theorem \ref{thm:MaximalLargestEigenvalue-all=00005CF_M} narrows down the search for $\bar{\mathbf{M}}$ on $S^2$. Note that $||\mathbb{F}_{\bar{\mathbf{M}}}\bar{\mathbf{m}}||=\bar{\lambda}$ means $|\lambda_{\bar{\mathbf{M}}}|=|\lambda_{\bar{\mathbf{m}}}|=\bar{\lambda}$ and  $\bar{\mathbf{M}},\bar{\mathbf{m}}\in\overline{\mathcal{M}}$, and therefore, $\overline{\mathcal{M}}\cap\mathrm{Det}^{\mathbb{F}}_0\neq\emptyset$ if and only if $\exists\;\bar{\mathbf{M}},\bar{\mathbf{m}}\in S^2$ s.t. $||\mathbb{F}_{\bar{\mathbf{M}}}\bar{\mathbf{m}}||=\bar{\lambda}$ and $\det\mathbb{F}_{\bar{\mathbf{M}}}\cdot\det\mathbb{F}_{\bar{\mathbf{m}}}=0$. The proof of Theorem \ref{thm:analytic-argmax} is then just to elaborate on the results about the three cases in Theorem \ref{thm:MaximalLargestEigenvalue-all=00005CF_M}, which is done through Theorems \ref{thm:F-Det==0}, \ref{thm:SetA} and \ref{thm:SetB-on-P}. 

 \subsection{Case 1: $\det\mathbb{F}_{\bar{\mathbf{M}}}\cdot\det\mathbb{F}_{\bar{\mathbf{m}}}=0$}

We first look at the case $\det\mathbb{F}_{\bar{\mathbf{M}}}\cdot\det\mathbb{F}_{\bar{\mathbf{m}}}=0$. Without loss of generality, suppose that $\det\mathbb{F}_{\bar{\mathbf{M}}}=0$, and then it follows from Theorem \ref{thm:=00005Cbar=00007B=00005Clambda=00007D=00003D=00005Clambda_=00007BM_F=00007D} that $\bar{\lambda}=\lambda_{\mathbf{M}_{\mathbb{F}}}$ holds for any (Remark \ref{rem:detF=0}) eigenvector $\mathbf{M}_{\mathbb{F}}$ of $\mathbb{F}^\intercal\mathbb{F}$ associated its largest eigenvalue $\lambda_{\mathbb{F}}$. To be clear, we introduce the set $\mathcal{E}_{\lambda_{\mathbb{F}}}:=\ker[\mathbb{F}^\intercal\mathbb{F}-\lambda_{\mathbb{F}}\mathbf{Id}]$ and recall that $\mathrm{Det}^{\mathbb{F}}_0:=\{\mathbf{M}\in \mathbb{R}^3\big|\det\mathbb{F}_{\mathbf{M}}=0\}$. The set $\mathcal{M}_\mathbb{F}$ in Theorem \ref{thm:analytic-argmax} is then $\mathcal{M}_\mathbb{F}=\mathcal{E}_{\lambda_{\mathbb{F}}}\cap S^2$.

The relation among these three sets is then described by the following lemma:
\begin{lem}
\label{lem:Det=0}
    If $\mathrm{Det}^{\mathbb{F}}_0\cap\overline{\mathcal{M}}\neq \emptyset$, then $\mathcal{M}_{\mathbb{F}}=\mathcal{E}_{\lambda_{\mathbb{F}}}\cap S^2=\mathrm{Det}^{\mathbb{F}}_0\cap\overline{\mathcal{M}}$, and, $\dim\mathcal{E}_{\lambda_{\mathbb{F}}}\geq2$.
\end{lem}
\begin{proof}
    Suppose that $\mathrm{Det}^{\mathbb{F}}_0\cap\overline{\mathcal{M}}\neq \emptyset$. Take arbitrary $\mathbf{M}\in\mathrm{Det}^{\mathbb{F}}_0\cap\overline{\mathcal{M}}$ and $\mathbf{M}_{\mathbb{F}}\in\mathcal{E}_{\lambda_{\mathbb{F}}}\cap S^2$. Applying \eqref{eq:=00005Clamda^2=000026|F_M|^2} we get
    \[
    ||\mathbb{F}_{\mathbf{M}_{\mathbb{F}}}||^2=\frac{3+r_{\mathbf{M}_\mathbb{F}}^2}{2}\lambda^2_{\mathbf{M}_{\mathbb{F}}}\geq||\mathbb{F}_{\mathbb{M}}||^2=2|\bar{\lambda}|^2.
    \]
    Since $r^2_{\mathbf{M}_{\mathbb{F}}}\in[0,1]$ and $\lambda^2_{\mathbf{M}_{\mathbb{F}}}\leq|\bar{\lambda}|^2$, the inequality above simply means $r^2_{\mathbf{M}_{\mathbb{F}}}=1$ and $|\lambda_{\mathbf{M}_{\mathbb{F}}}|=\bar{\lambda}$, and, $||\mathbb{F}_{\mathbb{M}}||=||\mathbb{F}_{\mathbf{M}_{\mathbb{F}}}||=\lambda_{\mathbb{F}}$. That is, $\mathbf{M}\in\mathcal{E}_{\lambda_{\mathbb{F}}}$, and, $\mathbf{M}_{\mathbb{F}}\in\overline{\mathcal{M}}$ and $\det\mathbb{F}_{\mathbf{M}_{\mathbb{F}}}=0$. 

It remains to show $\dim\mathcal{E}_{\lambda_{\mathbb{F}}}\geq2$ in this case. Again, take $\mathbf{M}\in\mathcal{E}_{\lambda_{\mathbb{F}}}\cap S^2$. Since $\det\mathbb{F}_\mathbf{M}=0$ and $\mathbb{F}_{\mathbf{M}}\in\mathrm{SyM}_3^{tr=0}$, the eigenvalues of $\mathbb{F}_{\mathbf{M}}$ are $0$ and $\pm\bar{\lambda}$.  Suppose that $\mathbf{m}_{\pm}$ are unit eigenvectors associated to $\pm\bar{\lambda}$, respectively. Since these two vectors are orthogonal, at least one of them is not parallel to $\mathbf{M}$. Therefore, by changing $\mathbf{M}$ to $-\mathbf{M}$  if necessary, there exists $\mathbf{m}\in S^2$ such that  $\mathbb{F}_{\mathbf{M}}\mathbf{m}=\bar{\lambda}\mathbf{m}$ while $\mathbf{M}\neq\pm\mathbf{m}$. By the reciprocity of $\mathbb{F}$, it holds $\mathbb{F}_{\mathbf{m}}\mathbf{M}=\bar{\lambda}\mathbf{m}$, and then $||\mathbb{F}_{\mathbf{m}}\mathbf{M}||=\bar{\lambda}\geq|\lambda_{\mathbf{m}}|$ , from which we deduce that $|\lambda_{\mathbf{m}}|=\bar{\lambda}$. As a result, $\mathbf{M}$ is not an eigenvector of the symmetric matrix $\mathbb{F}_{\mathbf{m}}$ while $||\mathbb{F}_{\mathbf{m}}\mathbf{M}||=|\lambda_{\mathbf{m}}|$, which can happen only if $\det\mathbb{F}_{\mathbf{m}}=0$. Therefore, $\mathbf{m}$ and $\mathbf{M}$ are both in $\mathrm{Det}^{\mathbb{F}}_0\cap\overline{\mathcal{M}}=\mathcal{E}_{\lambda_{\mathbb{F}}}\cap S^2$ while they are linearly independent, and the proof is concluded.
\end{proof}

Based on Lemma \ref{lem:Det=0}, we can prove the following proposition, which fully characterizes the algebra $(\mathbb{R}^3,\mathbb{F})$ in the case $\det\mathbb{F}_{\bar{\mathbf{M}}}=0$.

\begin{prop}
\label{prop:detF=0}
    Suppose that $\det\mathbb{F}_{\bar{\mathbf{M}}}=0$  (in other words, $\mathrm{Det}^{\mathbb{F}}_0\cap\overline{\mathcal{M}}\neq \emptyset$), and, $||\mathrm{P}||>0$. Then, $\bar{\lambda}=||\mathrm{P}||$, $\mathcal{E}_{\lambda_{\mathbb{F}}}=\mathrm{span}\{\bar{\mathbf{n}},\bar{\mathrm{P}}\}$, and,  
     \begin{equation}\label{eq:F-Det=0}
            \mathbb{F}_{\bar{\mathrm{Q}}}=[\mathbf{0}] \;\;\text{and}\;\;\mathbb{F}_{\bar{\mathrm{P}}}=||\mathrm{P}||\big(\bar{\mathbf{n}}\bar{\mathbf{n}}^\intercal-\bar{\mathrm{P}}\bar{\mathrm{P}}^\intercal\big).
     \end{equation}

\end{prop}
\begin{proof}
Take $\bar{\mathbf{M}}$ to be an arbitrary element of $\mathrm{Det}^{\mathbb{F}}_0\cap\overline{\mathcal{M}}\neq \emptyset$, and note that the eigenvalues of $\mathbb{F}_{\bar{\mathbf{M}}}$ are $\pm\bar{\lambda}$ and $0$.  
First of all, we show that $\bar{\lambda}$ has to be $||\mathrm{P}||$ in this case. Since $\lambda_{\bar{\mathbf{n}}}=||\mathrm{P}||$, we always have $\bar{\lambda}\geq||\mathrm{P}||$. Now assume that $\bar{\lambda}>||\mathrm{P}||$ and argue by contradiction. Under this assumption, $\bar{\mathbf{M}}$ and $\bar{\mathbf{n}}$ are eigenvectors of the symmetric matrix $\mathbb{F}^\intercal\mathbb{F}$ associated to different eigenvalues, and thus $\bar{\mathbf{M}}\perp\bar{\mathbf{n}}$, i.e., $\bar{\mathbf{M}}\in\mathcal{P}$. Consequently,  $\bar{\mathbf{n}}$ is an eigenvector of $\mathbb{F}_{\bar{\mathbf{M}}}$ with eigenvalue $\bar{\mathbf{M}}\cdot\mathrm{P}$, and so $\bar{\mathbf{M}}\cdot\mathrm{P}\in\{0,\pm\bar{\lambda}\}$. Since $|\bar{\mathbf{M}}\cdot\mathrm{P}|\leq|\mathrm{P}|\leq\bar{\lambda}$, either $|\bar{\mathbf{M}}\cdot\mathrm{P}|=\bar{\lambda}=|\mathrm{P}|$ and then $\bar{\mathbf{M}}=\pm\bar{\mathrm{P}}$, or $\bar{\mathrm{M}}\cdot\mathrm{P}=0$ and then $\bar{\mathbf{M}}=\bar{\mathrm{Q}}$, holds. Now that $\bar{\lambda}>|\mathrm{P}|$ is assumed, $\bar{\mathbf{M}}=\pm\bar{\mathrm{Q}}$, and the arbitrariness of $\bar{\mathbf{M}}$ as an element of $\mathrm{Det}^{\mathbb{F}}_0\cap\overline{\mathcal{M}}$ simply implies $\mathrm{Det}^{\mathbb{F}}_0\cap\overline{\mathcal{M}}=\mathcal{E}_{\lambda_{\mathbb{F}}}\cap S^2=\{\pm\bar{\mathrm{Q}}\}$, which contradicts with $\dim\mathcal{E}_{\lambda_{\mathbb{F}}}\geq2$ as Lemma \ref{lem:Det=0} implies. 

We proceed to show $\mathbb{F}_{\bar{\mathrm{Q}}}=[\mathbf{0}]$ and $\mathbb{F}_{\bar{\mathrm{P}}}=||\mathrm{P}||\big(\bar{\mathbf{n}}\bar{\mathbf{n}}^\intercal-\bar{\mathrm{P}}\bar{\mathrm{P}}^\intercal\big)$. Now that $\bar{\lambda}=|\mathrm{P}|=\lambda_{\bar{\mathbf{n}}}$ and $\det\mathbf{F}_{\bar{\mathbf{n}}}=0$, $\bar{\mathbf{n}}\in\mathrm{Det}^{\mathbb{F}}_0\cap\overline{\mathcal{M}}$.
Since the dimension of the linear subspace $\mathcal{E}_{\lambda_{\mathbb{F}}}$ is no less than $2$, $\mathrm{Det}^{\mathbb{F}}_0\cap\overline{\mathcal{M}}=\mathcal{E}_{\lambda_{\mathbb{F}}}\cap S^2$ at least contains a big circle on $S^2$ running through $\bar{\mathbf{n}}$, and therefore it contains at least one $\mathbf{M}\in\mathcal{P}\cap S^2$. 

In fact, such an $\mathbf{M}$ has to be $\pm\bar{\mathbf{P}}$, and, $\dim\mathcal{E}_{\lambda_{\mathbb{F}}}=2$. We first show $\dim\mathcal{E}_{\lambda_{\mathbb{F}}}=2$. Note that $\mathbb{F}_{\mathbf{M}}$ has eigenvalues  $0$ and  $\pm\bar{\lambda}=\pm|\mathrm{P}|$. On the other hand, since $\mathbf{M}\in\mathcal{P}\cap S^2$, $\mathbf{M}\cdot\mathrm{P}$ is an eigenvalue with the eigenvector $\bar{\mathbf{n}}$. Hence, either $\mathbf{M}\cdot\mathrm{P}=0$ and then $\mathbf{M}=\pm\bar{\mathbf{Q}}$, or $\mathbf{M}\cdot\mathrm{P}=\pm|\mathrm{P}|$ and then $\mathbf{M}=\pm\bar{\mathrm{P}}$. Note that $\mathbf{M}$ has been taken to be an arbitrary vector in $\mathcal{E}_{\lambda_{\mathbb{F}}}\cap S^2\cap\mathcal{P}$, and hence this means $\mathcal{E}_{\lambda_{\mathbb{F}}}\cap S^2\cap\mathcal{P}\subset\{\pm\bar{\mathrm{P}},\pm\bar{\mathrm{Q}}\}$. However,  if $\dim\mathcal{E}_{\lambda_{\mathbb{F}}}=3$, then $\mathcal{E}_{\lambda_{\mathbb{F}}}\cap S^2\cap\mathcal{P}=S^2\cap\mathcal{P}$, which would be a circle. Therefore, $\dim\mathcal{E}_{\lambda_{\mathbb{F}}}=2$, and, 
\[\mathcal{E}_{\lambda_{\mathbb{F}}}\cap S^2\cap\mathcal{P}=\,\textbf{either}\;\;\{\pm\bar{\mathrm{P}}\}\;\;\textbf{or}\;\;\{\pm\bar{\mathrm{Q}}\}.\]

Now we show $\mathbf{M}\neq\pm\bar{\mathrm{Q}}$, i.e., $\pm\bar{\mathrm{Q}}\notin\mathcal{E}_{\lambda_{\mathbb{F}}}\cap S^2\cap\mathcal{P}$. If $\mathbf{M}=\bar{\mathrm{Q}}$, then the plane $\mathrm{span}\{\bar{\mathbf{n}},\bar{\mathrm{Q}}\}$ is contained in $\mathcal{E}_{\lambda_{\mathbb{F}}}$, and then $\mathrm{span}\{\bar{\mathbf{n}},\bar{\mathrm{Q}}\}\cap S^2\subset\mathrm{Det}^{\mathbb{F}}_0\cap\overline{\mathcal{M}}$. For any $\beta\in(0,\pi)$ and $\mathbf{m}=\cos\beta\bar{\mathrm{Q}}+\sin\beta\bar{\mathbf{n}}$, it holds $\det\mathbb{F}_{\mathbf{m}}=0$, which means that $\mathbb{F}_{\mathbf{m}}\bar{\mathrm{P}}$, $\mathbb{F}_{\mathbf{m}}\bar{\mathrm{Q}}$, $\mathbb{F}_{\mathbf{m}}\bar{\mathbf{n}}$ are linearly dependent. However, check that
\[
\mathbb{F}_{\mathbf{m}}=\cos\beta\mathbb{F}_{\bar{\mathrm{Q}}}+\sin\beta\mathbb{F}_{\bar{\mathbf{n}}},
\]
and then 
\[
\mathbb{F}_{\mathbf{m}}\bar{\mathrm{P}}=\cos\beta\mathbb{F}_{\bar{\mathrm{Q}}}\bar{\mathrm{P}}+\sin\beta\mathbb{F}_{\bar{\mathbf{n}}}\bar{\mathrm{P}}=\cos\beta\mathbb{F}_{\bar{\mathrm{Q}}}\bar{\mathrm{P}}+\sin\beta||\mathrm{P}||\bar{\mathbf{n}},
\]
\[
\mathbb{F}_{\mathbf{m}}\bar{\mathrm{Q}}=\cos\beta\mathbb{F}_{\bar{\mathrm{Q}}}\bar{\mathrm{Q}}+\sin\beta\mathbb{F}_{\bar{\mathbf{n}}}\bar{\mathrm{Q}}=\cos\beta\mathbb{F}_{\bar{\mathrm{Q}}}\bar{\mathrm{Q}},
\]
\[
\mathbb{F}_{\mathbf{m}}\bar{\mathbf{n}}=\cos\beta\mathbb{F}_{\bar{\mathrm{Q}}}\bar{\mathbf{n}}+\sin\beta\mathbb{F}_{\bar{\mathbf{n}}}\bar{\mathbf{n}}=\sin\beta\mathrm{P}.
\]
Observe that $\mathbb{F}_{\mathbf{m}}\bar{\mathrm{Q}}$ and $\mathbb{F}_{\mathbf{m}}\bar{\mathbf{n}}$  both lie in the plane $\mathcal{P}$, while $\mathbb{F}_{\mathbf{m}}\bar{\mathrm{P}}$ has a component along $\bar{\mathbf{n}}$. Therefore, these three vectors are linearly dependent for any $\beta$ if and only if $\mathbb{F}_{\bar{\mathrm{Q}}}\bar{\mathrm{Q}}$ and $\mathrm{P}$ are collinear.  Since $\mathbb{F}_{\bar{\mathrm{Q}}}\bar{\mathbf{n}}=\mathbf{0}$,  it holds
\[
0=\mathrm{tr}\mathbb{F}_{\bar{\mathrm{Q}}}=\bar{\mathrm{P}}^\intercal\mathbb{F}_{\bar{\mathrm{Q}}}\bar{\mathrm{P}}+\bar{\mathrm{Q}}^\intercal\mathbb{F}_{\bar{\mathrm{Q}}}\bar{\mathrm{Q}}+\bar{\mathbf{n}}^\intercal\mathbb{F}_{\bar{\mathrm{Q}}}\bar{\mathbf{n}}=\bar{\mathrm{P}}^\intercal\mathbb{F}_{\bar{\mathrm{Q}}}\bar{\mathrm{P}}.
\]
Therefore, with respect to the orthonormal basis $\bar{\mathrm{P}},\bar{\mathrm{Q}},\bar{\mathbf{n}}$, the matrix of $\mathbb{F}_{\bar{\mathrm{Q}}}$ takes the form
\[
\big[\mathbb{F}_{\bar{\mathrm{Q}}}\big]=
\left[\begin{array}{ccc}
0 & \lambda_{\bar{\mathrm{Q}}} & 0\\
\lambda_{\bar{\mathrm{Q}}} & 0 & 0\\
0 & 0& 0
\end{array}\right]=
\left[\begin{array}{ccc}
0 & \pm||\mathrm{P}|| & 0\\
\pm||\mathrm{P}|| & 0 & 0\\
0 & 0& 0
\end{array}\right],
\]
from which we get  $\mathbb{F}_{\bar{\mathrm{P}}}\bar{\mathrm{Q}}=\mathbb{F}_{\bar{\mathrm{Q}}}\bar{\mathrm{P}}=\lambda_{\bar{\mathrm{Q}}}\bar{\mathrm{Q}}=\pm||\mathrm{P}||\bar{\mathrm{Q}}$. Now that we have, in addition, $\mathbb{F}_{\bar{\mathrm{P}}}\bar{\mathbf{n}}=||\mathrm{P}||\bar{\mathbf{n}}$ , $\mathrm{tr}\mathbb{F}_{\bar{\mathrm{P}}}=0$, and $||\mathbb{F}_{\bar{\mathrm{P}}}||^2\leq\lambda_{\mathbb{F}}=2||\mathrm{P}||^2$, we deduce that $\mathbb{F}_{\bar{\mathrm{P}}}\bar{\mathrm{P}}=\mathbf{0}$ and $\lambda_{\bar{\mathrm{Q}}}=-||\mathrm{P}||$. As a result, $\det\mathbb{F}_{\bar{\mathrm{P}}}=0$ and $||\mathbb{F}_{\bar{\mathrm{P}}}||^2=2||\mathrm{P}||^2=\lambda_{\mathbb{F}}$, and then, $\bar{\mathrm{P}}\in\mathcal{E}_{\lambda_{\mathbb{F}}}$. This implies $\dim\mathcal{E}_{\lambda_{\mathbb{F}}}=3$, which contradicts with the conclusion $\dim\mathcal{E}_{\lambda_{\mathbb{F}}}=2$ proved above. Therefore, $\mathbf{M}$ cannot be $\bar{\mathrm{Q}}$ (or $-\bar{\mathrm{Q}}$).

It remains to show that $\mathcal{E}_{\lambda_{\mathbb{F}}}\cap S^2\cap\mathcal{P}=\{\pm\bar{\mathrm{P}}\}$ is possible, and in this case,  $ \mathbb{F}_{\bar{\mathrm{Q}}}=[\mathbf{0}] $ and $\mathbb{F}_{\bar{\mathrm{P}}}=||\mathrm{P}||\big(\bar{\mathbf{n}}\bar{\mathbf{n}}^\intercal-\bar{\mathrm{P}}\bar{\mathrm{P}}^\intercal\big)$. We first check that, given any orthonormal basis $\bar{\mathbf{n}},\bar{\mathrm{P}}=\frac{\mathrm{P}}{||\mathrm{P}||},\bar{\mathrm{Q}}$ , and correspondingly, three matrices $ \mathbb{F}_{\bar{\mathrm{Q}}}=[\mathbf{0}] $,  $\mathbb{F}_{\bar{\mathrm{P}}}=||\mathrm{P}||\big(\bar{\mathbf{n}}\bar{\mathbf{n}}^\intercal-\bar{\mathrm{P}}\bar{\mathrm{P}}^\intercal\big)$ and $\mathbb{F}_{\bar{\mathbf{n}}}=\mathrm{P}\bar{\mathbf{n}}^\intercal+\bar{\mathbf{n}}\mathrm{P}$, the mapping
\begin{equation}
    \mathbf{M}\mapsto\mathbb{F}_{\mathbf{M}}=(\mathbf{M}\cdot\bar{\mathrm{P}})\mathbb{F}_{\bar{\mathrm{P}}}+(\mathbf{M}\cdot\bar{\mathbf{n}})\mathbb{F}_{\bar{\mathbf{n}}}
\end{equation}
defines a magnetic algebra with an invariant plane $\mathcal{P}=\mathrm{span}\{\bar{\mathrm{P}},\bar{\mathrm{Q}}\}$. Indeed, the $\mathbb{F}$ defined above is a linear map taking values in $\mathrm{SyM}_3^{tr=0}$, and, $\mathcal{P}$ is indeed an $\mathbb{F}$-invariant plane. Based on the linearity, for the reciprocity it suffices to check the commutativity on the basis $\mathbb{F}_{\bar{\mathrm{P}}}\bar{\mathbf{n}}=\mathbb{F}_{\bar{\mathbf{n}}}{\bar{\mathrm{P}}}$, $\mathbb{F}_{\bar{\mathrm{Q}}}\bar{\mathbf{n}}=\mathbb{F}_{\bar{\mathbf{n}}}{\bar{\mathrm{Q}}}$, and $\mathbb{F}_{\bar{\mathrm{P}}}\bar{\mathrm{Q}}=\mathbb{F}_{\bar{\mathrm{Q}}}{\bar{\mathrm{P}}}$. To see $\mathcal{E}_{\lambda_{\mathbb{F}}}\cap S^2\cap\mathcal{P}=\{\pm\bar{\mathrm{P}}\}$, we first check that besides $\bar{\mathbf{n}}$, $\bar{\mathrm{P}},\bar{\mathrm{Q}}$ are also eigenvalues of $\mathbb{F}^\intercal\mathbb{F}$. Since $\mathbb{F}_{\bar{\mathrm{Q}}}=[\mathbf{0}]$, $\bar{\mathrm{Q}}$ is an eigenvector with eigenvalue $0$ of  $\mathbb{F}^\intercal\mathbb{F}$, and then $\bar{\mathrm{P}}$ is another (unit) eigenvector. As a result,  $\bar{\mathbf{n}}$, $\bar{\mathrm{P}},\bar{\mathrm{Q}}$   constitute an orthonormal eigen-basis. The corresponding eigenvalues are 
\[
\bar{\mathbf{n}}^\intercal[\mathbb{F}^\intercal\mathbb{F}]\bar{\mathbf{n}}=||\mathbb{F}_{\bar{\mathbf{n}}}||^2=2||\mathrm{P}||^2,\;\;
\\
\\
\bar{\mathrm{P}}^\intercal[\mathbb{F}^\intercal\mathbb{F}]\bar{\mathrm{P}}=||\mathbb{F}_{\bar{\mathrm{P}}}||^2=2||\mathrm{P}||^2,\;\;\text{and},\;\;
\\
\\
\bar{\mathrm{Q}}^\intercal[\mathbb{F}^\intercal\mathbb{F}]\bar{\mathrm{Q}}=||\mathbb{F}_{\bar{\mathrm{Q}}}||^2=0,
\]
from which it is clear that $\lambda_{\mathbb{F}}=2||\mathrm{P}||^2$, and  $\mathcal{E}_{\lambda_{\mathbb{F}}}=\mathrm{span}\{\bar{\mathbf{n}},\bar{\mathrm{P}}\}$.

Now we conclude the whole proof by showing the converse, that when  $\mathcal{E}_{\lambda_{\mathbb{F}}}\cap S^2\cap\mathcal{P}=\{\pm\bar{\mathrm{P}}\}$, then  \eqref{eq:F-Det=0} holds. Since $\bar{\mathrm{P}}\in\mathcal{E}_{\lambda_{\mathbb{F}}}\cap S^2=\mathrm{Det}^{\mathbb{F}}_0\cap\overline{\mathcal{M}}$, $\mathbb{F}_{\bar{\mathrm{P}}}$ has eigenvalues $\pm||\mathrm{P}||=\pm\bar{\lambda}$ and $0$. By Theorem \ref{thm:StruturalThm-PlaneConfig} $\bar{\mathbf{n}}$ is an eigenvector of $\mathbb{F}_{\bar{\mathrm{P}}}$  with eigenvalue $||\mathrm{P||}$. Let  $\mathbf{m}_0$ be an unit eigenvector with eigenvalue $-||\mathrm{P}||$, and then  $\mathbf{m}_0\perp\bar{\mathbf{n}}$ with $\mathbb{F}_{\mathbf{m}_0}\bar{\mathrm{P}}=\mathbb{F}_{\bar{\mathrm{P}}}\mathbf{m}_0=||\mathrm{P}||\mathbf{m}_0$. As a result, 
\[
|\lambda_{\mathbf{m}_0}|\geq||\mathbb{F}_{\mathbf{m}_0}\bar{\mathrm{P}}||=\bar{\lambda}\geq|\lambda_{\mathbf{m}_0}|
\]
and consequently, $\lambda_{\mathbf{m}_0}=\bar{\lambda}$.  If $\bar{\mathrm{P}}\neq\mathbf{m}_0$, then this would mean $\det\mathbb{F}_{\mathbf{m}_0}=0$ and then 
\[
\mathbf{m}_0\in\mathrm{Det}^{\mathbb{F}}_0\cap\overline{\mathcal{M}}\cap\mathcal{P}=\mathcal{E}_{\lambda_{\mathbb{F}}}\cap S^2\cap\mathcal{P}=\{\pm\bar{\mathrm{P}}\}.
\]
Therefore, $\mathbf{m}_0=\bar{\mathrm{P}}$ and it is an eigenvector of $\mathbb{F}_{\bar{\mathrm{P}}}$ with eigenvalue $-||\mathrm{P}||$. This implies $\bar{\mathrm{Q}}$ to the eigenvector with eigenvalue $0$ of $\mathbb{F}_{\bar{\mathrm{P}}}$, and then $\mathbb{F}_{\bar{\mathrm{P}}}=||\mathrm{P}||\big(\bar{\mathbf{n}}\bar{\mathbf{n}}^\intercal-\bar{\mathrm{P}}\bar{\mathrm{P}}^\intercal\big)$ holds. To see $\mathbb{F}_{\bar{\mathrm{Q}}}=[\mathbf{0}]$, note that $\bar{\mathbf{n}}$ is an eigenvector with eigenvalue $0=\mathrm{P}\cdot\bar{\mathrm{Q}}$, and so is $\bar{\mathrm{P}}$ since $\mathbb{F}_{\bar{\mathrm{Q}}}\bar{\mathrm{P}}=\mathbb{F}_{\bar{\mathrm{P}}}\bar{\mathrm{Q}}=\mathbf{0}$. Since  $\mathbb{F}_{\bar{\mathrm{Q}}}$ is symmetric and has zero trace, this means it is the zero matrix, and the proof is done.

\end{proof}

    The theorem below and its corollary follow directly from Proposition \ref{prop:detF=0}
    \begin{thm}\label{thm:F-Det==0}
    Suppose that $(\mathbb{R}^3,\mathbb{F})$ is a magnetic algebra with planarity. Then,
         $\mathrm{Det}^{\mathbb{F}}_0\cap\overline{\mathcal{M}}\neq\emptyset$ if and only if $\mathbb{F}$ is given by
    \[
    \mathbb{F}_{\mathbf{M}}=(\mathbf{M}\cdot\bar{\mathbf{n}})[\mathrm{P}\bar{\mathbf{n}}^\intercal+\bar{\mathbf{n}}\mathrm{P}^\intercal]+(\mathrm{M}\cdot\mathrm{P})[\bar{\mathbf{n}}\bar{\mathbf{n}}^\intercal-\bar{\mathrm{P}}\bar{\mathrm{P}}^\intercal],\;\;\forall\mathbf{M}\in\mathbb{R}^3.
    \]

    \end{thm}
 \begin{cor}
       $\mathrm{Det}^{\mathbb{F}}_0\cap\overline{\mathcal{M}}\neq\emptyset$ if and only if $\mathrm{Det}^{\mathbb{F}}_0=\mathbb{R}^3$ (i.e., $\det\mathbb{F}_{\mathbf{M}}=0$ for all $\mathbf{M}\in\mathbb{R}^3$), and in this case, $\overline{\mathcal{M}}=\mathrm{span}\{\bar{\mathbf{n}},\mathrm{P}\}\cap S^2$.
  \end{cor}

Due to Theorem \ref{thm:MaximalLargestEigenvalue-all=00005CF_M},
$\mathrm{Det}^{\mathbb{F}}_0\cap\overline{\mathcal{M}}=\emptyset$ implies $\mathbb{F}_{\bar{\mathbf{M}}}\bar{\mathbf{M}}=\pm\bar{\lambda}\bar{\mathbf{M}}$ for any $\bar{\mathbf{M}}\in\overline{{\mathcal{M}}}$. In the meantime, from Theorem \ref{thm:F-Det==0}, we can also deduce that if $\mathrm{Det}^{\mathbb{F}}_0\cap\overline{\mathcal{M}}\neq\emptyset$, then 
\[
\mathbb{F}_{(-\bar{\mathrm{P}})}(-\bar{\mathrm{P}})=\mathbb{F}_{\bar{\mathrm{P}}}\bar{\mathrm{P}}=-||\mathrm{P}||\bar{\mathrm{P}}=\bar{\lambda}(-\bar{\mathrm{P}}).
\]
Therefore,
\begin{cor}
    Suppose that $(\mathbb{R}^3,\mathbb{F})$ is a magnetic algebra with two dimensional sub-algebra $\mathcal{P}$, and $\mathbb{F}_{\bar{\mathbf{n}}}\bar{\mathbf{n}}=\mathrm{P}\neq\mathbf{0}$. Then, there is always some $\bar{\mathbf{M}}$ s.t. $\mathbb{F}_{\bar{\mathbf{M}}}\bar{\mathbf{M}}=\bar{\lambda}\bar{\mathbf{M}}$, and then
    \[
    \bar{\lambda}=\max_{\mathbf{M}\in S^2}\mathbf{M}^\intercal\mathbb{F}_{\mathbf{M}}\mathbf{M}.
    \]
\end{cor}

\subsection{Case 2: $\mathbb{P}_{\bar{\mathbf{M}}}\bar{\mathbf{M}}=\mathrm{P}$}

Next, we shall discuss the case where $\mathbb{P}_{\bar{\mathbf{M}}}\bar{\mathbf{M}}=\mathrm{P}$ (with $\bar{\lambda}=|\lambda_{\bar{\mathbf{M}}}|=2|\bar{\mathbf{M}}\cdot\mathrm{P}|$). Since $\mathbb{P}_{\mathbf{M}}\mathbf{m}=\mathbb{P}_{\mathbf{m}}\mathbf{M}$ for any $\mathbf{M},\mathbf{m}\in \mathbb{R}^3$ and $\mathbb{P}_{\bar{\mathbf{n}}}=[\mathbf{0}]$, it holds
\begin{equation}\label{eq:[P_M]M=-P-on-P}
\mathbb{P}_{\bar{\mathbf{M}}}\bar{\mathbf{M}}=\mathbb{P}_{\bar{\mathbf{M}}_\perp}\bar{\mathbf{M}}_\perp=\mathrm{P},
\end{equation}
where $\bar{\mathbf{M}}_\perp:=\bar{\mathbf{M}}-(\bar{\mathbf{M}}\cdot\bar{\mathbf{n}})\bar{\mathbf{n}}$ is the orthogonal projection of $\bar{\mathbf{M}}$ on the plane $\mathcal{P}$. Therefore, we can narrow down the search for $\bar{\mathbf{M}}_\perp$ on $\mathcal{P}$ to solving on the plane the following quadratic equation:
\begin{equation}
\label{eq:QuaEq-P}
\bar{\mathrm{Q}}^\intercal\mathbb{P}_{\mathbf{m}}\mathbf{m}=\mathbf{m}^\intercal\mathbb{P}_{\bar{\mathrm{Q}}}\mathbf{m}=0,\;\mathbf{m}\in \mathcal{P}.
\end{equation}
Equation \ref{eq:QuaEq-P} has solutions since the $3\times3$ matrix $\mathbb{P}_{\bar{\mathrm{Q}}}$ is symmetric and $\mathrm{tr}\mathbb{P}_{\bar{\mathrm{Q}}}=\det\mathbb{P}_{\bar{\mathrm{Q}}}=0$, and every solution to the equation takes the form of 
\[
\mathbf{m}=r(\mathbf{m}_+\pm\mathbf{m}_-),
\]
where $r\in\mathbb{R}$ and $\mathbf{m}_+,\mathbf{m}_-\in S^2\cap\mathcal{P}$ are orthogonal unit eigenvectors of $\mathbb{P}_{\bar{\mathrm{Q}}}$. These results are formally stated in the proposition below:
\begin{prop} 
\label{prop:P_Q}
Suppose that $\mathbb{P}_{\bar{\mathrm{Q}}}\neq[\mathbf{0}]$.
    $\mathbb{P}_{\bar{\mathrm{Q}}}$ has eigenvectors $\mathbf{m}_+,\mathbf{m}_-\in S^2\cap\mathcal{P}$ with $\mathbf{m}_+\perp\mathbf{m}_-$ such that $\mathbb{P}_{\bar{\mathrm{Q}}}\mathbf{m}_+=\beta_{\bar{\mathrm{Q}}}\mathbf{m}_+$ and $\mathbb{P}_{\bar{\mathrm{Q}}}\mathbf{m}_-=-\beta_{\bar{\mathrm{Q}}}\mathbf{m}_-$. As a result, $\mathbf{m}\in\mathcal{P}$ solves Eq. \eqref{eq:QuaEq-P} if and only if $\mathbf{m}=r(\mathbf{m}_+\pm\mathbf{m}_-)$ with $r\in\mathbb{R}$.
\end{prop}
\begin{proof}
    $\mathbb{P}_{\bar{\mathrm{Q}}}=\mathbb{E}_{\bar{\mathrm{Q}}}-\mathbb{F}_{\bar{\mathrm{Q}}}$. Since $\mathbb{E}_{\bar{\mathrm{Q}}}=\bar{\mathrm{Q}}\mathrm{P}^\intercal+\mathrm{P}\bar{\mathrm{Q}}^\intercal$, it holds $\mathrm{tr}\mathbb{E}_{\bar{\mathrm{Q}}}=\mathrm{P}^\intercal\bar{\mathrm{Q}}+\bar{\mathrm{Q}}^\intercal\mathrm{P}=0$, and then $\mathrm{tr}\mathbb{P}_{\bar{\mathrm{Q}}}=\mathrm{tr}\mathbb{E}_{\bar{\mathrm{Q}}}-\mathrm{tr}\mathbb{F}_{\bar{\mathrm{Q}}}=0$. Check that $\bar{\mathbf{n}}$ is an eigenvector of $\mathbb{P}_{\bar{\mathrm{Q}}}$ with eigenvalue $0$:
    \[
\mathbb{P}_{\bar{\mathrm{Q}}}\bar{\mathbf{n}}=\mathbb{P}_{\bar{\mathbf{n}}}\bar{\mathrm{Q}}=[\mathbf0]\bar{\mathrm{Q}}=\mathbf{0}.
    \]
    Therefore, the other eigenvalues are $\pm\beta_{\bar{\mathrm{Q}}}$, and the associated eigenvectors lie on $\mathcal{P}$. Taking $\mathbf{m}_+, \mathbf{m}_-$, any vector on $\mathcal{P}$ takes the form $\mathbf{m}=a\mathbf{m}_++b\mathbf{m}_-$, and then $ \mathbf{m}^\intercal\mathbb{P}_{\bar{\mathrm{Q}}}\mathbf{m}=(a^2-b^2)\beta_{\bar{\mathrm{Q}}}=0$ if and only if $a^2=b^2$.
\end{proof}

Now we continue to discuss the solutions to the equation
\begin{equation}
\label{eq:[P_m]m=-P}
    \mathbb{P}_{\mathbf{M}}\mathbf{M}=\mathrm{P},\;\mathbf{M}\in S^2
\end{equation}
Take $\hat{\mathbf{m}}_+:=\frac{\mathbf{m}_++\mathbf{m}_-}{||\mathbf{m}_++\mathbf{m}_-||}=\frac{\mathbf{m}_++\mathbf{m}_-}{\sqrt{2}}$ and $\hat{\mathbf{m}}_-:=\frac{\mathbf{m}_+-\mathbf{m}_-}{||\mathbf{m}_+-\mathbf{m}_-||}=\frac{\mathbf{m}_+-\mathbf{m}_-}{\sqrt{2}}$. Both of $\mathbb{P}_{\hat{\mathbf{m}}_{\pm}}\hat{\mathbf{m}}_{\pm}$ are then collinear with $\mathrm{P}$, and
\begin{equation}\label{eq:Trace[P]=1/2Trace+Remainder}
\mathbb{P}_{\hat{\mathbf{m}}_{\pm}}\hat{\mathbf{m}}_{\pm}=\frac{1}{2}\big(\mathbb{P}_{\mathbf{m}_{+}}\mathbf{m}_{+}+\mathbb{P}_{\mathbf{m}_{-}}\mathbf{m}_{-}\big)\pm\mathbb{P}_{\mathbf{m_+}}\mathbf{m_-}.
\end{equation}
\eqref{eq:[P_m]m=-P} has a solution if and only if at least one of $\gamma_\pm=\bar{\mathrm{P}}^\intercal\mathbb{P}_{\hat{\mathbf{m}}_{\pm}}\hat{\mathbf{m}}_{\pm}$ is larger than $||\mathrm{P}||=\bar{\mathrm{P}}^\intercal\mathrm{P}$, and we have the following solutions:
\begin{equation}\label{eq:Solution-[P_m]m=P}
    \mathbf{M}:=\pm\sqrt{1-\frac{||\mathrm{P}||}{\gamma}}\bar{\mathbf{n}}\pm\sqrt{\frac{||\mathrm{P}||}{\gamma}}\hat{\mathbf{m}}.
\end{equation}

Note that $\mathrm{tr}\mathbb{P}_{\bar{\mathrm{P}}}=\mathrm{tr}\mathbb{E}_{\bar{\mathrm{P}}}=5||\mathrm{P}||$ while $\mathbb{P}_{\bar{\mathrm{P}}}\bar{\mathbf{n}}=\mathbb{P}_{\bar{\mathbf{n}}}\bar{\mathrm{P}}=\mathbf{0}$, and thence it holds
\begin{equation}\label{eq:Trace[P]}
\hat{\mathbf{m}}^\intercal_{+}\mathbb{P}_{\bar{\mathrm{P}}}\hat{\mathbf{m}}_{+}+\hat{\mathbf{m}}^\intercal_{-}\mathbb{P}_{\bar{\mathrm{P}}}\hat{\mathbf{m}}_{-}=\mathbf{m}^\intercal_{+}\mathbb{P}_{\bar{\mathrm{P}}}\mathbf{m}_{+}+\mathbf{m}^\intercal_{-}\mathbb{P}_{\bar{\mathrm{P}}}\mathbf{m}_{-}=5||\mathrm{P}||>2||\mathrm{P}||.
\end{equation}

As a consequence, at least one of $\hat{\mathbf{m}}^\intercal_{\pm}\mathbb{P}_{\bar{\mathrm{P}}}\hat{\mathbf{m}}_{\pm}=\gamma_\pm$ has to be larger than $||\mathrm{P}||$, and hence

\begin{thm}
\label{thm:SetA}
    The equation \eqref{eq:[P_m]m=-P} on $S^2$ always has solutions. Any solution to \eqref{eq:[P_m]m=-P} takes the form
    \[
\hat{\mathbf{M}}=\cos\hat{\theta}\cdot\bar{\mathbf{n}}+\sin{\hat{\theta}}\cdot\hat{\mathbf{m}}_\pm,
    \]
    where $\hat{\theta}$ and the sign $\pm$ are determined by
    \[\sin^2\hat{\theta}\cdot\big(\frac{5}{2}||\mathrm{P}||\pm\mathbf{m}^\intercal_+\mathbb{P}_{\bar{\mathrm{P}}}\mathbf{m}_-\big)=||\mathrm{P}||.\]
\end{thm}
\begin{proof}
   The form of $\hat{\mathbf{M}}$ is just from equation \eqref{eq:Solution-[P_m]m=P}, and hence
   \[
\bar{\mathrm{P}}^\intercal\mathbb{P}_{\hat{\mathbf{M}}}\hat{\mathbf{M}}=\sin^2\hat{\theta}\cdot\bar{\mathrm{P}}^\intercal\mathbb{P}_{\hat{\mathbf{m}}_\pm}\hat{\mathbf{m}}_\pm=||\mathrm{P}||.
   \]
   The proof is then completed by applying \eqref{eq:Trace[P]=1/2Trace+Remainder} and \eqref{eq:Trace[P]} to the above equation.
\end{proof}

\subsection{Case 3: $\mathbb{F}_{\bar{\mathbf{M}}}\bar{\mathbf{M}}=\bar{\lambda}\bar{\mathbf{M}}$ with $\bar{\mathbf{M}}\in\mathcal{P}$}

For the last case where $\bar{\mathbf{M}}\in\mathcal{P}$ and $\mathbb{F}_{\bar{\mathbf{M}}}\bar{\mathbf{M}}=\bar{\lambda}\bar{\mathbf{M}}=\lambda_{\mathcal{P}}\bar{\mathbf{M}}$, it is to solve the optimization problem on the circle $\mathcal{P}\cap S^2$:
\begin{equation}\label{eq:argmax}
    \underset{\mathbf{M}\in\mathcal{P}\cap S^2}{\arg\max}\; \mathbf{M}^\intercal\mathbb{F}_{\mathbf{M}}\mathbf{M}.
\end{equation}
Alternatively, with $\mathrm{Q}:=\bar{\mathbf{n}}\times\mathrm{P}$, one can also obtain the pair $\big(\mathbf{M},\lambda\big)=\big(\bar{\mathbf{M}},\bar{\lambda}\big)$ by solving the set of quadratic equations below:
\begin{equation}\label{eq:QuaEqs-on-P}
\begin{cases}
    \mathbf{M}^\intercal\mathbb{F}_{\mathrm{P}}\mathbf{M}-\lambda\mathrm{P}^\intercal\mathbf{M}=0;\\
    \mathbf{M}^\intercal\mathbb{F}_{\mathrm{Q}}\mathbf{M}-\lambda\mathrm{Q}^\intercal\mathbf{M}=0.
\end{cases}
\end{equation}
It is not heavy work to solve \eqref{eq:argmax} and \eqref{eq:QuaEqs-on-P} numerically on the circle $\mathcal{P}\cap S^2$. For better understanding \eqref{eq:QuaEqs-on-P} from an analytical point of view, divide both the equations by $\lambda^2$ and get
\begin{equation*}
\begin{cases}
    \frac{\mathbf{M^\intercal}}{\lambda}\mathbb{F}_{\mathrm{P}}\frac{\mathbf{M}}{\lambda}-\mathrm{P}^\intercal\frac{\mathbf{M}}{\lambda}=0;\\
    \frac{\mathbf{M^\intercal}}{\lambda}\mathbb{F}_{\mathrm{Q}}\frac{\mathbf{M}}{\lambda}-\mathrm{Q}^\intercal\frac{\mathbf{M}}{\lambda}=0.
\end{cases}
\end{equation*}
Taking $\mathbf{m}=\frac{\mathbf{M}}{\lambda}$, we see that, if a pair $(\mathbf{M}_0,\lambda_0)\in\big(\mathcal{P}\cap S^2\big)\times\mathbb{R}_+$ solves \eqref{eq:QuaEqs-on-P}, then $\mathbf{m}_0:=\frac{\mathbf{M}_0}{\lambda_0}$ lies on the intersection of the conic curves on $\mathcal{P}$ given by \eqref{eq:ConicCurves-P&Q}, that is,
\begin{equation*}
\begin{aligned}
\mathbf{m}^\intercal\mathbb{F}_{\mathrm{P}}\mathbf{m}=\mathrm{P}^\intercal\mathbf{m} \;\;\And\;\;
\mathbf{m}^\intercal\mathbb{F}_{\mathrm{Q}}\mathbf{m}=\mathrm{Q}^\intercal\mathbf{m}.
\end{aligned}
\end{equation*}
Conversely, if $\mathbf{m}_0$ lies on the curves given by \eqref{eq:ConicCurves-P&Q} on $\mathcal{P}$ and $||\mathbf{m}_0||>0$, then with $\lambda_0=\frac{1}{||\mathbf{m}_0||}$ and $\mathbf{M}_0=\frac{\mathbf{m}_0}{||\mathbf{m}_0||}=\lambda_0\mathbf{m}_0$, $(\mathbf{M}_0,\lambda_0)$ is a solution pair to \eqref{eq:QuaEqs-on-P}. This actually means that the equation $\mathbb{F}_{\mathbf{M}_0}\mathbf{M}_0=\lambda_0\mathbf{M}_0$ is satisfied, since $\mathbf{M}_0\in\mathcal{P}$ already implies
\[
\bar{\mathbf{n}}^\intercal\mathbb{F}_{\mathbf{M}_0}\mathbf{M}_0=\mathbf{M}_0^\intercal\mathbb{F}_{\bar{\mathbf{n}}}\mathbf{M}_0=0=\lambda_0\mathbf{M}_0\cdot\bar{\mathbf{n}}.
\]

\begin{thm}\label{thm:SetB-on-P}
    $\mathbf{m}\in\mathcal{P}\setminus\{\mathbf{0}\}$ solves \eqref{eq:ConicCurves-P&Q} if and only if $(\mathbf{M},\lambda)=(\frac{\mathbf{m}}{||\mathbf{m}||},\frac{1}{||\mathbf{m}||})$ solves \eqref{eq:QuaEqs-on-P}, and hence the solutions to \eqref{eq:ConicCurves-P&Q} on $\mathcal{P}\setminus\{\mathbf{0}\}$ are in one-to-one correspondence with the solutions to \eqref{eq:QuaEqs-on-P} on $\big(\mathcal{P}\cap S^2\big)\times\mathbb{R}_+$. Moreover, if $(\mathbf{M},\lambda)$ solves \eqref{eq:QuaEqs-on-P}, so does $(-\mathbf{M},-\lambda)$.
\end{thm}
\begin{rem}
    Note that $\mathbb{F}_{\mathrm{Q}}\bar{\mathbf{n}}=\mathbb{F}_{\bar{\mathbf{n}}}\mathrm{Q}=\mathbf{0}$, and therefore, the eigenvalues of $\mathbb{F}_{\mathrm{Q}}$ are $\pm\lambda_{\mathrm{Q}}$ and $0$, and the curve given by the second equation in \eqref{eq:ConicCurves-P&Q} is a hyperbola on $\mathcal{P}$. To solve \eqref{eq:ConicCurves-P&Q}, one can parameterize each component of the hyperbola and substitute the parametrization into the first equation in \eqref{eq:ConicCurves-P&Q}, reducing the problem to solving one single-variable equation.
\end{rem}

\section{Example: $\mathbb{F}^\intercal\mathbb{F}=\lambda_{\mathbb{F}}\mathbf{Id}$.}
\label{sec:Example-F*F=uniform}
The sets $\mathcal{M}_{\mathbf{F}}$, $\mathcal{A}$ and $\mathcal{B}$ in Theorem \ref{thm:analytic-argmax} may overlap with each other, and $\overline{\mathcal{M}}$ may intersect with more than one of them. For example, the algebra given by \eqref{eq:FwhenDet=0} has $\overline{\mathcal{M}}=\mathcal{M}_{\mathbb{F}}=\mathrm{span}\{\bar{\mathbf{n}},\mathrm{P}\}$. Moreover, it can be check that
\[
\begin{aligned}\mathbb{F}_{\mathbf{M}}= & (\mathbf{M}\cdot\bar{\mathbf{n}})[\bar{\mathbf{n}}\mathrm{P}^{\intercal}+\mathrm{P}\bar{\mathbf{n}}]+(\mathbf{M}\cdot\mathrm{P})[\bar{\mathbf{n}}\bar{\mathbf{n}}^{\intercal}-\bar{\mathrm{P}}\bar{\mathrm{P}}^{\intercal}]\\
= & \underbrace{\mathbf{M}\mathrm{P}^{\intercal}+\mathrm{P}\mathbf{M}^{\intercal}+(\mathbf{M}\cdot\mathrm{P})\mathbf{Id}}_{\mathbb{E}_{\mathbf{M}}}-\underbrace{\bigg(3(\mathbf{M}\cdot\mathrm{P})\bar{\mathrm{P}}\bar{\mathrm{P}}^{\intercal}+(\mathbf{M}\cdot\mathrm{P})[\bar{\mathrm{P}}\bar{\mathrm{P}}^{\intercal}+\bar{\mathrm{Q}}\bar{\mathrm{Q}}^{\intercal}]+(\mathbf{M}\cdot\mathrm{Q})[\bar{\mathrm{Q}}\bar{\mathrm{P}}^{\intercal}+\bar{\mathrm{P}}\bar{\mathrm{Q}}^{\intercal}]\bigg)}_{\mathbb{P}_{\mathbf{M}}},
\end{aligned}
\]
and therefore, $\mathbb{P}_{\bar{\mathbf{M}}_0}\bar{\mathbf{M}}_0=\mathrm{P}$ with
\[
\bar{\mathbf{M}}_0=\pm\frac{\sqrt{3}}{2}\bar{\mathbf{n}}\pm\frac{1}{2}\bar{\mathrm{P}}.
\]
Note that $\bar{\mathbf{M}}_0\in\mathrm{span}\{\bar{\mathbf{n}},\mathrm{P}\}=\overline{\mathcal{M}}$ while $\bar{\mathbf{M}}_0\notin\mathcal{P}$, and hence $\bar{\mathbf{M}}_0\in\mathcal{A}\setminus\mathcal{B}$.

Now we consider another specific case of $(\mathbb{R}^3,\mathbb{F})$ which has actually emerged in the argument for Proposition \ref{prop:detF=0}. We will see that, in this example, $\overline{\mathcal{M}}=\mathcal{A}\sqcup\mathcal{B}$ and $\mathcal{A}\cap\mathcal{P}=\emptyset$.
The construction of this algebra is specified with the following theorem:
\begin{thm}\label{thm:F*F=uniform}
    $\mathbb{F}^\intercal\mathbb{F}=\left[\begin{array}{ccc}
\lambda_{\mathbb{F}} & 0 & 0\\
0 & \lambda_{\mathbb{F}} & 0\\
0 & 0 & \lambda_{\mathbb{F}}
\end{array}\right]$ if and only if 
\begin{equation}\label{eq:Uniform-F*F}
    \mathbb{F}_{\mathbf{M}}=(\mathbf{M}\cdot\bar{\mathbf{n}})[\bar{\mathbf{n}}\mathrm{P}^\intercal+\mathrm{P}\bar{\mathbf{n}}^\intercal]+(\mathbf{M}\cdot\mathrm{P})[\bar{\mathbf{n}}\bar{\mathbf{n}}^\intercal-\bar{\mathrm{Q}}\bar{\mathrm{Q}}^\intercal]-(\mathbf{M}\cdot\mathrm{Q})[\bar{\mathrm{P}}\bar{\mathrm{Q}}^\intercal+\bar{\mathrm{Q}}\bar{\mathrm{P}}^\intercal].
\end{equation}
\end{thm}
\begin{rem}
    Like before, $\mathrm{Q}:=\bar{\mathbf{n}}\times\mathrm{P}$, and hence $||\mathrm{Q}||=||\mathrm{P}||$.
\end{rem}
\begin{rem}
    Since $\bar{\mathbf{n}}$ is always an eigenvector of $\mathbb{F}^\intercal\mathbb{F}$ with eigenvalue $2||\mathrm{P}||^2$, $\mathbb{F}^\intercal\mathbb{F}=\lambda_{\mathbb{F}}\mathbf{Id}$ implies $\lambda_{\mathbb{F}}=2||\mathrm{P}||^2$.
\end{rem}
\begin{proof}
    Supposing that \eqref{eq:Uniform-F*F} holds, we show that all the eigenvalues of $\mathbb{F}^\intercal\mathbb{F}$ are equal. It is straightforward to check that 
    \[
||\mathbb{F}_{\bar{\mathbf{n}}}||_{\mathbb{R}^9}^2=||\mathbb{F}_{\bar{\mathbf{P}}}||_{\mathbb{R}^9}^2=||\mathbb{F}_{\bar{\mathbf{Q}}}||_{\mathbb{R}^9}^2=2||\mathrm{P}||^2.
    \]
    Note that $(\bar{\mathbf{n}},\bar{\mathrm{P}},\bar{\mathrm{Q}})$ is an orthonormal basis of $\mathbb{R}^3$. Also note that $\bar{\mathbf{n}}$ is an eigenvalue of $\mathbb{F}^\intercal\mathbb{F}$, and also, $\mathbb{F}_{\bar{\mathrm{P}}}\perp\mathbb{F}_{\bar{\mathrm{Q}}}$ in $\mathbb{R}^9$ since
    \[
\bar{\mathrm{P}}^\intercal[\mathbb{F}^\intercal\mathbb{F}]\bar{\mathrm{Q}}=\mathrm{tr}[\mathbb{F}_{\bar{\mathrm{P}}}\mathbb{F}_{\bar{\mathrm{Q}}}]=\mathrm{tr}[-\mathrm{Q}\mathrm{P}^\intercal]=0.
    \]
    From these it can be deduced that 
    \[
\mathbf{M}^\intercal[\mathbb{F}^\intercal\mathbb{F}]\mathbf{M}=||\mathbb{F}_{\mathbf{M}}||^2_{\mathbb{R}^9}=2||\mathrm{\mathrm{P}}||^2,\;\forall\mathbf{M}\in S^2,
    \]
    and hence the symmetric matrix $\mathbb{F}^\intercal\mathbb
    F$ is positively definite and  all its eigenvalues of $\mathbb{F}^\intercal\mathbb{F}$ are equal.

    Conversely, suppose that all the eigenvalues of $\mathbb{F}^\intercal\mathbb{F}$ are equal, and we deduce \eqref{eq:Uniform-F*F}. Note that under this assumption, it holds $\lambda_{\mathbb{F}}=2||\mathrm{P}||^2$ and $\mathrm{span}\{\bar{\mathbf{n}},\bar{\mathrm{Q}}\}\subset\mathcal{E}_{\lambda_{\mathbf{F}}}$. The proof is then completed by Lemma \ref{lem:WeakUniform-F*F} given below with this seemly weaker condition. Note that due to the linearity of the map $\mathbb{F}$, \eqref{eq:Uniform-F*F} holds if and only if \eqref{eq:Uniform-F*F-Basis} holds (since we already know $\mathbb{F}_{\bar{\mathbf{n}}}=\bar{\mathbf{n}}\mathrm{P}^\intercal+\mathrm{P}\bar{\mathbf{n}}^\intercal$).
\end{proof}

\begin{lem}\label{lem:WeakUniform-F*F}
    Suppose that $\mathrm{span}\{\bar{\mathbf{n}},\bar{\mathrm{Q}}\}\subset\mathcal{E}_{\lambda_{\mathbf{F}}}$. It holds
    \begin{equation}\label{eq:Uniform-F*F-Basis}
\mathbb{F}_{\bar{\mathrm{P}}}=||\mathrm{P}||\cdot[\bar{\mathbf{n}}\bar{\mathbf{n}}^\intercal-\bar{\mathrm{Q}}\bar{\mathrm{Q}}^\intercal],\;\;\text{and},\;\;
\mathbb{F}_{\bar{\mathrm{Q}}}=-||\mathrm{P}||\cdot[\bar{\mathrm{P}}\bar{\mathrm{Q}}^\intercal+\bar{\mathrm{Q}}\bar{\mathrm{P}}^\intercal].
    \end{equation}
\end{lem}
\begin{proof}
    Note that $\mathrm{span}\{\bar{\mathbf{n}},\bar{\mathrm{Q}}\}\subset\mathcal{E}_{\lambda_{\mathbb{F}}}$ implies $\lambda_{\mathbb{F}}=\bar{\mathbf{n}}^\intercal\mathbb{F}^\intercal\mathbb{F}\bar{\mathbf{n}}=\mathrm{tr}\mathbb{F}^2_{\bar{\mathbf{n}}}=2||\mathrm{P}||^2$, and then also $\mathrm{tr}\mathbb{F}^2_{\bar{\mathrm{Q}}}=\lambda_{\mathbb{F}}=2||\mathrm{P}||^2$. Since $\mathbb{F}_{\bar{\mathrm{Q}}}\bar{\mathbf{n}}=\mathbf
    0$, it further implies $\lambda_{\bar{\mathrm{Q}}}=||\mathrm{P}||$, and also, $\mathbb{F}_{\bar{\mathrm{Q}}}$ is nondegenerate on $\mathcal{P}$.
    
   Now that the (nonzero) vectors in $\mathrm{span}\{\bar{\mathbf{n}},\bar{\mathrm{Q}}\}$  are all eigenvectors of $\mathbb{F}^\intercal\mathbb{F}$, it is an invariant subspace of $\mathbb{F}^\intercal\mathbb{F}$, and so is its orthogonal complement $\mathrm{span}\{\mathrm{P}\}$. As a result, 
    \[\mathrm{tr}[\mathbb{F}_{\bar{\mathrm{P}}}\mathbb{F}_{\bar{\mathrm{Q}}}]=\bar{\mathrm{P}}^\intercal[\mathbb{F}^\intercal\mathbb{F}]\bar{\mathrm{Q}}=0.\]
    Now that $\mathbb{F}_{\bar{\mathrm{Q}}}\bar{\mathbf{n}}=\mathbb{F}_{\bar{\mathbf{n}}}\bar{\mathrm{Q}}=\mathbf{0}$, it yields
    \[\bar{\mathrm{P}}^\intercal[\mathbb{F}_{\bar{\mathrm{P}}}\mathbb{F}_{\bar{\mathrm{Q}}}]\bar{\mathrm{P}}+\bar{\mathrm{Q}}^\intercal[\mathbb{F}_{\bar{\mathrm{P}}}\mathbb{F}_{\bar{\mathrm{Q}}}]\bar{\mathrm{Q}}=0.\]
    Note that $\mathbb{F}_{\bar{\mathrm{Q}}}\in\mathrm{SyM}^{tr=0}_3$ and $\mathbb{F}_{\bar{\mathrm{Q}}}\bar{\mathbf{n}}=\mathbf{0}$ imply $\mathbb{F}_{\bar{\mathrm{Q}}}^2=\left[\begin{array}{ccc}
0 & 0 & 0\\
0 & \lambda^2_{\bar{\mathrm{Q}}} & 0\\
0 & 0 & \lambda^2_{\bar{\mathrm{Q}}}
\end{array}\right]$, and then $\bar{\mathrm{Q}}^\intercal\mathbb{F}_{\bar{\mathrm{P}}}\mathbb{F}_{\bar{\mathrm{Q}}}\bar{\mathrm{Q}}=\bar{\mathrm{P}}^\intercal\mathbb{F}^2_{\bar{\mathrm{Q}}}\bar{\mathrm{Q}}=\lambda^2_{\bar{\mathrm{Q}}}\bar{\mathrm{P}}^\intercal\bar{\mathrm{Q}}=0$. Consequently, it holds
\[\bar{\mathrm{P}}^\intercal\mathbb{F}^\intercal_{\bar{\mathrm{P}}}\mathbb{F}_{\bar{\mathrm{P}}}\bar{\mathrm{Q}}=\bar{\mathrm{P}}^\intercal\mathbb{F}_{\bar{\mathrm{P}}}\mathbb{F}_{\bar{\mathrm{Q}}}\bar{\mathrm{P}}=-\bar{\mathrm{Q}}^\intercal\mathbb{F}_{\bar{\mathrm{P}}}\mathbb{F}_{\bar{\mathrm{Q}}}\bar{\mathrm{Q}}=0.\]
Therefore, $\mathbb{F}_{\bar{\mathrm{P}}}\bar{\mathrm{P}}$ is either $\mathbf{0}$ or perpendicular to $\mathbb{F}_{\bar{\mathrm{P}}}\bar{\mathrm{Q}}$, while $\mathbb{F}_{\bar{\mathrm{Q}}}\bar{\mathrm{Q}}$ is nonzero and perpendicular to $\mathbb{F}_{\bar{\mathrm{Q}}}\bar{\mathrm{P}}$, since $\mathbb{F}_{\bar{\mathrm{Q}}}$ is nondegenerate on $\mathcal{P}$. Since $\mathcal{P}$ is invariant under the mappings of both $\mathbb{F}_{\bar{\mathrm{Q}}}$ and $\mathbb{F}_{\bar{\mathrm{P}}}$, and, since $\mathbb{F}_{\bar{\mathrm{P}}}\bar{\mathrm{Q}}=\mathbb{F}_{\bar{\mathrm{Q}}}\bar{\mathrm{P}}$, the restriction of $\mathbb{F}_{\bar{\mathrm{Q}}}$ and $\mathbb{F}_{\bar{\mathrm{P}}}$ to $\mathcal{P}$ have matrix representations with respect to the orthonormal basis $(\bar{\mathrm{P}},\bar{\mathrm{Q}})$ as
\begin{equation}
    [\mathbb{F}_{\bar{\mathrm{P}}}]=
\left[\begin{array}{ccc}
 a & b\\
 b & c
\end{array}\right]
\;\;\text{and}\;\;\;
[\mathbb{F}_{\bar{\mathrm{Q}}}]=
\left[\begin{array}{ccc}
 b & c\\
 c & -b
\end{array}\right].
\end{equation}
in which either $a=-c$ or $b=0$ holds. Note that $\mathbb{F}_{\bar{\mathrm{P}}}$ has zero trace while $\bar{\mathbf{n}}^\intercal\mathbb{F}_{\bar{\mathrm{P}}}\bar{\mathbf{n}}=\bar{\mathbf{n}}^\intercal\mathbb{F}_{\bar{\mathbf{n}}}\bar{\mathrm{P}}=||\mathrm{P}||$, and hence $a+c=-||\mathrm{P}||<0$. This means $a\neq-c$ and thus $b=0$. Consequently, $2c^2=||\mathbb{F}_{\bar{\mathrm{Q}}}||^2_{\mathbb{R}^9}=2||\mathrm{P}||^2$ and thus $c=\pm||\mathrm
P||$, and then $a$ equals either $0$ or $-2||\mathrm{P}||$. Since  $||\mathbb{F}_{\bar{\mathrm{P}}}||^2_{\mathbb{R}^9}\leq\lambda_{\mathbb{F}}=2||\mathrm{P}||^2$, it has to be $a=0$ and thence $c=-||\mathrm{P}||$. As a result,
\begin{equation}
    [\mathbb{F}_{\bar{\mathrm{P}}}]=
\left[\begin{array}{ccc}
 0 & 0\\
 0 & -||\mathrm{P}||
\end{array}\right]
\;\;\text{and}\;\;\;
[\mathbb{F}_{\bar{\mathrm{Q}}}]=
\left[\begin{array}{ccc}
 0 & -||\mathrm{P}||\\
 -||\mathrm{P}|| & 0
\end{array}\right],
\end{equation}
from which \eqref{eq:Uniform-F*F-Basis} can be directly deduced.
\end{proof}

In practice, Theorem \ref{thm:F*F=uniform} means that, once we find $\mathbb{F}^\intercal\mathbb{F}=\lambda_{\mathbb{F}}\mathbf{Id}$, both $\bar{\lambda}$ and $\overline{\mathcal{M}}$ can be determined directly without further implementing Theorem \ref{thm:analytic-argmax}, since $\mathbb{F}$ is uniquely determined.
Now that in this case $\dim\mathcal{E}_{\lambda_{\mathbb{F}}}=3$,
with Proposition \ref{prop:detF=0} we deduce $\mathrm{Det}_{0}^{\mathbb{F}}\cap\overline{\mathcal{M}}=\emptyset$.
Taking an arbitrary $\bar{\mathbf{M}}\in\overline{\mathcal{M}}$,
from Theorem \ref{thm:MaximalLargestEigenvalue-all=00005CF_M} we know that $\mathbb{F}_{\bar{\mathbf{M}}}\bar{\mathbf{M}}=\pm\bar{\lambda}\bar{\mathbf{M}}$.
Since $(\bar{\mathbf{n}},\bar{\mathrm{P}},\bar{\mathrm{Q}})$ is an
orthonormal frame and $\bar{\mathbf{M}}\in S^{2}$, there exists $(a,b,c)\in S^{2}$
such that $\bar{\mathbf{M}}=a\bar{\mathbf{n}}+b\bar{\mathrm{P}}+c\bar{\mathrm{Q}}$.
Note that $||\mathrm{Q}||=||\mathrm{P}||$, and then we have
\[
\begin{aligned}\mathbb{F}_{\bar{\mathbf{M}}}\bar{\mathbf{M}}= & a[\mathrm{P}\bar{\mathbf{n}}^{\intercal}+\bar{\mathbf{n}}\mathrm{P}^{\intercal}]\bar{\mathbf{M}}+(b\cdot||\mathrm{P}||)[\bar{\mathbf{n}}\bar{\mathbf{n}}^{\intercal}-\bar{\mathrm{Q}}\bar{\mathrm{Q}}^{\intercal}]\bar{\mathbf{M}}-(c\cdot||\mathrm{Q}||)[\bar{\mathrm{P}}\bar{\mathrm{Q}}^{\intercal}+\bar{\mathrm{Q}}\bar{\mathrm{P}}^{\intercal}]\bar{\mathbf{M}}\\
= & a^{2}\mathrm{P}+ab||\mathrm{P}||\bar{\mathbf{n}}+(ab\cdot||\mathrm{P}||)\bar{\mathbf{n}}-(bc\cdot||\mathrm{P}||)\bar{\mathrm{Q}}-(c^{2}\cdot||\mathrm{Q}||)\bar{\mathrm{P}}-(bc\cdot||\mathrm{Q}||)\bar{\mathrm{Q}}\\
= & ||\mathrm{P}||\cdot\big(2ab\bar{\mathbf{n}}+(a^{2}-c^{2})\bar{\mathrm{P}}-2bc\bar{\mathrm{Q}}\big).
\end{aligned}
\]
Now that $\mathbb{F}_{\bar{\mathbf{M}}}\bar{\mathbf{M}}$ is parallel
to $\bar{\mathbf{M}}$, $(2ab,a^{2}-c^{2},-2bc)$ is colinar with
$(a,b,c)$. Check that $\pm\bar{\mathbf{n}},\pm\bar{\mathrm{P}},\pm\bar{\mathrm{Q}}\in\mathrm{Det}_{0}^{\mathbb{F}}$
and thus $\pm\bar{\mathbf{n}},\pm\bar{\mathrm{P}},\pm\bar{\mathrm{Q}}\neq\pm\bar{\mathbf{M}}$.
Therefore, it holds $b\neq0$, since otherwise $\mathbb{F}_{\bar{\mathbf{M}}}\bar{\mathbf{M}}=(a^{2}-c^{2})\mathrm{P}$
and then $\bar{\mathbf{M}}=\pm\bar{\mathrm{P}}$. Consequently, $(2a,-2c)$
is colinear with $(a,c)$, which means either $a=0$ or $c=0$ holds.
As a result, 
\[
\begin{aligned}\mathbb{F}_{\bar{\mathbf{M}}}\bar{\mathbf{M}}= & \text{either} & ||\mathrm{P}||\cdot\big(2ab\bar{\mathbf{n}}+a^{2}\bar{\mathrm{P}}\big)\\
\text{} & \text{or} & -||\mathrm{P}||\cdot\big(c^{2}\bar{\mathrm{P}}+2bc\bar{\mathrm{Q}}\big)
\end{aligned}
.
\]
In the first case, $c=0$ and $a\neq0$, and then $\bar{\mathbf{M}}=a\bar{\mathbf{n}}+b\bar{\mathrm{P}}$
with $a^{2}+b^{2}=1$. That $(2ab,a^{2})=a\cdot(2b,a)$ and $(a,b)$
are colinear leads to $a^{2}-2b^{2}=0$, and hence $b=\pm\frac{1}{\sqrt{3}}$
and $a=\pm\sqrt{\frac{2}{3}}$. For the second case, it can be obtained
in a similar way that $(b,c)=(\pm\frac{1}{\sqrt{3}},\pm\sqrt{\frac{2}{3}})$.
It is then straightforward to check that
\[
\begin{aligned}\bar{\lambda}=||\mathbb{F}_{\bar{\mathbf{M}}}\bar{\mathbf{M}}||= & ||\mathrm{P}||\cdot||\pm\frac{2\sqrt{2}}{3}\bar{\mathbf{n}}+\frac{2}{3}\bar{\mathrm{P}}||\\
= & ||\mathrm{P}||\cdot||\frac{2}{3}\bar{\mathrm{P}}\pm\frac{2\sqrt{2}}{3}\bar{\mathrm{Q}}|| & =\frac{2}{\sqrt{3}}||\mathrm{P}||\\
 &  & =2|\bar{\mathbf{M}}\cdot\mathrm{P}|.
\end{aligned}
\]
Note that both $\bar{\mathbf{M}}_{0}=\pm\sqrt{\frac{2}{3}}\bar{\mathbf{n}}\pm\frac{1}{\sqrt{3}}\bar{\mathrm{P}}$
and $\bar{\mathbf{M}}_{1}=\pm\frac{1}{\sqrt{3}}\bar{\mathrm{P}}\pm\sqrt{\frac{2}{3}}\bar{\mathrm{Q}}$
are elements in $\overline{\mathcal{M}}$, and, $\bar{\mathbf{M}}_{0}\in\mathcal{A}$
and $\bar{\mathbf{M}}_{1}\in\mathcal{B}$. $\pm\frac{1}{\sqrt{3}}\bar{\mathrm{P}}$ then solves \eqref{eq:QuaEq-P}. Since $\bar{\mathbf{M}}_{1}\in\mathcal{P}$ and is neither collinear with nor perpendicular to $\pm\frac{1}{\sqrt{3}}\bar{\mathrm{P}}$, it is not a solution to \eqref{eq:QuaEq-P}, and thence $\bar{\mathbf{M}}_{1}\notin\mathcal{A}$. Meanwhile, $\bar{\mathbf{M}}_{0}\notin\mathcal{P}$ and hence $\bar{\mathbf{M}}_{0}\notin\mathcal{B}$.



While Theorem \ref{thm:F*F=uniform} shows that $\mathbb{F}^\intercal\mathbb{F}$ having identical eigenvalues uniquely determines the structure of $\mathbb{F}$, Lemma \ref{lem:WeakUniform-F*F} indicates that it suffices to know only part of $\mathbb{F}^\intercal\mathbb{F}$. 
It turns out that the piece of information about $\mathbb{F}^\intercal\mathbb{F}$ in Lemma \ref{lem:WeakUniform-F*F} can be replaced by another piece as shown in  
Theorem \ref{thm:F*F-semiUniform} below.
%
\begin{rem}
    $\mathbb{F}_{\bar{\mathrm{Q}}}$ has eigenvalues $0$ and $\pm\lambda_{\bar{\mathrm{Q}}}$. While both $\pm\lambda_{\bar{\mathrm{Q}}}$ are the principal eigenvalues, it is our convention to take $\lambda_{\bar{\mathrm{Q}}}$.
\end{rem}

\begin{thm}\label{thm:F*F-semiUniform}
    If $\mathcal{P}\subset\mathcal{E}_{\lambda_{\mathbb{F}}}$, then $\lambda_{\mathbb{F}}=2||\mathrm{P}||^2$, and consequently, all the eigenvalues of $\mathbb{F}^\intercal\mathbb{F}$ are identical.
\end{thm}
The proof of this theorem shares the core argument with that for Theorem \ref{thm:F*F=uniform}, which is formally stated as Lemma \ref{lem:tr[F_P][F_Q]=0} below:
\begin{lem}\label{lem:tr[F_P][F_Q]=0}
   If $\mathrm{tr}[\mathbb{F}_{\bar{\mathrm{Q}}}\mathbb{F}_{\bar{\mathrm{P}}}]=0$, then it holds with some $\lambda\in\mathbb{R}^3$ that 
    \begin{equation}\label{eq:tr[F_P][F_Q]=0}
\mathbb{F}_{\bar{\mathrm{Q}}}=\lambda\bigg(\bar{\mathrm{P}}\bar{\mathrm{Q}}^\intercal+\bar{\mathrm{Q}}\bar{\mathrm{P}}^\intercal\bigg)
\;\;\;\;\;\And\;\;\;\;
\mathbb{F}_{\bar{\mathrm{P}}}=||\mathrm{P}||\bar{\mathbf{n}}\bar{\mathbf{n}}^\intercal-(||\mathrm{P}||+\lambda)\bar{\mathrm{P}}\bar{\mathrm{P}}^\intercal+\lambda\bar{\mathrm{Q}}\bar{\mathrm{Q}}^\intercal.
    \end{equation}
\end{lem}
\begin{rem}
    $\lambda$ in \eqref{eq:tr[F_P][F_Q]=0} is in fact $\pm\lambda_{\bar{\mathrm{Q}}}$.
\end{rem}
\begin{proof}
Since $\mathbb{F}_{\bar{\mathrm{Q}}}\bar{\mathbf{n}}=\mathrm{P}\cdot\bar{\mathrm{Q}}=\mathbf{0}$, we have
\[
\bar{\mathrm{P}}^\intercal[\mathbb{F}_{\bar{\mathrm{Q}}}\mathbb{F}_{\bar{\mathrm{P}}}]\bar{\mathrm{P}}+\bar{\mathrm{Q}}^\intercal[\mathbb{F}_{\bar{\mathrm{Q}}}\mathbb{F}_{\bar{\mathrm{P}}}]\bar{\mathrm{Q}}=0.
\]
Now that $\mathbb{F}_{\bar{\mathrm{Q}}}$ has eigenvalues $\pm\lambda_{\bar{\mathrm{Q}}}$ on $\mathcal{P}$, it holds $\mathbb{F}^2_{\bar{\mathrm{Q}}}\mathbf{m}=\lambda^2_{\mathbf{Q}}\mathbf{m}$ for all $\mathbf{m}\in\mathcal{P}$, and thence
\[\bar{\mathrm{Q}}^\intercal\mathbb{F}^2_{\bar{\mathrm{P}}}\bar{\mathrm{P}}=\bar{\mathrm{P}}^\intercal[\mathbb{F}_{\bar{\mathrm{Q}}}\mathbb{F}_{\bar{\mathrm{P}}}]\bar{\mathrm{P}}=-\bar{\mathrm{Q}}^\intercal[\mathbb{F}_{\bar{\mathrm{Q}}}\mathbb{F}_{\bar{\mathrm{P}}}]\bar{\mathrm{Q}}=\bar{\mathrm{Q}}^\intercal\mathbb{F}^2_{\bar{\mathrm{Q}}}\bar{\mathrm{P}}=0.\]
This means $\mathbb{F}^2_{\bar{\mathrm{P}}}\bar{\mathrm{P}}$ to be perpendicular to $\bar{\mathrm{Q}}$, while both of them are in the two dimensional subspace $\mathcal{P}$. As a result, $\mathbb{F}^2_{\bar{\mathrm{P}}}\bar{\mathrm{P}}$ is collinear with $\bar{\mathrm{P}}$, that is, $\bar{\mathrm{P}}$ is an eigenvector of $\mathbb{F}^2_{\bar{\mathrm{P}}}$. By Theorem \ref{thm:StruturalThm-PlaneConfig}, $\bar{\mathbf{n}}$ is an eigenvector of $\mathbb{F}_{\bar{\mathrm{P}}}$ with eigenvalue $||\mathrm{P}||$ and hence also an eigenvector of $\mathbb{F}^2_{\bar{\mathrm{P}}}$  with $||\mathrm{P}||^2$. Now that $\mathbb{F}^2_{\bar{\mathrm{P}}}$ is symmetric with both $\bar{\mathbf{n}}$ and $\bar{\mathrm{P}}$ being its eigenvectors, $\bar{\mathrm{Q}}=\bar{\mathbf{n}}\times\bar{\mathrm{P}}$ is also an eigenvector of $\mathbb{F}^2_{\bar{\mathrm{P}}}$. We shall show that $\bar{\mathrm{P}},\bar{\mathrm{Q}}$ are in fact eigenvectors of $\mathbb{F}_{\bar{\mathrm{P}}}$, which implies $\mathbb{F}_{\bar{\mathrm{P}}}$ to take the form
\[\mathbb{F}_{\bar{\mathrm{P}}}=||\mathrm{P}||\bar{\mathbf{n}}\bar{\mathbf{n}}^\intercal+\lambda_0\bar{\mathrm{P}}\bar{\mathrm{P}}^\intercal+\lambda\bar{\mathrm{Q}}\bar{\mathrm{Q}}^\intercal
\]

Let $\bar{\mathbf{n}},\mathbf{m}_0,\mathbf{m}_1$ be orthonormal eigenvectors of  $\mathbb{F}_{\bar{\mathrm{P}}}$, and then $\mathbf{m}_0,\mathbf{m}_1\in\mathcal{P}$. Let $\lambda_0,\lambda$ be the eigenvalues of $\mathbb{F}_{\bar{\mathrm{P}}}$ associated to $\mathbf{m}_0,\mathbf{m}_1$, respectively. Since $\mathrm{tr}\mathbb{F}_{\bar{\mathrm{P}}}=0$, we have $\lambda_0+\lambda=-||\mathrm{P}||$.  $\lambda^2_0,\lambda^2$ are then eigenvalues of  $\mathbb{F}^2_{\bar{\mathrm{P}}}$, and,  $\mathbf{m}_0,\mathbf{m}_1$ are associated eigenvectors, respectively.
Therefore, either $\lambda_0=\lambda=-\frac{||\mathrm{P}||}{2}$, or $|\lambda_0|\neq|\lambda|$ . If $\lambda_0=\lambda=-\frac{||\mathrm{P}||}{2}$, then plane $\mathcal{P}$ is the corresponding eigen-space and hence all the vectors therein, including  $\bar{\mathrm{P}},\bar{\mathrm{Q}}$, are eigenvectors of  $\mathbb{F}_{\bar{\mathrm{P}}}$.  If $|\lambda_0|\neq|\lambda|$, then $\lambda^2_0,\lambda^2$ are distinct eigenvalues of $\mathbb{F}^2_{\bar{\mathrm{P}}}$ , and in particular, they are also the (distinct) eigenvalues of the restriction of  $\mathbb{F}^2_{\bar{\mathrm{P}}}$ to the invariant plane $\mathcal{P}$. This means that the corresponding eigen-spaces are both one dimensional. Since  $\bar{\mathrm{P}},\bar{\mathrm{Q}}$  are also eigenvectors of the restricted $\mathbb{F}^2_{\bar{\mathrm{P}}}$ , we know that either   $\bar{\mathrm{P}}=\pm\mathbf{m}_0,\bar{\mathrm{Q}}=\pm\mathbf{m}_1$, or $\bar{\mathrm{P}}=\pm\mathbf{m}_1,\bar{\mathrm{Q}}=\pm\mathbf{m}_0$.
Without loss of generality, we assume $\bar{\mathrm{P}}=\pm\mathbf{m}_0,\bar{\mathrm{Q}}=\pm\mathbf{m}_1$. As a result, $\mathbb{F}_{\bar{\mathrm{Q}}}\bar{\mathrm{P}}=\mathbb{F}_{\bar{\mathrm{P}}}\bar{\mathrm{Q}}=\lambda\bar{\mathrm{Q}}$. Since the restriction of $\mathbb{F}_{\bar{\mathrm{Q}}}$ to $\mathcal{P}$ is symmetry and has zero trace, its matrix representation w.r.t. the orthonormal basis $(\bar{\mathrm{P}},\bar{\mathrm{Q}})$ is  
$\left[\begin{array}{ccc}
0 & \lambda\\
\lambda & 0
\end{array}\right]$, which means $\mathbb{F}_{\bar{\mathrm{Q}}}=\lambda\bigg(\bar{\mathrm{P}}\bar{\mathrm{Q}}^\intercal+\bar{\mathrm{Q}}\bar{\mathrm{P}}^\intercal\bigg)$. The remaining coefficient $\lambda_0$ in \eqref{eq:tr[F_P][F_Q]=0} for $\mathbb{F}_{\bar{\mathrm{P}}}$ is then fixed through 
\[
\mathrm{tr}\mathbb{F}_{\bar{\mathrm{P}}}=||\mathrm{P}||+\lambda_0+\lambda=0.
\]
\end{proof}
\begin{rem}
    For the sake of rigor, we also want to make sure that magnetic algebras described in Lemma \ref{lem:tr[F_P][F_Q]=0} really exist.
    Let $(\bar{\mathbf{n}},\bar{\mathrm{P}},\bar{\mathrm{Q}})$ be an arbitrary orthonormal basis of $\mathbb{R}^3$, and, let $\mathbb{F}_{\bar{\mathrm{Q}}}$, $\mathbb{F}_{\bar{\mathrm{Q}}}$ be matrices given by \eqref{eq:tr[F_P][F_Q]=0} and $\mathbb{F}_{\bar{\mathbf{n}}}$ with \eqref{eq:[F_n]}. It is then straightforward to check that the algebra $(\mathbb{R}^3,\mathbb{F})$ with $\mathbb{F}$ defined as below is a magnetic algebra with invariant plane $\mathcal{P}=\mathrm{span}\{\bar{\mathbf{n}},\bar{\mathrm{P}}\}$, and also, $\mathrm{tr}[\mathbb{F}_{\bar{\mathrm{Q}}}\mathbb{F}_{\bar{\mathrm{P}}}]=0$:
    \[\mathbb{F}_{\mathbf{M}}:=(\mathbf{M}\cdot\bar{\mathbf{n}})\mathbb{F}_{\bar{\mathbf{n}}}+(\mathbf{M}\cdot\bar{\mathrm{P}})\mathbb{F}_{\bar{\mathrm{P}}}+(\mathbf{M}\cdot\bar{\mathrm{Q}})\mathbb{F}_{\bar{\mathrm{Q}}}.\]
\end{rem}

\textbf{Proof of Theorem \ref{thm:F*F-semiUniform}}: Since $\mathcal{E}_{\lambda_{\mathbb{F}}}=\mathcal{P}$, it holds $\mathrm{tr}[\mathbb{F}_{\bar{\mathrm{Q}}}\mathbb{F}_{\bar{\mathrm{P}}}]=\bar{\mathrm{Q}}^\intercal[\mathbb{F}^\intercal\mathbb{F}]\bar{\mathrm{P}}=0$. By Lemma \ref{lem:tr[F_P][F_Q]=0}, $\mathbb{F}_{\bar{\mathrm{P}}},\mathbb{F}_{\bar{\mathrm{Q}}}$ take their form in \eqref{eq:tr[F_P][F_Q]=0}, which implies $\lambda$ to be either $\lambda_{\bar{\mathrm{Q}}}$ or $-\lambda_{\bar{\mathrm{Q}}}$, and,  $||\mathrm{P}||+\lambda_0+\lambda=\mathrm{tr}\mathbb{F}_{\bar{\mathrm{P}}}=0$.


Note that $\mathrm{tr}\mathbb{F}^2_{\bar{\mathrm{Q}}}=2\lambda^2_{\bar{\mathrm{Q}}}=\lambda_{\mathbb{F}}$, and also,
\[\mathrm{tr}\mathbb{F}^2_{\bar{\mathrm{P}}}
=||\mathrm{P}||^2+\lambda^2_0+\lambda^2
=||\mathrm{P}||^2+\lambda^2_0+\lambda_{\bar{\mathrm{Q}}}^2
=\lambda_{\mathbb{F}}.
\]
Therefore, $||\mathrm{P}||^2+\lambda^2_0=\frac{\lambda_{\mathbb{F}}}{2}=\lambda^2_{\bar{\mathrm{Q}}}=\lambda^2$. On the other hand, we have $\lambda=-(||\mathrm{P}||+\lambda_0)$ since $||\mathrm{P}||+\lambda_0+\lambda=0$, and therefore, 
\[\lambda^2=||\mathrm{P}||^2+2\lambda_0||\mathrm{P}||+\lambda^2_0.\]
These together imply $\lambda_0||\mathrm{P}||=0$ and thus $\lambda_0=0$, and consequently, $\lambda=-||\mathrm{P}||$ and 
\[\lambda_{\mathbb{F}}=2\lambda_{\bar{\mathrm{Q}}}=2\lambda^2=2||\mathrm{P}||^2.\]

   
\section{Bounds for $\bar{\lambda}$}
\label{sec:GeometricBound}

Let $(\mathbb{R}^{3},\mathbb{F})$ be a magnetic algebra with an $\mathbb{F}$-invariant
plane $\mathcal{P}$, and $\bar{\mathbf{n}}$ be a normal direction
of $\mathcal{P}$. In this section we prove Theorem \ref{thm:=00005BGeometric-Bound=00005D},
which relates $\bar{\lambda}$ to $\lambda_{\mathcal{P}}$, the maximal
translational force with $\mathbf{M}\in\mathcal{P}$. To this end,
we break the theorem into smaller ones: Theorems \ref{thm:=00005Clambda_=00007B=00005CP=00007D>|=00005Clambda_=00007BM_F=00007D|>|P|},
\ref{thm:1stGeometricThm}, and \ref{thm:2ndGeometricThm}.

\subsection{$\bar{\lambda}$ and $\lambda_{\mathcal{P}}$}

Due to continuity and compactness, the maximum of $||\mathbb{F}_{\mathbf{M}}\mathbf{m}||$
for $\mathbf{M}\in S^{2}\cap\mathcal{P},\mathbf{m}\in S^{2}$ can
be achieved at some point $(\mathbf{M}_{\mathcal{P}},\tilde{\mathbf{m}})$,
and hence
\[
\underset{\mathbf{M}\in S^{2}\cap\mathcal{P},\mathbf{m}\in S^{2}}{\max}||\mathbb{F}_{\mathbf{M}}\mathbf{m}||=||\mathbb{F}_{\mathbf{M}_{\mathcal{P}}}\tilde{\mathbf{m}}||=\max_{\mathbf{m}\in S^{2}}||\mathbb{F}_{\mathbf{M}_{\mathcal{P}}}\mathbf{m}||=|\lambda_{\mathbf{M}_{\mathcal{P}}}|.
\]
As a consequence, 
\[
|\lambda_{\mathbf{M}_{\mathcal{P}}}|=\underset{\mathbf{M}\in S^{2}\cap\mathcal{P}}{\max}\lambda_{\mathbf{M}}=:\lambda_{\mathcal{P}}.
\]
Note that $||\mathrm{P}||$ is an eigenvalue of $\mathbb{F}_{\mathbf{M}_{\bar{\mathrm{P}}}}$
with $\bar{\mathrm{P}}=\frac{\mathrm{P}}{||\mathrm{P}||}$ and hence
$\lambda_{\mathcal{P}}\geq||\mathrm{P}||$. Meanwhile, by Theorem
\ref{thm:=00005BStructuralThm-PlaneConfig=00005D=00005Cn-eigenVectorF*F},
$\bar{\mathbf{n}}$ is always an eigenvector of $\mathbb{F}^{\intercal}\mathbb{F}$,
and hence $\mathbf{M}_{\mathbb{F}}$ can be always taken from the
subset $\{\bar{\mathbf{n}}\}\cup\mathcal{P}$. If $\mathbf{M}_{\mathbb{F}}=\bar{\mathbf{n}}$,
then $|\lambda_{\mathbf{M}_{\mathbb{F}}}|=|\lambda_{\bar{\mathbf{n}}}|=||\mathrm{P}||\leq\lambda_{\mathcal{P}}$.
If $\mathbf{M}_{\mathbb{F}}\in\mathcal{P}$, then $\lambda_{\mathcal{P}}\geq|\lambda_{\mathbf{M}_{\mathbb{F}}}|$
just holds by definition (of $\lambda_{\mathcal{P}}$). This means
\begin{thm}
\label{thm:=00005Clambda_=00007B=00005CP=00007D>|=00005Clambda_=00007BM_F=00007D|>|P|}For
a magentic algebra $(\mathbb{R}^{3},\mathbb{F})$ with planarity,
it always holds
\begin{equation}
\lambda_{\mathcal{P}}\geq|\lambda_{\mathbf{M}_{\mathbb{F}}}|\geq||\mathrm{P}||.\label{eq:=00005Clambda_=00005CP>|=00005Clambda_=00007BM_F=00007D|>|P|}
\end{equation}
\end{thm}
For a general pair $(\mathbf{M},\mathbf{m})\in S^{2}\times S^{2}$,
consider the orthogonal decomposition
\[
\mathbf{M}=\cos\beta\bar{\mathbf{n}}+\sin\beta\mathbf{M}_{\perp}\ \ \text{ and }\ \ \mathbf{m}=\cos\alpha\bar{\mathbf{n}}+\sin\alpha\mathbf{m}_{\perp}
\]
where $\mathbf{M}_{\perp},\mathbf{m}_{\perp}\in\mathcal{P}\cap S^{2}$.
By linearity it yields
\begin{equation}
\mathbb{F}_{\mathbf{M}}\mathbf{m}=\cos\beta\cos\alpha\mathbb{F}_{\bar{\mathbf{n}}}\bar{\mathbf{n}}+\cos\beta\sin\alpha\mathbb{F}_{\bar{\mathbf{n}}}\mathbf{m}_{\perp}+\sin\beta\cos\alpha\mathbb{F}_{\mathbf{M}_{\perp}}\bar{\mathbf{n}}+\sin\beta\sin\alpha\mathbb{F}_{\mathbf{M}_{\perp}}\mathbf{m}_{\perp}.\label{eq:OrthogonalDecomposition-=00005BF_M=00005Dm}
\end{equation}

Due to the planar configuration,
\begin{equation}
\mathbb{F}_{\bar{\mathbf{n}}}=\bar{\mathbf{n}}\mathrm{P}^{\intercal}+\mathrm{P}\bar{\mathbf{n}}^{\intercal},\label{eq:=00005CF_=00005Cn}
\end{equation}
and then $\mathbb{F}_{\bar{\mathbf{n}}}\bar{\mathbf{n}}=\mathrm{P}$,
$\mathbb{F}_{\bar{\mathbf{n}}}\mathbf{m}_{\perp}=(\mathrm{P}\cdot\mathbf{m}_{\perp})\bar{\mathbf{n}}$,
and,
\[
\mathbb{F}_{\mathbf{M}_{\perp}}\bar{\mathbf{n}}=\mathbb{F}_{\bar{\mathbf{n}}}\mathbf{M}_{\perp}=(\mathrm{P}\cdot\mathbf{M}_{\perp})\bar{\mathbf{n}}.
\]
The sum in (\ref{eq:OrthogonalDecomposition-=00005BF_M=00005Dm})
then becomes ($\mathbb{F}_{\mathbf{M}}\mathbf{m}=$)
\[
\cos\beta\cos\alpha\cdot\mathrm{P}+\cos\beta\sin\alpha(\mathrm{P}\cdot\mathbf{m}_{\perp})\bar{\mathbf{n}}+\sin\beta\cos\alpha(\mathrm{P}\cdot\mathbf{M}_{\perp})\bar{\mathbf{n}}+\sin\beta\sin\alpha\mathbb{F}_{\mathbf{M}_{\perp}}\mathbf{m}_{\perp},
\]
yielding the orthogonal decomposition $\mathbb{F}_{\mathbf{M}}\mathbf{m}=\big(\mathbb{F}_{\mathbf{M}}\mathbf{m}\big)_{\mathcal{P}}+\big(\mathbb{F}_{\mathbf{M}}\mathbf{m}\big)_{\bar{\mathbf{n}}}$
\[
\big(\mathbb{F}_{\mathbf{M}}\mathbf{m}\big)_{\mathcal{P}}=\cos\beta\cos\alpha\cdot\mathrm{P}+\sin\beta\sin\alpha\mathbb{F}_{\mathbf{M}_{\perp}}\mathbf{m}_{\perp}
\]
and
\[
\big(\mathbb{F}_{\mathbf{M}}\mathbf{m}\big)_{\bar{\mathbf{n}}}=\cos\beta\sin\alpha(\mathrm{P}\cdot\mathbf{m}_{\perp})\bar{\mathbf{n}}+\sin\beta\cos\alpha(\mathrm{P}\cdot\mathbf{M}_{\perp})\bar{\mathbf{n}}.
\]
As a result, $\big|\big|\mathbb{F}_{\mathbf{M}}\mathbf{m}\big|\big|^{2}=\big|\big|\big(\mathbb{F}_{\mathbf{M}}\mathbf{m}\big)_{\mathcal{P}}\big|\big|^{2}+\big|\big|\big(\mathbb{F}_{\mathbf{M}}\mathbf{m}\big)_{\bar{\mathbf{n}}}\big|\big|^{2}$
and then
\[
\begin{aligned}\big|\big|\mathbb{F}_{\mathbf{M}}\mathbf{m}\big|\big|^{2} \leq & ||\mathrm{P}||^{2}\big(\cos\beta\sin\alpha+\sin\beta\cos\alpha\big)^{2}+\\&\big(|\sin\beta\sin\alpha|\cdot|\lambda_{\mathbf{M}_{\perp}}|+|\cos\beta\cos\alpha|\cdot||\mathrm{P}||\big)^{2}\\
 \leq & ||\mathrm{P}||^{2}\big(1-\sin^{2}\beta\sin^{2}\alpha)+\frac{|\sin2\alpha\sin2\beta|}{2}||\mathrm{P}||\big(||\mathrm{P}||+|\lambda_{\mathbf{M}_{\perp}}|\big)+\\&\sin^{2}\beta\sin^{2}\alpha\cdot|\lambda_{\mathbf{M}_{\perp}}|^{2}\\
 = & (1-x^{2}y^{2})||\mathrm{P}||^{2}+2xy\sqrt{1-x^{2}}\sqrt{1-y^{2}}||\mathrm{P}||\big(||\mathrm{P}||+|\lambda_{\mathbf{M}_{\perp}}|\big)+\\&x^{2}y^{2}\cdot|\lambda_{\mathbf{M}_{\perp}}|^{2}\\
 \leq & (1-x^{2}y^{2})||\mathrm{P}||^{2}+2xy(1-xy)||\mathrm{P}||\big(||\mathrm{P}||+|\lambda_{\mathbf{M}_{\perp}}|\big)+x^{2}y^{2}\cdot|\lambda_{\mathbf{M}_{\perp}}|^{2}
\end{aligned}
\]
where $(x,y)=(|\sin\alpha|,|\sin\beta|)$. Let $t\in xy$ and $\lambda_{\mathcal{P}}=\underset{\mathbf{M}\in\mathcal{P}\cap S^{2}}{\max}|\lambda_{\mathbf{M}}|$.
From the inequation above we get

\[
\begin{aligned}\max_{\mathbf{M},\mathbf{m}\in S^{2}}\big|\big|\mathbb{F}_{\mathbf{M}}\mathbf{m}\big|\big|^{2}\leq & \max_{t\in[0,1]}(1-t^{2})||\mathrm{P}||^{2}+2t(1-t)||\mathrm{P}||\cdot\big(||\mathrm{P}||+\lambda_{\mathcal{P}}\big)+t^{2}\cdot\lambda_{\mathcal{P}}^{2}.\\
 \stackrel{i.e.}{=} & \max_{t\in[0,1]}||\mathrm{P}||^{2}+2||\mathrm{P}||\big(||\mathrm{P}||+\lambda_{\mathcal{P}}\big)t+(\lambda_{\mathcal{P}}+||\mathrm{P}||)(\lambda_{\mathcal{P}}-3||\mathrm{P}||)t^{2}.
\end{aligned}
\]
Let $g(t)=||\mathrm{P}||^{2}+2||\mathrm{P}||\big(||\mathrm{P}||+\lambda_{\mathcal{P}}\big)t+(\lambda_{\mathcal{P}}+||\mathrm{P}||)(\lambda_{\mathcal{P}}-3||\mathrm{P}||)t^{2}$.
Then $g(t)$ takes its maximum on $[0,1]$ at either $t=0\text{ or }1$,
or, $t=\frac{||\mathrm{P}||}{3||\mathrm{P}||-\lambda_{\mathcal{P}}}$
when
\[
\frac{||\mathrm{P}||}{3||\mathrm{P}||-\lambda_{\mathcal{P}}}\in(0,1).
\]
Check that $g(0)=||\mathrm{P}||^{2}$ and $g(1)=\lambda_{\mathcal{P}}^{2}.$
When $\frac{||\mathrm{P}||}{3||\mathrm{P}||-\lambda_{\mathcal{P}}}\in(0,1)$,
i.e., $\lambda_{\mathcal{P}}<2||\mathrm{P}||$, it holds 
\[
\max_{t\in[0,1]}g(t)=g(\frac{||\mathrm{P}||}{3||\mathrm{P}||-\lambda_{\mathcal{P}}})=\frac{4||\mathrm{P}||^{2}\cdot||\mathrm{P}||}{3||\mathrm{P}||-\lambda_{\mathcal{P}}}\leq4||\mathrm{P}||^{2}.
\]

These results are summarized in the following theorem:
\begin{thm}
\label{thm:1stGeometricThm}If $\lambda_{\mathcal{P}}\geq2||\mathrm{P}||$,
then $\bar{\lambda}=\lambda_{\mathcal{P}}=|\lambda_{\mathbf{M}_{\mathcal{P}}}|$
and hence $\bar{\mathbf{M}}$ can be taken as $\mathbf{M}_{\mathcal{P}}\in\mathcal{P}$.
If $\lambda_{\mathcal{P}}<2||\mathrm{P}||$, then
\[
\bar{\lambda}\leq2||\mathrm{P}||\cdot\sqrt{\frac{||\mathrm{P}||}{3||\mathrm{P}||-\lambda_{\mathcal{P}}}}.
\]
\end{thm}
\begin{rem}
From Theorem \ref{thm:MaximalLargestEigenvalue-all=00005CF_M}, we
can also obtain part of Theorem \ref{thm:1stGeometricThm} : when
$\bar{\lambda}>2||\mathrm{P}||$, $\bar{\lambda}=\lambda_{\mathcal{P}}$.
According to Theorem \ref{thm:MaximalLargestEigenvalue-all=00005CF_M},
either $\bar{\lambda}=\lambda_{\mathbf{M}_{\mathbb{F}}}$, or $\bar{\lambda}=2\bar{\mathbf{M}}\cdot\mathrm{P}\leq2||\mathrm{P}||$,
or $\bar{\mathbf{M}}\subset\mathcal{P}$. So, if $\bar{\lambda}>2||\mathrm{P}||$, then either $\bar{\lambda}=\lambda_{\mathbf{M}_{\mathbb{F}}}$ 
or $\bar{\mathbf{M}}\subset\mathcal{P}$ holds. Furthermore, if $\bar{\lambda}=\lambda_{\mathbf{M}_{\mathbb{F}}}$,
either $\mathbf{M}_{\mathbb{F}}=\bar{\mathbf{n}}$ and then $\bar{\lambda}=\lambda_{\bar{\mathbf{n}}}=||\mathrm{P}||$,
or, $\mathbf{M}_{\mathbb{F}}$ can be taken from $\mathcal{P}$.
\end{rem}
\vspace{0.2cm}
Theorem \ref{thm:1stGeometricThm} indicates to bound $\lambda_{\mathcal{P}}=\underset{\mathbf{M}\in\mathcal{P}\cap S^{2}}{\max}|\lambda_{\mathbf{M}}|$. The corollary below follows
directly from Theorem \ref{thm:1stGeometricThm}, which is convenient
to use for latter discussion.
\begin{cor}
\label{cor:1stGeometricCorollary}$\bar{\lambda}>2||\mathrm{P}||\iff\lambda_{\mathcal{P}}>2||\mathrm{P}||$,
$\bar{\lambda}=2||\mathrm{P}||\iff\lambda_{\mathcal{P}}=2||\mathrm{P}||$,
and $\bar{\lambda}<2||\mathrm{P}||\iff\lambda_{\mathcal{P}}<2||\mathrm{P}||$
\end{cor}
\vspace{0.2cm}
\subsection{Bounds for $\lambda_{\mathcal{P}}=\protect\underset{\mathbf{M}\in\mathcal{P}\cap S^{2}}{\max}|\lambda_{\mathbf{M}}|$ }

Consider $\mathbf{M}\in\mathcal{P}$, according to Theorem \ref{thm:StruturalThm-PlaneConfig},
$\mathrm{P}\cdot\mathbf{M}$ is an eigenvalue of $\mathbb{F}_{\mathbf{M}}$
and hence it equals to one of $\lambda_{\mathbf{M}}$ and $-\frac{\lambda_{\mathbf{M}}}{2}\pm\delta_{\mathbf{M}}$,
that is,
\[
\text{either }\ \lambda_{\mathbf{M}}=\mathrm{P}\cdot\mathbf{M}\ \text{ or \ }\frac{|\lambda_{\mathbf{M}}|}{2}(1\pm r_{\mathbf{M}})=|\mathrm{P}\cdot\mathbf{M}|.
\]
Due to Inequation \ref{eq:=00005Clambda_=00005CP>|=00005Clambda_=00007BM_F=00007D|>|P|},
for studying $\lambda_{\mathcal{P}}=|\lambda_{\mathbf{M}_{\mathcal{P}}}|$
we only need to focus on those $\mathbf{M}\in\mathcal{P}\cap S^{2}$
with $|\lambda_{\mathbf{M}}|>||\mathrm{P}||$. For such an $\mathbf{M}$,
we have 
\[
\frac{|\lambda_{\mathbf{M}}|}{2}(1\pm r_{\mathbf{M}})=|\mathrm{P}\cdot\mathbf{M}|
\]
and hence
\[
r_{\mathbf{M}}^{2}=\big(\frac{2|\mathrm{P}\cdot\mathbf{M}|}{|\lambda_{\mathbf{M}}|}-1\big)^{2},
\]
plugging which into (\ref{eq:=00005Clamda^2=000026|F_M|^2}) yields
the equation
\[
\text{tr}\mathbb{F}_{\mathbf{M}}^{2}=\frac{3+\big(\frac{2|\mathrm{P}\cdot\mathbf{M}|}{|\lambda_{\mathbf{M}}|}-1\big)^{2}}{2}\lambda_{\mathbf{M}}^{2}=\frac{3|\lambda_{\mathbf{M}}|^{2}+\big(2|\mathrm{P}\cdot\mathbf{M}|-|\lambda_{\mathbf{M}}|\big)^{2}}{2},
\]
i.e., $t=|\lambda_{\mathbf{M}}|$ satisfies the quadratic equation
$t^{2}-|\mathrm{P}\cdot\mathbf{M}|t+|\mathrm{P}\cdot\mathbf{M}|^{2}=\frac{\text{tr}\mathbb{F}_{\mathbf{M}}^{2}}{2}$,
solving which returns
\[
\frac{|\mathrm{P}\cdot\mathbf{M}|\pm\sqrt{2\text{tr}\mathbb{F}_{\mathbf{M}}^{2}-3|\mathrm{P}\cdot\mathbf{M}|^{2}}}{2}.
\]

For an arbitrary $\mathbf{M}\in\mathcal{P}$, $|\lambda_{\mathbf{M}}|\geq|\mathrm{P}\cdot\mathbf{M}|$
since $\mathrm{P}\cdot\mathbf{M}$ is an eigenvalue, and therefore,
\begin{equation}
|\lambda_{\mathbf{M}}|=\text{either }\ \frac{|\mathrm{P}\cdot\mathbf{M}|+\sqrt{2\text{tr}\mathbb{F}_{\mathbf{M}}^{2}-3|\mathrm{P}\cdot\mathbf{M}|^{2}}}{2}\ \text{ or }\ |\mathrm{P}\cdot\mathbf{M}|,\ \ \forall\mathbf{M}\in\mathcal{P}.\label{eq:=00005Clambda-on-PLane}
\end{equation}
The theorem below gives a slight more accurate statement about this
issue:
\begin{thm}
\label{thm:|=00005Clambda|-on-PLane}For any $\mathbf{M}\in\mathcal{P}$,
\begin{equation}
|\lambda_{\mathbf{M}}|=\max\bigg\{\frac{|\mathrm{P}\cdot\mathbf{M}|+\sqrt{2\text{tr}\mathbb{F}_{\mathbf{M}}^{2}-3|\mathrm{P}\cdot\mathbf{M}|^{2}}}{2},|\mathrm{P}\cdot\mathbf{M}|\bigg\}.\label{eq:|=00005Clambda|-on-PL=00003Dmax=00007B=00007D}
\end{equation}
\end{thm}
\begin{proof}
The discussion above has shown that if $-\frac{\lambda_{\mathbf{M}}}{2}\pm\delta_{\mathbf{M}}=\mathrm{P}\cdot\mathbf{M}$,
then 
\[
|\lambda_{\mathbf{M}}|=\frac{|\mathrm{P}\cdot\mathbf{M}|+\sqrt{2\text{tr}\mathbb{F}_{\mathbf{M}}^{2}-3|\mathrm{P}\cdot\mathbf{M}|^{2}}}{2}\geq|\mathrm{P}\cdot\mathbf{M}|.
\]
For the case $\lambda_{\mathbf{M}}=\mathrm{P}\cdot\mathbf{M}$, we
resort to (\ref{eq:BasicInequation-=00005Clambda_M}):
\[
\begin{aligned}
\frac{|\mathrm{P}\cdot\mathbf{M}|+\sqrt{2\text{tr}\mathbb{F}_{\mathbf{M}}^{2}-3|\mathrm{P}\cdot\mathbf{M}|^{2}}}{2}&\leq\frac{|\mathrm{P}\cdot\mathbf{M}|+\sqrt{4|\mathrm{P}\cdot\mathbf{M}|^{2}-3|\mathrm{P}\cdot\mathbf{M}|^{2}}}{2}\\ &=|\mathrm{P}\cdot\mathbf{M}|=|\lambda_{\mathbf{M}}|. 
\end{aligned}
\]

\end{proof}
Since $\lambda_{\mathcal{P}}\geq||\mathrm{P}||$ and is achieved at
some $\mathbf{M}_{\mathcal{P}}\in\mathcal{P}\cap S^{2}$, Theorem
\ref{thm:|=00005Clambda|-on-PLane} implies
\begin{equation}
\lambda_{\mathcal{P}}=\max\bigg\{\max_{\mathbf{M}\in\mathcal{P}\cap S^{2}}\frac{|\mathrm{P}\cdot\mathbf{M}|+\sqrt{2\text{tr}\mathbb{F}_{\mathbf{M}}^{2}-3|\mathrm{P}\cdot\mathbf{M}|^{2}}}{2},||\mathrm{P}||\bigg\}.\label{eq:Formula=00005B=00005Clambda_=00007B=00005CP=00007D=00005D}
\end{equation}
Now that $\text{tr}\mathbb{F}_{\mathbf{M}_{\mathbb{F}}}^{2}\geq\text{tr}\mathbb{F}_{\bar{\mathbf{n}}}^{2}=2||\mathrm{P}||^{2}$,
\[
\begin{aligned}
   &\max\bigg\{\max_{\mathbf{M}\in\mathcal{P}\cap S^{2}}\frac{|\mathrm{P}\cdot\mathbf{M}|+\sqrt{2\text{tr}\mathbb{F}_{\mathbf{M}}^{2}-3|\mathrm{P}\cdot\mathbf{M}|^{2}}}{2},||\mathrm{P}||\bigg\}\\
\leq&\frac{|\mathrm{P}\cdot\mathbf{M}|+\sqrt{2\text{tr}\mathbb{F}_{\mathbf{M}_{\mathbb{F}}}^{2}-3|\mathrm{P}\cdot\mathbf{M}|^{2}}}{2},
\end{aligned}
\]
As a result, for an upper bound to $\lambda_{\mathcal{P}}$, it suffices
to bound \[h(t)=\frac{t||\mathrm{P}||+\sqrt{2\text{tr}\mathbb{F}_{\mathbf{M}_{\mathbb{F}}}^{2}-3t^{2}||\mathrm{P}||^{2}}}{2}\]
for $t\in[0,1]$. Note that the function $h$ is well defined for all $t\in[-\sqrt{\frac{2\text{tr}\mathbb{F}_{\mathbf{M}_{\mathbb{F}}}^{2}}{3\cdot||\mathrm{P}||^{2}}},\sqrt{\frac{2\text{tr}\mathbb{F}_{\mathbf{M}_{\mathbb{F}}}^{2}}{3\cdot||\mathrm{P}||^{2}}}]$.
Solving $\frac{dh}{dt}=0$ returns
\[
t^{2}=\frac{\text{tr}\mathbb{F}_{\mathbf{M}_{\mathbb{F}}}^{2}}{6\cdot||\mathrm{P}||^{2}}\leq\frac{2\text{tr}\mathbb{F}_{\mathbf{M}_{\mathbb{F}}}^{2}}{3\cdot||\mathrm{P}||^{2}}.
\]
Also, check that $t^{2}<\frac{\text{tr}\mathbb{F}_{\mathbf{M}_{\mathbb{F}}}^{2}}{6\cdot||\mathrm{P}||^{2}}$
implies $\frac{dh}{dt}>0$, and $t^{2}>\frac{\text{tr}\mathbb{F}_{\mathbf{M}_{\mathbb{F}}}^{2}}{6\cdot||\mathrm{P}||^{2}}$
implies $\frac{dh}{dt}<0$, which means $t=\sqrt{\frac{\text{tr}\mathbb{F}_{\mathbf{M}_{\mathbb{F}}}^{2}}{6\cdot||\mathrm{P}||^{2}}}$
is the maximum point of $h$. Moreover, the maximum
\[
\begin{aligned}h\bigg(\sqrt{\frac{\text{tr}\mathbb{F}_{\mathbf{M}_{\mathbb{F}}}^{2}}{6\cdot||\mathrm{P}||^{2}}}\bigg)= & \frac{\sqrt{\frac{\text{tr}\mathbb{F}_{\mathbf{M}_{\mathbb{F}}}^{2}}{6\cdot||\mathrm{P}||^{2}}}||\mathrm{P}||+\sqrt{2\text{tr}\mathbb{F}_{\mathbf{M}_{\mathbb{F}}}^{2}-3\frac{\text{tr}\mathbb{F}_{\mathbf{M}_{\mathbb{F}}}^{2}}{6\cdot||\mathrm{P}||^{2}}||\mathrm{P}||^{2}}}{2}\\
= & \sqrt{\frac{2}{3}\text{tr}\mathbb{F}_{\mathbf{M}_{\mathbb{F}}}^{2}},
\end{aligned}
\]
which exactly equals the upper bound given in (\ref{eq:(1st estimate)<=00005Cbar=00007B=00005Clambda=00007D<}).
Also, $h(0)=\sqrt{\frac{\text{tr}\mathbb{F}_{\mathbf{M}_{\mathbb{F}}}^{2}}{2}}\geq\sqrt{\frac{\text{tr}\mathbb{F}_{\bar{\mathbf{n}}}^{2}}{2}}\geq||\mathrm{P}||^{2}.$
As a result,
\[
\max_{\mathbf{M}\in S^{2}\cap\mathcal{P}}|\lambda_{\mathbf{M}}|\leq\max_{t\in[0,1]}h(t)\leq\sqrt{\frac{2}{3}\text{tr}\mathbb{F}_{\mathbf{M}_{\mathbb{F}}}^{2}}.
\]
Furthermore, in the case $\frac{\text{tr}\mathbb{F}_{\mathbf{M}_{\mathbb{F}}}^{2}}{6\cdot||\mathrm{P}||^{2}}>1$
( i.e. $\text{tr}\mathbb{F}_{\mathbf{M}_{\mathbb{F}}}^{2}>6\cdot||\mathrm{P}||^{2}$),
since $\underset{t\in[0,1]}{\max}h(t)=h(1)<h(\sqrt{\frac{\text{tr}\mathbb{F}_{\mathbf{M}_{\mathbb{F}}}^{2}}{6\cdot||\mathrm{P}||^{2}}})$,
the upper bound is strictly shaper than $\sqrt{\frac{2}{3}\text{tr}\mathbb{F}_{\mathbf{M}_{\mathbb{F}}}^{2}}$:
\begin{equation}
\max_{\mathbf{M}\in S^{2}\cap\mathcal{P}}|\lambda_{\mathbf{M}}|\leq h(1)=\frac{||\mathrm{P}||+\sqrt{2\text{tr}\mathbb{F}_{\mathbf{M}_{\mathbb{F}}}^{2}-3||\mathrm{P}||^{2}}}{2}<\sqrt{\frac{2}{3}\text{tr}\mathbb{F}_{\mathbf{M}_{\mathbb{F}}}^{2}}.\label{Ineq:UpperBound-onPLane}
\end{equation}

Note that when $\bar{\lambda}>2||\mathrm{P}||$, (\ref{eq:(1st estimate)<=00005Cbar=00007B=00005Clambda=00007D<})
implies $\text{tr}\mathbb{F}_{\mathbf{M}_{\mathbb{F}}}^{2}>6\cdot||\mathrm{P}||^{2}$.
Combining the conclusion above with Theorem \ref{thm:1stGeometricThm}
yields
\begin{thm}
\label{thm:2ndGeometricThm}When $\text{tr}\mathbb{F}_{\mathbf{M}_{\mathbb{F}}}^{2}>6\cdot||\mathrm{P}||^{2}$,
$\lambda_{\mathcal{P}}:=\underset{\mathbf{M}\in S^{2}\cap\mathcal{P}}{\max}\lambda_{\mathbf{M}}$
has an upper bound $\frac{||\mathrm{P}||+\sqrt{2\text{tr}\mathbb{F}_{\mathbf{M}_{\mathbb{F}}}^{2}-3||\mathrm{P}||^{2}}}{2}$.
In particular, if $\lambda_{\mathcal{P}}\geq2||\mathrm{P}||$, then
\[
\bar{\lambda}=\lambda_{\mathcal{P}}<\frac{||\mathrm{P}||+\sqrt{2\text{tr}\mathbb{F}_{\mathbf{M}_{\mathbb{F}}}^{2}-3||\mathrm{P}||^{2}}}{2}<\sqrt{\frac{2}{3}\text{tr}\mathbb{F}_{\mathbf{M}_{\mathbb{F}}}^{2}}.
\]
\end{thm}

\subsection{Proof of Theorem \ref{thm:=00005BGeometric-Bound=00005D}}

The lower bounds to $\bar{\lambda}$ in \ref{eq:BasicGeoInequation=00005BwithPlanarity=00005D}
have been given by Theorem \ref{thm:=00005Clambda_=00007B=00005CP=00007D>|=00005Clambda_=00007BM_F=00007D|>|P|}. 

For the upper bound in \ref{eq:BasicGeoInequation=00005BwithPlanarity=00005D},
it suffices to prove (\ref{eq:RefinedUpperBound-=00005B=00005Clambda<2|P|=00005D})
and (\ref{eq:RefinedUpperBound-=00005B=00005Clambda>2|P|=00005D})
directly. This is because, 
\[
|\lambda_{\mathbf{M}_{\mathbb{F}}}|+\frac{||\mathrm{P}||}{3}\leq|\lambda_{\mathbf{M}_{\mathbb{F}}}|+\frac{||\mathrm{P}||}{2},
\]
 and, by Inequation (\ref{eq:BasicInequation-=00005Clambda_M}),
\[
\frac{||\mathrm{P}||+\sqrt{2\lambda_{\mathbb{F}}-3||\mathrm{P}||^{2}}}{2}\leq\sqrt{\frac{\text{tr}\mathbb{F}_{\mathbf{M}_{\mathbb{F}}}^{2}}{2}}+\frac{||\mathrm{P}||}{2}\leq|\lambda_{\mathbf{M}_{\mathbb{F}}}|+\frac{||\mathrm{P}||}{2}.
\]

The refined upper bound for the case $\lambda_{\mathcal{P}}\geq2||\mathrm{P}||$
is simply given by Theorem \ref{thm:2ndGeometricThm}.

For (\ref{eq:RefinedUpperBound-=00005B=00005Clambda<2|P|=00005D}),
note that when $\lambda_{\mathcal{P}}\leq2||\mathrm{P}||$, (\ref{eq:=00005Clambda_=00005CP>|=00005Clambda_=00007BM_F=00007D|>|P|})
yields $|\lambda_{\mathbf{M}_{\mathbb{F}}}|\leq\lambda_{\mathcal{P}}\leq2||\mathrm{P}||$,
and then Inequation (\ref{Ineq:1stErrEsti-GenericPosition}) leads
to 
\[
\bar{\lambda}\leq|\lambda_{\mathbf{M}_{\mathbb{F}}}|+\frac{r_{\mathbf{M}_{\mathbb{F}}}^{2}}{6}|\lambda_{\mathbf{M}_{\mathbb{F}}}|\leq|\lambda_{\mathbf{M}_{\mathbb{F}}}|+\frac{||\mathrm{P}||}{3},
\]
combining which with Theorem \ref{thm:1stGeometricThm} gives
\[
\bar{\lambda}\leq\min\bigg\{|\lambda_{\mathbf{M}_{\mathbb{F}}}|+\frac{||\mathrm{P}||}{3},2||\mathrm{P}||\cdot\sqrt{\frac{||\mathrm{P}||}{3||\mathrm{P}||-\lambda_{\mathcal{P}}}}\bigg\}.
\]

\pagebreak{}

\section*{Acknowledgements}
This work has been supported by Research Development Fund of Xi\text{'}an Jiaotong-Liverpool University (RDF-24-01-086),
Wenzhou Science and Technology Major Project (ZG2024005), Wenzhou Municipal Key Science and Research Program
(ZG2023023) and Wenling Municipal Key Science and Research Program (2023G00009).

\vspace{0.6cm}

{\large\textbf{Conflict of Interests}}\\
The author declares no conflict of interests.

\vspace{0.2cm}

{\large\textbf{Data Availability Statement}}\\
Data sharing is not applicable to this article as no data sets were generated or analyzed during the current study.

\vspace{2cm}


\bibliographystyle{plain}
\bibliography{MagAlgebra}

\pagebreak{}

\section*{Appendix}

\subsection*{A. Formula for $|\lambda_{\mathbf{M}_{\mathbb{F}}}|$ }

Given an abstract magnetic algebra $(\mathbb{R}^{3},\mathbb{F})$
with invariant plane $\mathcal{P}$, we already know that $\mathbf{M}_{\mathbb{F}}$
can always be taken from the subset $\{\pm\bar{\mathbf{n}}\}\cup\mathcal{P}$. 

In the case $\mathbf{M}_{\mathbb{F}}\in\mathcal{P}$, $\mathbf{M}_{\mathbb{F}}\cdot\mathrm{P}$
is always an eigenvector of $\mathbb{F}_{\mathbf{M}_{\mathbb{F}}}$,
and therefore,
\[
\text{either }\ |\lambda_{\mathbf{M}_{\mathbb{F}}}|=|\mathbf{M}_{\mathbb{F}}\cdot\mathrm{P}|\ \text{ or }\ \frac{|\lambda_{\mathbf{M}_{\mathbb{F}}}|}{2}\pm\delta_{\mathbf{M}_{\mathbb{F}}}=|\mathbf{M}_{\mathbb{F}}\cdot\mathrm{P}|.
\]

If $\frac{|\lambda_{\mathbf{M}_{\mathbb{F}}}|}{2}\pm\delta_{\mathbf{M}_{\mathbb{F}}}=|\mathbf{M}_{\mathbb{F}}\cdot\mathrm{P}|$,
from (\ref{eq:=00005Clambda-on-PLane}) it holds
\begin{equation}
|\lambda_{\mathbf{M}_{\mathbb{F}}}|=\frac{|\mathrm{P}\cdot\mathbf{M}_{\mathbb{F}}|+\sqrt{2\text{tr}\mathbb{F}_{\mathbf{M}_{\mathbb{F}}}^{2}-3|\mathrm{P}\cdot\mathbf{M}_{\mathbb{F}}|^{2}}}{2}.\label{eq:=00005Clambda_MF-onPLane}
\end{equation}

It turns out that (\ref{eq:=00005Clambda_MF-onPLane}) also holds
when $|\lambda_{\mathbf{M}_{\mathbb{F}}}|=|\mathbf{M}_{\mathbb{F}}\cdot\mathrm{P}|$.
In this case, (\ref{Ineq:LowBound-=00005CLambda_=00005CMF}) implies
\[
|\lambda_{\mathbf{M}_{\mathbb{F}}}|=|\mathbf{M}_{\mathbb{F}}\cdot\mathrm{P}|=||\mathrm{P}||
\]
and hence $\mathbf{M}_{\mathbb{F}}=\pm\bar{\mathrm{P}}=\pm\frac{\mathrm{P}}{|\mathrm{P}|}$.
Moreover, note that
\[
2||\mathrm{P}||^{2}=\text{tr}\mathbb{F}_{\bar{\mathbf{n}}}^{2}\leq\text{tr}\mathbb{F}_{\mathbf{M}_{\mathbb{F}}}^{2}\leq2\lambda_{\mathbf{M}_{\mathbb{F}}}^{2}=2||\mathrm{P}||^{2},
\]
implying $\text{tr}\mathbb{F}_{\mathbf{M}_{\mathbb{F}}}^{2}=2\lambda_{\mathbf{M}_{\mathbb{F}}}^{2}=2||\mathrm{P}||^{2}$,
and hence
\[
\frac{|\mathrm{P}\cdot\mathbf{M}_{\mathbb{F}}|+\sqrt{2\text{tr}\mathbb{F}_{\mathbf{M}_{\mathbb{F}}}^{2}-3|\mathrm{P}\cdot\mathbf{M}_{\mathbb{F}}|^{2}}}{2}=\frac{|\mathrm{P}|+\sqrt{4||\mathrm{P}||^{2}-3||\mathrm{P}||^{2}}}{2}=||\mathrm{P}||=|\lambda_{\mathbf{M}_{\mathbb{F}}}|.
\]

It is worth noting that (\ref{eq:=00005Clambda_MF-onPLane}) is also
valid in the case $\mathbf{M}_{\mathbb{F}}=\bar{\mathbf{n}}$, since
then $|\lambda_{\mathbf{M}_{\mathbb{F}}}|=|\lambda_{\bar{\mathbf{n}}}|=||\mathrm{P}||$
and $\text{tr}\mathbb{F}_{\mathbf{M}_{\mathbb{F}}}^{2}=\text{tr}\mathbb{F}_{\bar{\mathbf{n}}}^{2}=2||\mathrm{P}||^{2}$
hold, and thus
\[
\frac{|\mathrm{P}\cdot\mathbf{M}_{\mathbb{F}}|+\sqrt{2\text{tr}\mathbb{F}_{\mathbf{M}_{\mathbb{F}}}^{2}-3|\mathrm{P}\cdot\mathbf{M}_{\mathbb{F}}|^{2}}}{2}=\frac{|\mathrm{P}\cdot\bar{\mathbf{n}}|+\sqrt{4||\mathrm{P}||^{2}-3|\mathrm{P}\cdot\bar{\mathbf{n}}|^{2}}}{2}=||\mathrm{P}||=|\lambda_{\mathbf{M}_{\mathbb{F}}}|.
\]
We summarize these results with the following lemma:
\begin{thm*}
For a magnetic algebra $(\mathbb{R}^{3},\mathbb{F})$
with planarity, the indentity (\ref{eq:=00005Clambda_MF-onPLane})
holds.
\end{thm*}
\vspace{0.6cm}

\subsection*{B. Translational Forces by Magnetic Dipoles}

In this appendix we give a brief introduction to the translation force
generated by a dipole system. For more details in physics, we refer
to \cite{Pollack2001Electromagnetism}.

Consider a single magnetic dipole $\mathfrak{M}$ (with dipole moment
$\mathbf{M}$) and place it at $\mathrm{o}=(x_{\mathrm{o}},y_{\mathrm{o}},z_{\mathrm{o}})\in\mathbb{R}^{3}$.
Denote by $\mathbf{B}$ the magnetic field generated by $\mathfrak{M}$.
At a field point $p=(x,y,z)\in\mathbb{R}^{3}$, 
\[
\mathbf{B}(p)=\frac{\mu_{0}}{4\pi}\big[\frac{3\hat{\mathbf{p}}\hat{\mathbf{p}}^{\intercal}-\mathbf{I}}{||\mathbf{p}||^{3}}\big]\mathbf{M}
\]
with $\mathbf{p}=p-\mathrm{o}$ and $\hat{\mathbf{p}}=\frac{\mathbf{p}}{||\mathbf{p}||}$. 

Given a test magnetic diple $\mathfrak{m}$ (with moment $\mathbf{m}$)
at $p$, the torque exerted on $\mathfrak{m}$ by $\mathbf{B}$ is
\[
\tau=\mathbf{m}\times\mathbf{B},
\]
and the associated orientational potential energy is then
\[
\begin{aligned}\mathrm{U} & = & -\mathbf{m}\cdot\mathbf{B} &  & = & -\frac{\mu_{0}}{4\pi}\mathbf{m}^{\intercal}\big[\frac{3\hat{\mathbf{p}}\hat{\mathbf{p}}^{\intercal}-\mathbf{I}}{||\mathbf{p}||^{3}}\big]\mathbf{M}\\
 &  &  &  & = & \frac{\mu_{0}}{4\pi}\bigg[\frac{\mathbf{M}^{\intercal}\mathbf{m}}{||\mathbf{p}||^{3}}-\frac{3(\mathbf{M}^{\intercal}\mathbf{p})(\mathbf{p}^{\intercal}\mathbf{m})}{||\mathbf{p}||^{5}}\bigg]
\end{aligned}
\]
The translational force experienced by $\mathfrak{m}$ due to the
variation of $\mathbf{B}$ is then
\[
\mathbf{F}=-\nabla\mathrm{U}=-\big(\frac{\partial\mathrm{U}}{\partial x},\frac{\partial\mathrm{U}}{\partial y},\frac{\partial\mathrm{U}}{\partial z}\big)^{\intercal}.
\]
Direct computation shows for each positive integer $k$,
\[
\nabla\frac{1}{||\mathbf{p}||^{k}}=-\frac{k}{||\mathbf{p}||^{k+1}}\hat{\mathbf{p}}.
\]
As a result, 
\[
\begin{aligned}\nabla\mathrm{U} &  & = & -\frac{\mu_{0}}{4\pi}\bigg[\frac{3\mathbf{M}^{\intercal}\mathbf{m}}{||\mathbf{p}||^{4}}\hat{\mathbf{p}}-\frac{15(\mathbf{M}^{\intercal}\mathbf{p})(\mathbf{p}^{\intercal}\mathbf{m})}{||\mathbf{p}||^{6}}\hat{\mathbf{p}}+3\frac{\mathbf{M}\mathbf{p}^{\intercal}+(\mathbf{M}^{\intercal}\mathbf{p})\cdot\mathbf{I}}{||\mathbf{p}||^{5}}\mathbf{m}\bigg]\\
 &  & = & -\frac{\mu_{0}}{4\pi}\bigg[\frac{3\mathbf{M}^{\intercal}\mathbf{m}}{||\mathbf{p}||^{4}}\hat{\mathbf{p}}-\frac{15(\mathbf{M}^{\intercal}\hat{\mathbf{p}})(\hat{\mathbf{p}}^{\intercal}\mathbf{m})}{||\mathbf{p}||^{4}}\hat{\mathbf{p}}+3\frac{\mathbf{M}\hat{\mathbf{p}}^{\intercal}+(\mathbf{M}^{\intercal}\mathbf{\hat{p}})\cdot\mathbf{I}}{||\mathbf{p}||^{4}}\mathbf{m}\bigg]\\
 &  & = & -\frac{3\mu_{0}}{4\pi}\bigg[\frac{\hat{\mathbf{p}}\mathbf{M}^{\intercal}+\mathbf{M}\hat{\mathbf{p}}^{\intercal}+(\mathbf{M}^{\intercal}\mathbf{\hat{p}})\cdot\mathbf{I}}{||\mathbf{p}||^{4}}\mathbf{m}-\frac{5(\mathbf{M}^{\intercal}\hat{\mathbf{p}})}{||\mathbf{p}||^{4}}\hat{\mathbf{p}}\hat{\mathbf{p}}^{\intercal}\bigg]\mathbf{m}.
\end{aligned}
\]

Now consider a system of $n$ magnets $\mathfrak{M}_{i}$, each of
which has (the same) dipole moment $\mathbf{M}$ and is placed at
$\mathrm{o}_{i}=(x_{\mathrm{o}_{i}},y_{\mathrm{o}_{i}},z_{\mathrm{o}_{i}})$,
respectively. Then, with $\mathbf{p}_{i}=p-\mathrm{o}_{i}$ and $\hat{\mathbf{p}}_{i}=\frac{\mathbf{p}_{i}}{||\mathbf{p}_{i}||}$,
the gradient of the potential energy by $\mathfrak{M}_{i}$ is
\[
\nabla\mathrm{U}_{i}=-\frac{3\mu_{0}}{4\pi}\bigg[\frac{\hat{\mathbf{p}}_{i}\mathbf{M}^{\intercal}+\mathbf{M}\hat{\mathbf{p}}_{i}^{\intercal}+(\mathbf{M}^{\intercal}\mathbf{\hat{p}}_{i})\cdot\mathbf{I}}{||\mathbf{p}_{i}||^{4}}\mathbf{m}-\frac{5(\mathbf{M}^{\intercal}\hat{\mathbf{p}}_{i})}{||\mathbf{p}_{i}||^{4}}\hat{\mathbf{p}}_{i}\hat{\mathbf{p}}_{i}^{\intercal}\bigg]\mathbf{m},
\]
and hence the total translational force on $\mathfrak{m}$ is
\[
\begin{aligned}\mathbf{F} &= -\sum_{i}\nabla\mathrm{U}_{i}=\frac{3\mu_{0}}{4\pi}\mathbb{F}_{\mathbf{M}}\mathbf{m}\\&
 =\frac{3\mu_{0}}{4\pi}\bigg[\underset{i}{\sum}\frac{\hat{\mathbf{p}}_{i}}{||\mathbf{p}_{i}||^{4}}\mathbf{M}^{\intercal}+\mathbf{M}\underset{i}{\sum}\frac{\hat{\mathbf{p}}_{i}^{\intercal}}{||\mathbf{p}_{i}||^{4}}+(\mathbf{M}^{\intercal}\underset{i}{\sum}\frac{\hat{\mathbf{p}}_{i}}{||\mathbf{p}_{i}||^{4}})\cdot[\mathbf{I}]-\underset{i}{\sum}\frac{5(\mathbf{M}^{\intercal}\hat{\mathbf{p}}_{i})}{||\mathbf{p}_{i}||^{4}}\hat{\mathbf{p}}_{i}\hat{\mathbf{p}}_{i}^{\intercal}\bigg]\mathbf{m},
\end{aligned}
\]
which verifies (\ref{eq:F(M)}) with $\mathrm{P}=\underset{i}{\sum}\frac{\hat{\mathbf{p}}_{i}}{||\mathbf{p}_{i}||^{4}}$.

\end{document}